\documentclass[11pt, reqno]{amsart}
\usepackage{textcmds}  
\usepackage{amsmath, amssymb, amsfonts, amstext, verbatim, amsthm, mathrsfs}
\usepackage[mathcal]{eucal}
\usepackage{microtype}
\usepackage{graphicx} 
\usepackage[all]{xy}
\usepackage[modulo]{lineno}
\usepackage[usenames]{color}
\usepackage{aliascnt} 
\usepackage[colorlinks=true,linkcolor=blue,citecolor=blue,urlcolor=blue,citebordercolor={0 0 1},urlbordercolor={0 0 1},linkbordercolor={0 0 1}]{hyperref} 
\usepackage[shortalphabetic]{amsrefs} 
\usepackage{enumitem}
\usepackage{xspace}
\usepackage[normalem]{ulem}
\usepackage{stmaryrd}

\usepackage{fullpage}

\theoremstyle{plain}

\newcommand{\refnewtheoremn}[4]{%
\newaliascnt{#1}{#2}
\newtheorem{#1}[#1]{#3}
\aliascntresetthe{#1}
\expandafter\providecommand\csname #1autorefname\endcsname{#4}}

\newcommand{\refnewtheorem}[3]{\refnewtheoremn{#1}{#2}{#3}{#3}}

\newtheorem{thm}{Theorem}[section]

\refnewtheorem{lem}{thm}{Lemma}
\refnewtheorem{claim}{thm}{Claim}
\refnewtheorem{cor}{thm}{Corollary}
\refnewtheorem{conj}{thm}{Conjecture}
\refnewtheorem{prop}{thm}{Proposition}
\refnewtheorem{quest}{thm}{Question}
\refnewtheorem{amplif}{thm}{Amplification}
\refnewtheorem{scholium}{thm}{Scholium}

\refnewtheorem{tablethm}{thm}{Table}
\refnewtheorem{assumption}{thm}{Assumption}
\refnewtheorem{conclusion}{thm}{Conclusion}
\refnewtheorem{problem}{thm}{Problem}

\newtheorem{thmx}{Theorem}

\theoremstyle{definition}
\refnewtheorem{rem}{thm}{Remark}
\refnewtheorem{defn}{thm}{Definition}
\refnewtheorem{construction}{thm}{Construction}
\refnewtheorem{notn}{thm}{Notation}
\refnewtheorem{ex}{thm}{Example}
\newtheorem*{ex*}{Example}
\refnewtheorem{fact}{thm}{Fact}
\newtheorem*{rem*}{Remark}

{\end{minipage} \end{center}}

\newcommand{\bA}{\mathbf{A}}
\newcommand{\cA}{\mathcal{A}}

\newcommand{\cB}{\mathcal{B}}

\newcommand{\bC}{\mathbf{C}}
\newcommand{\cC}{\mathcal{C}}

\newcommand{\cD}{\mathcal{D}}

\newcommand{\cE}{\mathcal{E}}

\newcommand{\cF}{\mathcal{F}}

\newcommand{\bG}{\mathbf{G}}
\newcommand{\cG}{\mathcal{G}}

\newcommand{\bL}{\mathbf{L}}

\newcommand{\sM}{\mathscr{M}}
\newcommand{\bN}{\mathbf{N}}

\newcommand{\cO}{\mathcal{O}}

\newcommand{\bP}{\mathbf{P}}
\newcommand{\cP}{\mathcal{P}}

\newcommand{\bQ}{\mathbf{Q}}

\newcommand{\cS}{\mathcal{S}}

\newcommand{\cU}{\mathcal{U}}

\newcommand{\cV}{\mathcal{V}}

\newcommand{\cW}{\mathcal{W}}

\newcommand{\cX}{\mathcal{X}}

\newcommand{\cY}{\mathcal{Y}}

\newcommand{\bZ}{\mathbf{Z}}
\newcommand{\cZ}{\mathcal{Z}}

\newcommand{\fg}{\mathfrak{g}}

\newcommand{\fm}{\mathfrak{m}}

\newcommand{\fz}{\mathfrak{z}}

\newcommand{\id}{\mathop{{\rm id}}\nolimits}

\renewcommand{\bC}{{\mathbb C}}
\renewcommand{\bZ}{{\mathbb Z}}
\renewcommand{\bQ}{{\mathbb Q}}

\newcommand{\op}[1]{\!\!\mathop{\rm ~#1}\nolimits}

\newcommand{\tinyop}[1]{\!\!\mathop{\mbox{\rm \tiny ~#1}}\nolimits}

\renewcommand{\bA}{\mathbb{A}}
\renewcommand{\bP}{\mathbb{P}}

\newcommand{\Gm}{\mathbb{G}_m}

\newcommand{\dual}{\vee}
\newcommand{\sod}[1]{\left\langle #1 \right\rangle}

\newcommand{\X}{\cX}
\newcommand{\Y}{\cY}
\newcommand{\Z}{\cZ}

\newcommand{\U}{\cU}

\newcommand{\Crit}{\op{Crit}_0}
\newcommand{\Perf}{\op{Perf}}
\newcommand{\APerf}{\operatorname{D}^{-}\operatorname{Coh}}
\newcommand{\A}{\cA}
\newcommand{\Aff}{\op{Aff}}
\newcommand{\RHom}{\op{RHom}}
\newcommand{\Hom}{\op{Hom}}
\newcommand{\Coh}{\op{Coh}}

\newcommand{\Rep}{\op{Rep}}

\newcommand{\Mod}{\!\operatorname{-Mod}}

\newcommand{\Map}{\op{Map}}

\newcommand{\Sp}{\op{Sp}}

\newcommand{\Space}{\mathbf{Spaces}}

\newcommand{\Madd}{\sM}
\newcommand{\DSing}{\operatorname{D^b Sing}}

\newcommand{\MF}{\operatorname{MF}}
\newcommand{\PreMF}{\op{PreMF}}

\newcommand{\pt}{\ast}

\newcommand{\DCoh}{\operatorname{D^b Coh}}
\newcommand{\hoch}{C_\bullet}
\newcommand{\Tate}{{\rm Tate}}
\newcommand{\Sone}{{S^1}}
\newcommand{\QC}{\op{QC}}
\newcommand{\IC}{\op{IndCoh}}

\newcommand{\hocolim}{\op{hocolim}}

\newcommand{\RGamma}{\operatorname{R\Gamma}}
\newcommand{\iHom}{\mathcal{H}om}
\newcommand{\HKR}{\operatorname{HKR}}
\newcommand{\Exp}{\operatorname{exp}}
\newcommand{\Sym}{\operatorname{Sym}}
\newcommand{\Fun}{\operatorname{Fun}}

\newcommand{\ps}[1]{[\![#1]\!]}
\newcommand{\ls}[1]{(\!(#1)\!)}
\newcommand{\itimes}{\dot{\otimes}}

\begin{document}

\title{Equivariant Hodge theory and noncommutative geometry}
\author{Daniel Halpern-Leistner}
\author{Daniel Pomerleano}

\begin{abstract}
We develop a version of Hodge theory for a large class of smooth formally proper quotient stacks $X/G$ analogous to Hodge theory for smooth projective schemes. We show that the noncommutative Hodge-de Rham sequence for the category of equivariant coherent sheaves degenerates. This spectral sequence converges to the periodic cyclic homology, which we canonically identify with the topological equivariant $K$-theory of $X$ with respect to a maximal compact subgroup of $G$, equipping the latter with a canonical pure Hodge structure. We also establish Hodge-de Rham degeneration for categories of matrix factorizations for a large class of equivariant Landau-Ginzburg models.
\end{abstract}

\maketitle
\tableofcontents

\addtocontents{toc}{\protect\setcounter{tocdepth}{-1}}

If $X$ is a smooth projective variety over $\bC$, then the cohomology groups $H^n(X;\bC)$ can be equipped with a pure Hodge structure of weight $n$. The theory of Hodge structures then allows one to ``linearize" many important problems in algebraic geometry. Our goal is to develop such a linearization for the equivariant algebraic geometry of a locally closed algebraic submanifold $X \subset \bP^n$ which is equivariant with respect to an action of a compact Lie group, $M$. Note that the complexification $G$ of $M$, a reductive algebraic group, acts on $X$ as well, and it is natural to ask for a Hodge theory associated intrinsically to the algebraic stack $\X:= X/G$.

One such linearization follows from the results of \cite{deligne1974theorie}, which establishes a canonical mixed Hodge structure on the cohomology of any smooth simplicial scheme and in particular on the equivariant cohomology, $H_G^\ast(X)$, which is the cohomology of the simplicial nerve of the action groupoid of $G$ on $X$. Building on these ideas, one can even associate a motive to the stack $X/G$ as a colimit of motives of schemes as in \cite{morel1999A1homotopy}*{Section 4.2}.

The present paper lays the groundwork for an alternative approach to equivariant Hodge theory based on equivariant topological $K$-theory of the underlying analytic variety $K^\ast_M(X^{an})$. The equivariant $K$-theory is a module over the representation ring $\Rep(G)$, and the Atiyah-Segal completion theorem canonically identifies its completion at the augmentation ideal of $\Rep(G)$ with the $\bZ/2\bZ$-graded equivariant cohomology, $$K^{\rm even/odd}_M(X^{an})^\wedge \cong H_G^{\rm even/odd}(X^{an}).$$ If one completes  $K_G^\ast(X^{an})$ with respect to the evaluation ideal at other $g\in M$, then one recovers the $\bZ/2\bZ$-graded equivariant cohomology of the fixed locus of $g$ with respect to the centralizer of $g$ \cite{MR2365650}*{Theorem 3.9}. Thus $K_G^\ast(X^{an})$ encodes all of these equivariant cohomology groups, as well as the data of how to ``spread'' them out into a single finitely generated $\Rep(G)$-module. In this sense $K^\ast_M(X^{an})$ is a much richer invariant than $H_G^\ast(X)$.


The main challenge is that unlike cohomology, equivariant $K$-theory is not simply the $K$-theory of the simplicial scheme arising from the action of $G$ on $X$, so Deligne's approach to equivariant Hodge theory does not generalize to $K$-theory. Instead, our Hodge structures originate in \emph{noncommutative} algebraic geometry, which views dg-categories as ``noncommutative spaces." We ultimately show that in many cases one can use the dg-enhanced derived category of $G$-equivariant perfect complexes of coherent sheaves on $X$, $\Perf(X/G)$, to construct a pure Hodge structure on $K^\ast_M(X^{an})$.


\subsection*{Noncommutative Hodge-de-Rham degeneration}


If $\cA$ is a dg-category over a field $k$, the Hochschild chain complex, $C_{\bullet}(\cA)$ plays the role of the Hodge cohomology in noncommutative algebraic geometry. The periodic cyclic complex $C_\bullet^{per}(\cA)$, which is a dg-module over $k(\!(u)\!)$ where $u$ has homological degree $-2$, behaves like noncommutative de Rham cohomology. There is a canonical Hodge filtration of the complex $C_\bullet^{per}(\cA)$ whose associated graded is $C_{\bullet}(\cA) \otimes k(\!(u)\!)$, which leads to a noncommutative Hodge-de Rham spectral sequence converging to $H_\ast(C_\bullet^{per}(\cA))$ whose first page is $H_\ast(C_\bullet(\cA)\otimes k(\!(u)\!))$. The Hodge filtration in our theory will be the filtration on $H_\ast(C_\bullet^{per}(\cA))$ arising from the degeneration of this spectral sequence.

Geometrically, our main tool for establishing degeneration for the category $\Perf(X/G)$ will be certain stratifications of $X$ into $G$-stable locally closed subvarieties. It is a classical observation in geometric invariant theory \cites{hesselink, ness} that projective varieties with a linearized $G$-action inherit a canonical stratification. Kirwan used this stratification to deduce many beautiful results concerning the equivariant topology of projective varieties with $G$-action, including a computation of the Betti numbers and Hodge numbers for a variety obtained as a GIT quotient \cite{MR766741}.

In this paper, we will consider a certain abstraction of this canonical stratification. These are the (semi-)complete Kirwan-Ness(KN)-stratifications of a $G$-variety, introduced in \cite{teleman2000quantization} and recalled in \autoref{def:semicomplete_KN}. The chief benefit of this more abstract definition is that it applies in many cases when the ambient variety $X$ is merely quasi-projective. The main classes of examples of $G$-varieties which admit semi-complete KN-stratifications to keep in mind are:
\begin{enumerate}[label=(\roman*)]
\item any $G$-variety $X$ which is projective over an affine $G$-variety; and
\item any $G$-variety $X$ such that $X/G$ admits a good quotient that is projective-over-affine.\footnote{Recall that $X/G$ admits a good quotient if there is an algebraic space $Y$ and a $G$-invariant map $\pi : X \to Y$ such that $\pi_\ast : \op{QCoh}(X/G) \to \op{QCoh}(Y)$ is exact and $(\pi_\ast \cO_X)^G \simeq \cO_Y$.}
\end{enumerate}
In the case (i), the KN stratification is complete if and only if $\dim \Gamma(X,\cO_X)^G < \infty$. In the case (ii) the stratification is trivial, and it is complete if and only if the good quotient of $X/G$ is projective.

\begin{thmx}[\autoref{cor:degeneration_complete}] \label{thm:intro_degeneration}
If $G$ is a reductive group and $X$ is a smooth $G$-quasi-projective variety which admits a complete KN-stratification, then the noncommutative Hodge-de Rham sequence for $\Perf(X/G)$ degenerates on the first page.
\end{thmx}

\begin{rem*}
This builds on \cite{teleman2000quantization}, which shows that a version of the Hodge-de Rham spectral sequence for $H_G^\ast(X)$ degenerates for such $G$-schemes and that the (a priori mixed) Hodge structure on $H_G^\ast(X)$ is pure in this case. Note that in these examples, the scheme $X$ is not proper, and neither is the quotient stack $X/G$, so degeneration of the Hodge-de Rham sequence is somewhat unexpected. Likewise from the noncommutative perspective, D. Kaledin's recent resolution \cite{kaledin2017spectral} of a well-known conjecture of M. Kontsevich and Y. Soibelman \cite{kontsevich2009notes} shows that the noncommutative Hodge-de Rham sequence degenerates for dg-categories which are smooth and proper. However, the categories $\Perf(X/G)$ are typically not smooth even when $X$ is smooth, and they are typically not generated by a single compact generator, so degeneration is again somewhat unexpected.
\end{rem*}

The key observation in establishing the degeneration property for $\Perf(X/G)$ is that the formation of the Hochschild complex takes semiorthogonal decompositions of dg-categories to direct sums, and its formation commutes with filtered colimits. Thus if $\cA$ is a retract of a dg-category which can be built from the derived category of smooth and proper DM stacks via an infinite semiorthogonal decomposition, then the degeneration property holds for $\cA$.

\begin{ex*}
One simple example is the quotient stack $\bA^n / \Gm$, where $\Gm$ acts with positive weights. Then the objects $\cO_{\bA^n}\{w\} \in \Perf(\bA^n / \Gm)$, which denote the twist of the structure sheaf by a character of $\Gm$, form an infinite full exceptional collection. Therefore the Hochschild complex of $\Perf(\bA^n/\Gm)$ is quasi-isomorphic to a countable direct sum of copies of $C_\bullet(\Perf(\op{Spec}(k)))$, and the degeneration property follows.
\end{ex*}

We can formulate this most cleanly in terms of G. Tabuada's universal additive invariant $\cU_k : \op{dgCat}_k \to \Madd_k$ \cites{tabuada2008higher, blumberg2013universal}. Here $\Madd_k$ is the $\infty$-category which is the localization of the $\infty$-category of small dg-categories which formally splits all semiorthogonal decompositions into direct sums, and $\cU_k$ is the localization map. The following is the main technical result of the paper, and we believe it is of independent interest.

\begin{thmx}[See \autoref{thm:motivic_main}] \label{thm:intro_C}
Let $X/G$ be a smooth quotient stack over a field $k$ of characteristic $0$ that admits a complete KN stratification. Then there is a smooth projective variety $Y$ such that $\cU_k(X/G)$ is a direct summand of $\cU_k(D^b(Y))^{\oplus \bN}$ in $\Madd_k$.
\end{thmx}

\subsection*{Connections with (classical) equivariant topology and purity}


If $G$ is the complexification of a compact Lie group $M$ as above, we show that one can recover the equivariant topological $K$-theory of the underlying complex analytic space $K_M(X^{an})$, as defined in \cites{atiyah1969equivariant,segal1968equivariant}, from the dg-category $\Perf(X/G)$.

The first ingredient is the recent construction by A. Blanc of a topological $K$-theory spectrum $K^{top}(\cA)$ for any dg-category $\cA$ over $\bC$ \cite{blanc2012topological}. Blanc constructs a Chern character natural transformation $\op{ch} : K^{top}(\cA) \to HP(\cA)$, shows that $\op{ch} \otimes \bC$ is an equivalence for $\Perf$ of a finite type $\bC$-scheme, and conjectures this property for any smooth and proper dg-category $\cA$. We show that $\op{ch} \otimes \bC$ is an isomorphism for all categories of the form $\Perf(\Y)$, where $\Y$ is a smooth DM stack or a smooth quotient stack admitting a semi-complete KN stratification. In fact, we expect that this ``lattice conjecture" should hold for a much larger class of dg-categories, such as the categories $D^b(\X)$ for any finite type $\bC$-stack and $\operatorname{Perf}(X/G)$ for any quotient stack. Following some ideas of Thomason in \cite{thomason1988equivariant}, we next construct a natural ``topologization'' map $\rho_{G,X} : K^{top}(\Perf(X/G)) \to K_M(X^{an})$ for any smooth $G$-quasiprojective scheme $X$ and show:

\begin{thmx}[See \autoref{thm:topological_comparison} and \autoref{thm:lattice_conjecture}] \label{thm:intro_A}
For any smooth quasi-projective $G$-scheme $X$ which admits a semi-complete KN stratification, the topologization map and the Chern character provide equivalences \footnote{We will see that these homology level equivalences are induced by suitable chain maps.} 
$$\xymatrix{ K^\ast_M(X^{an}) \otimes \bC & \pi_\ast K^{top}(\Perf(X/G)) \otimes \bC \ar[r]^{\op{ch}} \ar[l]_-{\rho_{G,X}} & H_*C_\bullet^{per}(\Perf(X/G))}$$ 
\end{thmx}

\begin{rem*}
In fact, \autoref{thm:topological_comparison} shows a bit more. $\rho_{G,X}$ is an equivalence for any smooth $G$-quasi-projective scheme. For an arbitrary $G$-quasi-projective scheme $X$ we construct an equivalence of spectra $\rho_{G,X} : K^{top}(\DCoh(X/G)) \to K_M^{c,\dual}(X^{an})$, where the latter denotes the $M$-equivariant Spanier-Whitehead dual of the spectrum $K_M(X^{an})$, sometimes referred to as the equivariant Borel-Moore $K$-homology of $X$. $\rho_{G,X}$ is compatible up to homotopy with natural pullback and pushforward maps (to be explained below). This result is of independent interest, and it allows one to ``decategorify" theorems regarding equivariant derived categories in a precise way.
\end{rem*}

Note that the groups $K_M^n(X^{an})$ are modules over $\Rep(M)$, the representation ring of $M$. We say that a $\Rep(M)$-linear Hodge structure of weight $n$ is a finite $\Rep(M)$-module $E$ along with a finite filtration of the finite $\Rep(M)_\bC$-module $E \otimes \bC$ inducing a Hodge structure of weight $n$ on the underlying abelian group $E$. Using the previous identification $K_M^n(X^{an}) \otimes \bC \simeq H_{-n}C_\bullet^{per}(\Perf(X/G))$, we will show

\begin{thmx}[See \autoref{thm:Hodge_main}] \label{thm:intro_B}
For any smooth $M$-quasiprojective scheme admitting a complete KN stratification, the noncommutative Hodge-de Rham sequence for $K_M^n(X^{an}) \otimes \bC$ degenerates on the first page, equipping $K_M^n(X^{an})$ with a pure $\Rep(M)$-linear Hodge structure of weight $n$, functorial in $X$. There is a canonical isomorphism
$$\op{gr}^p K_M^n(X^{an}) \simeq H^{n-2p}(R\Gamma(I^{der}_{\X},\cO_{I^{der}_\X})).$$
\end{thmx}
In this theorem, $I_\X^{der}$ denotes the derived inertia stack, sometimes referred to as the ``derived loop stack.'' As we will see in \autoref{lem:derivedint} below, we can express this more concretely as
$$R\Gamma(I^{der}_\X,\cO_{I^{der}_\X}) \simeq R\Gamma(G\times X \times X, \cO_\Gamma \otimes^L \cO_{\bar{\Delta}})^G,$$
where $G$ acts on $G \times X \times X$ by $g \cdot (h,x,y) = (g h g^{-1},gx,gy)$ and the two $G$-equivariant closed subschemes of $G \times X \times X$ are defined as $\Gamma = \{(g,x,gx)\}$ and $\bar{\Delta} = \{(g,x,x)\}$ respectively.

\begin{ex*}
Along the way, show that the lattice conjecture holds for an arbitrary smooth DM stack, and explicitly compute the Hochschild invariants of $\Perf(\X)$. For a smooth and proper DM stack, we construct an isomorphism of Hodge structures
$$\pi_{n} K^{top}(\Perf(\X)) \otimes \bQ \simeq \bigoplus_k H^{2k-n}_{Betti}(I^{cl}_\X;\bQ\langle k\rangle).$$
\end{ex*}

It should be noted that the motivic decompositions of Theorem \autoref{thm:intro_C} play a key role in the proof of Theorem \autoref{thm:intro_B}, but these decompositions do not respect the $\operatorname{Rep}(M)$-linear nature of the Hodge structure on $K_M^n(X^{an})$.

In \autoref{sect:explicit} we spend some time discussing more explicit models for the Hochschild homology and periodic cyclic homology for quotient stacks. For example, we show that when $X$ is smooth and affine, there is an explicit bar-type complex computing the Hochschild homology of $\Perf(X/G)$. As an application of Theorem \ref{thm:intro_B}, we prove an HKR type theorem for the completion of this bar complex at various points of $\op{Spec}(\Rep(G))$ when $X/G$ is formally proper. A corollary of this theorem is a description of the completed Hochschild homology modules equipped with the Connes operator in terms of differential forms equipped with the de Rham differential.

\subsection*{Extensions to categories of singularities}

Another major source of Hodge structures in algebraic geometry comes from singularity theory. For instance in \cite{MR723468}, Kyoji Saito constructs analogues of Hodge theoretic structures on the universal unfolding of an isolated singularity. More precisely, he describes analogues of the Gauss-Manin connection and period mappings as well as canonical coordinates on the base space of the universal unfolding. Motivated in part by Saito's work, Katzarkov, Kontsevich, and Pantev \cites{KKP,KKP2} have proposed a vast generalization of Hodge theory which they call \emph{noncommutative} (nc) Hodge structures. As the name suggests, they envision that nc Hodge structures should arise naturally from smooth and proper dg-categories(``nc spaces").

Let $(X,W)$ be a \emph{Landau-Ginzburg} (LG) model, that is a pair $(X,W)$ consisting of a smooth quasi-projective variety $X$ and a regular function $W: X \to \mathbb{A}^1.$ To any LG model, one may associate the category of \emph{matrix factorizations} $\MF(X,W)$, which is a 2-periodic (meaning $k\ls{\beta}$-linear where $\beta$ has homological degree $-2$) dg-category. Applying the theory of nc Hodge structures to these nc spaces is expected to yield a vast generalization of Saito's theory to pairs $(X,W)$ with proper critical locus, $\operatorname{Crit}(W)$.\footnote{Noncommutative Hodge structures are also expected to exist in other contexts, notably on the quantum cohomology of a compact symplectic manifold.} 

For any LG model $(X,W)$, there is a ``$dW$-twisted" Hodge de-Rham spectral sequence which relates the hypercohomologies of the complexes $(\Omega_X^{\bullet}, dW \wedge)$ and $(\Omega_X {^\bullet}, d+ dW \wedge)$(see \cite{KKP}*{Section 3.2} for details). Similarly to the classical case, this spectral sequence is known to degenerate when $W: X \to \mathbb{A}^1$ is proper by work of Ogus-Vologodsky \cite{ogus2007nonabelian}. This degeneration result plays a central role in the noncommutative Hodge theory of LG pairs---for example, a version of this result has been used to establish the smoothness of versal deformation spaces of (compactified) LG models (generalizing the universal unfolding space of a singularity)\cite{KKP2}.  Efimov and Preygel \cites{Preygel, efimov_cyclic} have independently identified the $dW$-twisted Hodge de-Rham spectral sequence with the $k\ls{\beta}$-linear noncommutative Hodge-de Rham spectral sequence for the category $\MF(X,W)$(for closely related results, see also \cites{dyckerhoff2011compact,esegalclosed,caldararu2013curved, lin2013global, shklyarov}). It follows that the result of Ogus-Vologodsky can be recast as establishing the degeneration of this spectral sequence of noncommutative origin.

We prove the following generalization of this degeneration result, which suggests that in the equivariant context nc Hodge theory should extend to certain dg-categories which are not smooth. We let $\Crit(W)$ denote the critical locus with critical value $0$.

\begin{thmx}(\autoref{prop:MFdeg}) \label{thm:intro_mf_degen}
If $X$ is a smooth $G$-quasi-projective scheme which admits a semi-complete KN stratification, and $W : X \to \bA^1$ is a $G$-invariant function such that $\Perf(\Crit(W)/G)$ is a proper dg-category, then the $k\ls{\beta}$-linear noncommutative Hodge-de Rham sequence for $\MF(X/G,W)$ degenerates on the first page.
\end{thmx}

Note that by \autoref{lem:properness}, the condition on $\Crit(W)$ in the theorem is equivalent to the induced KN-stratification on $\Crit(W)$ being complete. In the case where $X$ is projective over an affine $G$-variety, the condition is equivalent to $\dim \Gamma(\Crit(W),\cO_{\Crit(W)})^G<\infty$, and if $X/G$ admits a good quotient, the condition is equivalent to the condition that $\Crit(W)/G$ admits a projective good quotient. 

We prove Theorem \autoref{thm:intro_mf_degen} by establishing an analog of Theorem \autoref{thm:intro_C} for the $k\ls{\beta}$-linear category $MF(X/G,W)$ in \autoref{thm:motivic_main_lg}. The proof is somewhat more subtle than the case of $\Perf(X/G)$, and its formulation is a little more complicated, because at the time of this writing we are not aware of a construction of $k\ls{\beta}$-linear additive noncommutative motives. Along the way, we also establish the degeneration property for $\MF(\X,W)$ in the case that $\X$ is a smooth quasi-projective DM stack and $\Crit(W)$ is proper (see \autoref{sect:singularities}).\footnote{As mentioned above, after the first draft of this paper was circulated, D. Kaledin proved the degeneration conjecture. To the authors' knowledge, however, the version of the degeneration conjecture for $k(\!(\beta)\!)$-linear categories, which is the one which applies to categories of the form $\MF(X,W)$, does not appear in the literature.} 



\subsection*{Further questions}

\subsubsection*{The notion of properness in equivariant geometry}

Our result on noncommutative Hodge-de Rham degeneration adds to the list of ways in which certain equivariant geometries behave as if they are proper despite not being proper in the sense of algebraic stacks. The intrinsic characterization of which smooth algebraic stacks behave as if they are proper from the perspective of Hodge theory, such as quotient stacks with a complete KN stratification, and which do not, such as $B \bG_a$ or $BU$ for a unipotent group $U$ (see \autoref{ex:bga}), is still somewhat fuzzy.

The paper \cite{halpern2014mapping} studies these properness phenomena systematically by introducing the class of \emph{formally proper} stacks, with the primary application being the algebraicity of the mapping stack out of a formally proper stack. The examples and counterexamples above for stacks exhibiting noncommutative Hodge-de Rham degeneration are also important examples and counterexamples for stacks which are formally proper in the sense of \cite{halpern2014mapping}. This raises the natural question.
\begin{quest}
Do there exist examples of perfect, smooth, and formally proper $k$ stacks $\cX$ for which the Hodge-de Rham sequence associated to $\Perf(\cX)$ does not degenerate?
\end{quest}


\subsubsection*{Hodge structures on equivariant K-theory}

We believe that our main theorem for Hodge structures on $K_M^n(X^{an})$ raises many questions for further inquiry into the role of Hodge theory in equivariant algebraic geometry. For example, it is plausible that the results above could be extended to construct mixed Hodge structures on some version of $K$-theory for arbitrary finite type stacks. In a different direction, one of the central notions in Hodge theory is that of a variation of Hodge structure. For simplicity, let $S$ be an affine scheme and suppose further that $\pi: X/G \to S$ is a smooth equivariant family over $S$ such that all of the fibers $X_s/G$ admit complete KN stratifications. Most of the techniques that we have developed work in families, which allows one to establish the existence of suitable Hodge filtrations on the quasi-coherent sheaf $H_*C_{S}^{per}(\Perf(X/G))$. We therefore believe it is quite likely that one can develop a theory of equivariant period maps. Finally, Theorem \autoref{thm:intro_mf_degen} suggests that it may be possible to develop a version of noncommutative Hodge theory which applies in the equivariant context.  

\subsubsection*{Noncommutative equivariant geometry}

Although we make use of noncommutative algebraic geometry, all of the differential graded categories in this paper are of commutative origin. It is interesting to try to formulate in noncommutative terms a criterion for the Hodge-de Rham spectral sequence to degenerate. Theorem \ref{thm:intro_B} suggests the following concrete question: Let $\cA$ be a proper dg-category which is a module over $\operatorname{Perf}(BG)$. Suppose that $\cA \otimes_{\operatorname{Perf}(BG)} k \cong \operatorname{Perf}(\mathcal{R})$, where $\mathcal{R}$ is a dg- algebra which is homotopically finitely presented, homologically bounded and such that $H_*(\mathcal{R})$ is a finitely generated module over  $HH^0(\mathcal{R})$. 

\begin{quest} 
Does the Hodge-de Rham spectral sequence always degenerate for such $\cA$?
\end{quest} 

\subsubsection*{Context} \label{S:preliminaries}
Throughout this work, unless explicitly stated otherwise, we work over a fixed subfield $k \subset \bC$. All of our functors are understood to be derived, so we write $i_\ast$ for $Ri_\ast$, $i^\ast$ for $Li^\ast$, $\Hom$ for $R\Hom$, etc. We will work with stacks over the \'{e}tale site of $k$-schemes. By convention, unless otherwise indicated the term \emph{quotient stack} will denote a quotient of a quasi-projective $k$-scheme by a linearizable action of an algebraic $k$-group $G$,\footnote{This is sometimes referred to as a $G$-quasi-projective scheme.} and we denote it $X/G$.

Our stacks will be classical whenever we are studying the derived category of coherent sheaves $\DCoh(X/G)$ and its relatives ($\QC(\X)$,$\Perf(\X)$, etc.), but when we discuss categories of matrix factorizations $\MF(\X,W)$ and its relatives ($\IC(\X)$, $\PreMF(\X,W)$, $\PreMF^\infty(\X,W)$, etc.), it will be convenient to work with derived stacks.


We will work with $k$-linear dg-categories. For some of the more abstract arguments involving homotopy limits and colimits and symmetric monoidal structures, it will be more convenient to replace them with equivalent stable (i.e. pre-triangulated) dg-categories in the Morita model structure on dg-categories, then to regard them as $k$-linear stable $\infty$-categories via the equivalence of \cite{cohn2013differential}. We permit ourselves a bit of fluidity on this point, in that we refer both to the literature on dg-categories and stable $\infty$-categories as needed for constructions which evidently make sense in either context.

\subsection*{Acknowledgements}
We would like to thank D. Ben-Zvi, D. Kaledin, F. Morel and A. Preygel for useful conversations. Our thesis advisor C. Teleman explained to us many of the ideas in equivariant K-theory, geometric invariant theory, and noncommutative geometry which this paper builds upon. We are grateful to B. Toen for explaining his computation of the periodic cyclic homology of a Deligne-Mumford stack, and to A. Blanc for useful discussions and for his careful reading of the first draft of this manuscript. Finally, we would like to thank the anonymous referee for several suggestions which improved this article. The first named author was supported by an NSF fellowship, Columbia University, and the Institute for Advanced Study.  The second named author was supported by Kavli IPMU, Imperial College, and an EPSRC fellowship. 

\addtocontents{toc}{\protect\setcounter{tocdepth}{2}}


\section{The noncommutative motive of a quotient stack}

In this section, we show that the noncommutative Hodge-de Rham spectral sequence degenerates for $\Perf(\X)$ for a large class of smooth quotient stacks subject to a properness condition. Our method for establishing the degeneration property will be to systematically realize the derived category of a smooth quotient stack as being ``glued together'' from (typically infinitely many) copies of the derived category of smooth Deligne-Mumford stacks. This method will be used several times throughout this paper, so we formulate our main result in a way that can be applied directly in different contexts.

We work with the category $\Madd_k$ of $k$-linear additive motives in the sense of \cite{tabuada2008higher} (see also \cite{blumberg2013universal} for a construction using the framework of $\infty$-categories). This is the $\infty$-category obtained as the left Bousfield localization of the $\infty$-category of small $k$-linear dg-categories localized at the class of morphisms $\cC \to \cA \oplus \cB$ coming from split exact sequences of small dg-categories $\cA \to \cC \to \cB$. In other words, objects of $\Madd_k$ are dg-categories $[\cC]$, where we have formally adjoined the relation $[\cC] = [\cA] \oplus [\cB]$ whenever we have a semiorthogonal decomposition $\cC = \sod{\cA,\cB}$. We denote the localization functor $\cU_k : \op{dgCat}_k \to \Madd_k$.

\subsection{Recollections on KN-stratifications}

Our primary geometric tool will be a ``KN stratification'' of a quotient stack, as defined in \cite{teleman2000quantization}*{(1.1)} or \cite{halpern2015derived}*{Definition 2.2}. This is a decomposition of $X$ as a union of $G$-equivariant, smooth, locally closed subschemes
\begin{equation} \label{eqn:KN_stratification}
X / G = X^{ss} / G \cup \bigcup_i {S_i / G}.
\end{equation}
For instance, when $X$ is projective-over-affine and $G$ is reductive, a KN-stratification of $X/G$ is induced by a choice of $G$-linearized ample line bundle $L$ and a Weyl-invariant inner product on the cocharacter lattice of $G$. Throughout our discussion, we will assume that we have fixed a choice of inner product on the cocharacter lattice of $G$, and we will refer to the KN-stratification induced by $L$ as the $L$-stratification.

For each $i$ there is a distinguished one parameter subgroup $\lambda_i$ of $G$. If we let $L_i$ be the centralizer of $\lambda_i$, then there is a smooth open subvariety $Z_i \subset X^{\lambda_i}$ which is $L_i$-invariant. Then by definition we have
$$S_i := G \cdot \left\{ x \in X | \lim_{t\to 0} \lambda_i(t) \cdot x \in Z_i \right\}$$
When the KN stratification arises from GIT, then in fact $Z_i$ is the semistable locus for the action of $L_i' = L_i / \lambda_i(\Gm)$ on the closure of $Z_i$.

The main object of study in this paper will be quotients stacks admitting a KN stratification of the following form:
\begin{defn} \label{def:semicomplete_KN}
A KN-stratification of a quotient stack $X/G$ is \emph{semi-complete} if $X^{ss}/G$ and $Z_i / L_i'$ all admit good quotients which are projective-over-affine. We say that the KN-stratification is \emph{complete} if all of the qood quotients are projective.
\end{defn}
\begin{rem} \label{rem:refining_strats}
Given a KN-stratification of a $G$-scheme $X$, if $X^{ss}/G$ and $Z_i/L'$ all admit semi-complete (resp. complete) KN stratifications, then the stratification of $X$ can be refined to a semi-complete (resp. complete) KN stratification by replacing each stratum with the preimage of the strata of $Z_i/L_i'$ under the projection $S_i/G \to Z_i/L_i'$ and taking the distinguishing one-parameter subgroup of each of these new strata to be $\lambda_i$ plus a very small rational multiple of the distinguished one-parameter subgroup of the corresponding stratum in $Z_i/L_i'$ (which can be lifted to $L$ rationally).
\end{rem}
In a sense the main theorem of GIT is the following:
\begin{thm}
Given a reductive $G$ and any $G$-ample bundle on a projective-over-affine $G$-scheme $X$, the $L$-stratification is semi-complete.
\end{thm}

Semi-complete KN stratifications are important because they lead to direct sum decompositions of noncommutative motives. If $I$ is a (possibly infinite) totally ordered set and $\A$ is a pretriangulated dg-category, we say that $\A = \sod{\A_i ; i\in I}$ forms a semiorthogonal decomposition if objects of the full pre-triangulated dg-subcategories $\A_i$ generate $\A$ under cones and shifts, and $\RHom(\A_i,\A_j) = 0$ for $i > j$. In other words, a semiorthogonal decomposition of a pre-triangulated dg-category is by definition a semiorthogonal decomposition of its homotopy category.

\begin{lem} \label{lem:KN_stratification}
If $X$ is a smooth $G$-scheme with a KN stratification, we have an equivalence in $\Madd_k$
$$\cU_k(\Perf(X/G)) \simeq \cU_k(\Perf(X^{ss}/G)) \oplus \bigoplus_i \cU_k(\Perf(Z_i/L_i)).$$
\end{lem}

\begin{proof}
The main theorem of \cite{halpern2015derived} provides an infinite semiorthogonal decomposition of $\Perf(X/G)$ under these hypotheses. One factor of the semiorthogonal decomposition is equivalent to $\Perf(X^{ss}/G)$, and the rest are of the form $\Perf(Z_i / L_i)_w$, where the subscript denotes the full subcategory of objects whose homology sheaves are concentrated in weight $w$ with respect to $\lambda$. The fact that $\cU_k$ commutes with filtered colimits implies that the infinite semiorthogonal decomposition maps to an infinite direct sum decomposition of $\cU_k(\Perf(X/G)) \in \Madd_{k}$. On the other hand, the category $\Perf(Z_i/L_i)$ decomposes as a direct sum of the subcategories $\Perf(Z_i/L_i)_w$ over all $w \in \bZ$, so $\bigoplus_w \cU_k(\Perf(Z_i/L_i)_w) \simeq \cU_k(\Perf(Z_i/L_i)) \in \Madd_k$.
\end{proof}

We will also use KN stratifications to compare properness of the dg-category $\Perf(X/G)$ to properness of the dg-category $\Perf(X^{ss}/G)$ and $\Perf(Z_i/L_i)$ for all $i$.

\begin{lem} \label{lem:properness}
Let $\X$ be a perfect derived $k$-stack of finite cohomological dimension, and let $\Y$ be another perfect derived $k$ stack such that $\Y^{cl,red} \simeq \X^{cl,red}$ and $\cO_\Y$ is eventually co-connective. Then the following are equivalent
\begin{enumerate}
\item $H_i \RGamma(\X,F)$ is finite dimensional for all $i$ and all $F \in \APerf(\X)$,
\item $\RGamma(\X,F)$ is finite dimensional for all $F \in \Coh(\X)$, 
\item $\Perf(\Y)$ is a proper dg-category.\footnote{We will need to consider the derived critical loci below, which is why we have introduced derived stacks here. If $\X$ is classical, then there is no need to replace $\X$ by an eventually co-connective approximation in (3), but the example $\Y=\X=\op{Spec}(k[t])$, where $t$ is a variable of homological degree $2$, shows that $(2)$ does not imply $(3)$ without the eventually co-connective hypothesis.}
\end{enumerate}
Furthermore, if $\X$ is a separated DM stack then this is equivalent to $\X$ being proper.
\end{lem}

\begin{proof}
Finite cohomological dimension implies that for any $F \in \APerf(\X)$ and all $i \in \bZ$, there is a sufficiently high $n$ such that $H_i \RGamma(\X,\tau_{\leq n} F) \simeq H_i \RGamma(\X,F)$, so $(2) \Rightarrow (1)$. Also, $(1) \Rightarrow (2)$ because $\Coh(\X) \subset \APerf(\X)$ and $H_i \RGamma(F)$ vanishes in all but finitely many degrees. It is clear that $(2)$ can be checked on $\X^{cl,red}$ because every $F \in \Coh(\X)$ is pushed forward from $\X^{cl}$, and any $F \in \Coh(\X^{cl})$ has a finite filtration whose associated graded is pushed forward from $\X^{cl,red}$.

To show that $(2) \Leftrightarrow (3)$, it thus suffices to show that $(2)$ is equivalent to $\Perf(\X)$ being a proper dg-category in the case when $\X$ is eventually co-connective. Because $\X$ is perfect, for any $F \in \DCoh(\X)$ and any $n$ we can find a perfect complex $P$ such that $F$ is a retract of $\tau_{\leq n}P$, so choosing $n$ large enough shows that $H_i \RGamma(\X,F)$ is a retract of $H_i \RGamma(\X,P)$, which is finite if $\Perf(\X)$ is a proper dg-category. On the other hand, $\Perf(\X) \subset \DCoh(\X)$ if $\X$ is eventually co-connective, so $\Hom_\X(E,F) = \RGamma(E^\dual \otimes F)$ is finite dimensional for perfect complexes $E$ and $F$.

For the further claim, it suffices to assume that $\X$ is classical. In this case if $\X$ is a separated DM stack, one may find a proper surjection from a quasi-projective scheme $X \to \X$ \cite{olsson2005proper}, and then deduce that $X$ is proper from property $(2)$, and hence $\X$ is proper.
\end{proof}

\begin{lem} \label{lem:KN_strat_properness}
Let $X/G$ be a quotient stack with a KN stratification. Then $\Perf(X/G)$ is a proper dg-category if and only if $\Perf(X^{ss}/G)$ and $\Perf(Z_i/L'_i)$ are proper dg-categories for all $i$.
\end{lem}

\begin{proof}
It suffices to consider the case of a single closed stratum $S \subset X$ with center $Z \subset S$ and with open complement $U$.

First assume that $\Perf(X/G)$ is a proper dg-category.  \cite{halpern2015remarks}*{Theorem 2.1} a fully faithful embedding $\Perf(U/G) \subset \APerf(X/G)$ (in fact one for each choice of $w \in \bZ$), and to prove the lemma it will suffice by \autoref{lem:properness} to show that this embedding preserves $\RGamma$. We will adopt the notation of \cite{halpern2015remarks}: this amounts to showing that we can choose a $w$ such that for $F \in \cG^w \subset \APerf(\X)$, which is identified with $\APerf(\X^{ss})$ under restriction, we have $\RGamma(\X,F) \simeq \RGamma(\X^{ss},F)$. This holds for $w=0$ by \cite{halpern2015remarks}*{Lemma 2.8}

Regarding $X$ as a derived stack, we may define the derived fixed locus $\tilde{Z} / L$, whose underlying classical stack is $Z/L$. Then \cite{halpern2015remarks}*{Theorem 2.1} shows that the functor
$$i_\ast \pi^\ast : \APerf(\tilde{Z}/L') \simeq \APerf(\tilde{Z}/L)^0 \to \APerf(X/G)$$
is fully faithfull. By \autoref{lem:properness} the dg-category $\Perf(\tilde{Z}/L')$ is proper, and thus so is $\Perf(Z/L')$.

Conversely, assume that $\Perf(Z/L')$ and $\Perf(U/G)$ are both proper dg-categories. We will show that $\Perf(X/G)$ is proper by invoking \autoref{lem:properness} and showing that $H_n \RGamma(X,F)^G$ is finite dimensional for any $n$ and any coherent sheaf $F$. Again by \cite{halpern2015remarks}*{Theorem 2.1}, we can functorially write $F$ as a finite extension of an object $F' \in \cG^0$ and two objects supported on the unstable stratum $\cS = S/G$, one in $\APerf_{\cS}(\X)^{\geq 0}$ and one in $\APerf_{\cS}(\X)^{<0}$. In particular as noted above we have $\RGamma(\X,F') \simeq \RGamma(U/G,F')$, which has finite dimensional homology.

Thus it suffices to show that $\RGamma(\X,F'')$ has finite dimensional homology for any $F'' \in \APerf(\X)$ which is set theoretically supported on $\cS$. Because $X$ has finite cohomological dimension, we may truncate $F''$ so that it lies in $\DCoh(\X)$, and then in can be built out of a sequence of extensions of shifts of objects of the form $i_\ast E$ for $E \in \Coh(S/G)$. Thus it suffices to show that $\Perf(S/G)$ is proper. A similar filtration argument using the baric decomposition of \cite{halpern2015remarks}*{Lemma 2.2} can be used to deduce that $\Perf(S/G)$ is proper because $\Perf(Z/L)$ is proper. Finally, the projection $Z/L \to Z/L'$ is a $\Gm$-gerbe, so the pushforward preserves perfect complexes, and thus $\Perf(Z/L)$ is proper if $\Perf(Z/L')$ is proper.
\end{proof}

\begin{cor} \label{cor:KN_strat_properness}
Let $X/G$ be a quotient stack with a semi-complete KN stratification. Then $\Perf(X/G)$ is a proper dg-category if and only if the stratification is complete.
\end{cor}
\begin{proof}
Combine \autoref{lem:KN_strat_properness} with \autoref{lem:properness}.
\end{proof}

\subsection{Motivic decompositions via KN stratifications}

We will consider the class of stacks which have semi-complete KN stratifications as in \autoref{def:semicomplete_KN}. We use the notation $\cC^{\oplus \bN}$ to denote the direct sum of countably many copies of the dg-category $\cC$. Recall also the definition of \cite{kresch2009geometry} that a DM stack of finite type with finite inertia over a field of characteristic $0$ is \emph{quasi-projective} if $\cX$ is a global quotient stack and has a quasi-projective coarse moduli space. We will say that $\cX$ is furthermore \emph{projective-over-affine} if its coarse moduli space is projective over an affine variety.

\begin{thm} \label{thm:motivic_main}
Let $G$ be an algebraic group over a field $k$ of characteristic $0$. Let $X$ be a smooth $G$-quasiprojective $k$-scheme with a semi-complete KN stratification. Then there is a smooth projective-over-affine Deligne-Mumford stack $\Y$ such that $\cU_k( \Perf(X/G))$ is a direct summand of $\cU_k(\Perf(\Y))^{\oplus \bN}$ in $\Madd_k$. Furthermore if $\Perf(X/G)$ is a proper dg-category, then $\Y$ can be chosen to be a smooth projective scheme.\footnote{The original version of this paper had $\Y$ as a projective DM stack, but by subsequent work of Bergh, Lunts, and Schnuerer \cite{bergh2016geometricity} one can further reduce to a smooth projective scheme.}
\end{thm}

Note that by \autoref{cor:KN_strat_properness}, $\Perf(X/G)$ is proper if and only if the KN stratification is complete.

\begin{rem} \label{rem:proof}
The proof is constructive, and actually produces something a bit stronger: if $\cC$ is the $\infty$-category of small dg-categories, then $\Perf(X/G)$ lies in the smallest subcategory containing $\Perf(\Y)$ and closed under countable semiorthogonal gluings and passage to semiorthogonal factors.
\end{rem}

\begin{ex} \label{ex:projective_over_affine}
If $X$ is projective-over-affine with a linearizable $G$-action, then the condition that $\Perf(X/G)$ is a proper dg-category is equivalent to the condition that $H_0 \RGamma(X,\cO_X)^G$ is finite dimensional, by \cite{halpern2014mapping}*{Proposition 4.2.3}.
\end{ex}

\begin{ex}
We can write any algebraic $k$-group $G$ as a semidirect product $G = U \rtimes L$, where $U$ is its unipotent radical and $L$ its reductive quotient. Assume that there is a one-parameter subgroup $\lambda : \Gm \to L$ which is central in $L$ and acts with positive weights on $\op{Lie}(U)$ in the adjoint representation of $G$. Then this one-parameter subgroup defines a single KN stratum $S = X = \{\pt \}$, and $Z/L = \pt / L' \to \pt$ is a good quotient. Thus \autoref{thm:motivic_main} applies to a large class of categories of the form $\Perf(BG)$, including when $G$ is a parabolic subgroup of a reductive group.
\end{ex}

\begin{ex}
If $G$ is as in the previous example, and $X$ is a smooth projective-over-affine $G$-scheme, then one can consider the Bialynicki-Birula stratification of $X$ under the action of $\lambda(\Gm)$, which is a KN stratification. If this is exhaustive, and $\Gamma(X^{\lambda(\Gm)}, \cO_{X^{\lambda(\Gm)}})^L$ is finite dimensional, then the Bialynicki-Birula stratification can be refined to a complete KN stratification of $X$ as in \autoref{rem:refining_strats}.
\end{ex}

Our proof of \autoref{thm:motivic_main} will proceed by a delicate inductive argument. One of the key tools is the following:

\begin{lem} \label{lem:rational_morphism}
Let $\pi : \Y \to \X$ be a rational morphism of finite-type $k$-stacks, meaning $R\pi_\ast \cO_\Y \simeq \cO_\X$. Assume that $\X$ is smooth and $\pi_\ast$ preserves $\DCoh$. Then $\cU_k(\Perf(\X))$ is a summand of $\cU_k(\Perf(\Y))$ in $\Madd_k$.
\end{lem}

\begin{proof}
First consider the categories $\Perf(\Y)$ and $\Perf(\X)$. The unit of adjunction $\op{id}_{X/G} \to \pi_\ast \pi^\ast$ is an equivalence in $\Perf(\X)$, hence $\pi^\ast$ is fully faithful and admits a right adjoint. Hence $\Perf(\X)$ is a semiorthogonal factor of $\Perf(\Y)$.
\end{proof}

We will apply \autoref{lem:rational_morphism} in three different situations.

\begin{ex}
If $\pi : \Y \to \X$ is a flat morphism of algebraic stacks such that for every $k$-point $p$ of $\X$ the fiber $\Y_p$ satisfies $\RGamma(\Y_p,\cO_{\Y_p}) \simeq k$, then $\pi$ is rational. If $\pi$ is not flat, then the same is true if we take $\Y_p$ to refer to the derived fiber.
\end{ex}

\begin{ex}
Any representable birational morphism of smooth $k$-stacks is rational. Indeed we can reduce this to the case for schemes, as birational morphisms are preserved by flat base change and the property of a morphism being rational is fppf-local on the base.
\end{ex}

\begin{ex}
Let $G \to H \to K$ be an extension of linearly reductive groups, and let $K$ act on a scheme, $X$. Then the morphism $p : X/H \to X/K$ is a $G$-gerbe -- after base change to $X$ this morphism becomes the projection $X \times BG \to X$. Thus because $G$ is linearly reductive $Rp_\ast \cO_{X/H} \simeq \cO_{X/K}$.
\end{ex}

Let $\pi : X' \to X$ be a projective morphism of smooth projective-over-affine varieties which is equivariant with respect to the action of a reductive group, $G$. For a $G$-ample invertible sheaf $L$ on $X$ and a relatively $G$-ample invertible sheaf $M$ on $X'$, we consider the fractional polarization $L_\epsilon = L + \epsilon M$ for $\epsilon \in \bQ$. We will need the following:

\begin{lem} \cite{teleman2000quantization}*{Lemma 1.2} \label{lem:refinements}
For any small positive $\epsilon \in \bQ$, the $L_\epsilon$-stratification of $X'$ refines the preimage of the $L$-stratification of $X$.
\end{lem}

Finally, we need another GIT lemma:
\begin{lem} \label{lem:reduce_to_GIT}
Let $X$ be a $G$-quasi-projective scheme which admits a good quotient $\pi : X \to Y$ such that $Y$ is projective-over-affine. Then $X = \bar{X}^{ss}$ for some linearized projective-over-affine $G$-scheme $\bar{X}$, which can be chosen to be smooth if $X$ is smooth.
\end{lem}
\begin{proof}
The proof of \cite{teleman2000quantization}*{Lemma 6.1} applies verbatim: one constructs a relative $G$-compactification for $X \to Y$ by choosing a sufficiently large coherent $F \subset \pi_\ast \cO_X$ so that $X$ embeds in the projectivization of $\op{Spec}_Y \op{Sym}(F)$. The closure of $X$ is projective over $Y$, and hence projective-over-affine, and it has a linearization for which $\bar{X}^{ss} = X$ by the cited argument. Furthermore, one can equivariantly resolve any singularities occuring in $\bar{X} \setminus X$ if $X$ is smooth.
\end{proof}

\begin{proof}[Proof of \autoref{thm:motivic_main}]
Over the course of the proof, we will actually construct a finite set of smooth quasi-projective DM stacks $\Y_1,\ldots,\Y_N$ such that $\cU_k(\Perf(X/G))$ is a retract of $\cU_k(\Perf(\Y_1))^{\oplus \bN} \oplus \cdots \oplus \cU_k(\Perf(\Y_N))^{\oplus \bN}$, and then we may take $\Y = \Y_1 \bigsqcup \cdots \bigsqcup \Y_N$ at the end. Also, when $\Perf(X/G)$ is a proper dg-category, it suffices to show that the $\Y_j$ can be chosen to be smooth proper DM stacks, because \cite{bergh2016geometricity}*{Theorem 6.6} implies that for any smooth proper DM stack $\Y$, $\Perf(\Y)$ is geometric, meaning a semiorthogonal factor of $\Perf$ of a smooth projective variety. We shall prove the theorem by induction on the rank of $G$.

Note that by \autoref{lem:KN_stratification} and the definition of a semi-complete KN-stratification, it suffices to prove this for quotient stacks which have projective-over-affine good quotients. For our purposes, it will be more convenient to consider smooth $G$-schemes which are projective-over-affine, and by \autoref{lem:reduce_to_GIT} and \autoref{lem:KN_stratification} it suffices to prove the claim for open unions of KN strata in a quotient stack of this form. We fix a $G$-ample bundle $L$ on $X$ and consider the $L$-stratification as in \autoref{eqn:KN_stratification}.

\medskip
\noindent  \textit{Case $X^{ss} = \emptyset$:}
\medskip

By \autoref{lem:KN_stratification} we must verify the claims for $\cU_k (\Perf(Z_i/L_i))$ for all $i$ for which $Z_i \subset U$. First assume that the inclusion $\lambda(\Gm) \subset L_i$ admits a splitting $L_i \to \Gm$, so that $L_i \simeq \Gm \times L_i'$ where the left factor is $\lambda(\Gm)$. Then $Z_i/L_i \simeq B\Gm \times Z_i/L'_i$, so $\cU_k(\Perf(Z_i/L_i))$ is a direct sum of copies of $\cU_k(\Perf(Z_i/L'_i))$. This is the only point of the proof at which an infinite direct sum enters, and it is an infinite direct sum of copies of the same category, hence throughout the proof we will only encounter a finite set of distinct DM stacks.

If $\lambda(\Gm) \subset L_i$ is not split, then we can choose a surjective homomorphism $\tilde{L} \to L_i$ with finite kernel, where $\tilde{L} \simeq \Gm \times L^\prime$ and $\Gm \times \{1\} \to L_i$ factors through $\lambda(\Gm)$. The morphism $p : Z_i / \tilde{L} \to Z_i / L_i$ is rational, hence \autoref{lem:rational_morphism} reduces the problem to showing the claim for $Z_i/\tilde{L}$. By the argument of the previous paragraph it again suffices to prove the claim for $Z_i/L_i'$.

Let $\bar{Z}_i$ be the closure of $Z_i$, which is a connected component of $X^{\lambda_i}$ and hence smooth and projective-over-affine. Then $Z_i$ is the semistable locus for the action of $L_i'$ on $\bar{Z}_i$, and $L'_i$ has lower rank than $L_i$, so the first claim of the theorem follows from the inductive hypothesis. Note that \autoref{lem:KN_strat_properness} implies that $\Perf(Z_i/L_i')$ is a proper dg-category if $\Perf(U/G)$ is, so in this case the inductive hypothesis implies that one can choose the DM stacks $\Y_j$ to be projective.

\medskip
\noindent \textit{Case $X^s = X^{ss} \neq \emptyset$:}
\medskip

The argument in the case where $X^{ss} = \emptyset$ applies here as well, so the inductive hypothesis implies that the conclusion of the theorem holds for $\Perf(Z_i / L_i)$ for all $Z_i \subset U$. By \autoref{lem:KN_stratification} it suffices to show that the claims hold for $X^{ss} / G$. In this case $X^{ss}/G$ is a smooth separated Deligne-Mumford stack whose coarse moduli space is projective-over-affine. Furthermore if $\Perf(U/G)$ is a proper dg-category then so is $\Perf(X^{ss}/G)$ by \autoref{lem:KN_strat_properness}, and hence $X^{ss}/G$ is projective by \autoref{lem:properness}.

\medskip
\noindent \textit{Case $X^{ss} \neq \emptyset$, $X^{ss} \neq X^s$, and $\op{codim}(X^{ss} \setminus X^s, X^{ss}) \geq 2$:}
\medskip

As in the previous case, it suffices to show the claims for $X^{ss}/G$. Here we use the main result of \cite{kirwan1985partial}, which says that there is a birational morphism $\pi : X' \to X$ such that $(X')^{ss}(L_\epsilon) = (X')^s(L_\epsilon)$, where $L_\epsilon = \pi^\ast L + \epsilon M$ for a suitable relatively $G$-ample $M$. By \autoref{lem:refinements} the open subset $U' := \pi^{-1}(X^{ss}(L))$ is a union of KN strata, and $\pi : U' \to X^{ss}(L)$ is rational, so by \autoref{lem:rational_morphism} we may reduce the main statement of the theorem to the corresponding claim for $U'/G$, which falls under the previous case. Furthermore, if $\Perf(U/G)$ is a proper dg-category the fact that $U'/G \to U/G$ is proper implies that $\Perf(U'/G)$ is a proper dg-category, so again we may reduce to the previous case to show that the $\Y_j$ can be chosen to be projective.

\medskip
\noindent \textit{Case $X^{ss} \neq \emptyset$ but $\op{codim}(X^{ss}\setminus X^{s},X^{ss}) \leq 1$:}
\medskip

Let $Y$ be a smooth projective variety with a $G$-action such that $\RGamma(Y,\cO_Y) \simeq k$ and for some linearization $M$ we have $\op{codim}(Y \setminus Y^s,Y) \geq 2$. For instance, $Y$ could be a suitable product of flag varieties, or a large projective space with a suitable linear $G$ action. We linearize the $G$ action on $X \times Y$ with $L_\epsilon = L + \epsilon M$. By \autoref{lem:refinements}, the open subvariety $U \times Y \subset X \times Y$ is a union of KN strata for the $L_\epsilon$ stratification. The projection $U \times Y / G \to U / G$ is rational, and so by \autoref{lem:rational_morphism} it suffices to prove the claims for $U \times Y / G$. Note that $\Perf(U \times Y /G)$ is a proper dg-category if $\Perf(U/G)$ is. Finally we have a sequence of inclusions
\[
X^{ss} \times Y^{s} \subset (X\times Y)^s \subset (X \times Y)^{ss} \subset X^{ss} \times Y,
\]
where the first inlcusion is due to the fact that points in $Y^{s}$ have finite stabilizers and thus so do points in $X^{ss} \times Y^s$, and the last inclusion follows from \autoref{lem:refinements}. This implies that
\[
\op{codim}((X \times Y)^{ss} \setminus (X \times Y)^{s}, (X \times Y)^{ss})  \geq \op{codim}(X^{ss}\times (Y \setminus Y^s) , X^{ss} \times Y) = \op{codim}(Y \setminus Y^s,Y) \geq 2,
\]
which reduces us to the previous case.

\end{proof}

\subsection{\texorpdfstring{$\Lambda$}{Lambda}-modules and the noncommutative Hodge-de Rham sequence} \label{sect:hochschild}

Let us recall the negative cyclic and periodic cyclic homology of a small $k$-linear dg-category, $\A$. We let $\hoch(\A) \in D(\Lambda)$ denote the (mixed) Hochschild complex of $\A$, regarded as a dg-module over $\Lambda = k[B] / B^2$ where $B$ has homological degree $1$ and acts on $\hoch(\A)$ by the Connes differential. We have
\begin{align*}
\hoch^{(n)} (\A) &:= \hoch(\A) \otimes k[u] / (u^n) \\
\hoch^{-} (\A) &:= \varprojlim_n \hoch^{(n)} \\
\hoch^{per} (\A) &:= \hoch^-(\A) \otimes_{k\ps{u}} k\ls{u}
\end{align*}
where $u$ is a variable of homological degree $-2$. The differential on each complex is given by $d + u B$, where $d$ is the differential on $\hoch(\A)$. In fact, these constructions make sense for any $\Lambda$-module $M$. We sometimes denote the negative cyclic construction $M^{\Sone}$ and the periodic cyclic construction $M^\Tate$. See \autoref{lem:Tate} below.

\begin{defn}[\cite{kontsevich2009notes}]
The category $\A$ is said to have the \emph{degeneration property} if $H_\ast(\hoch^{(n)}(\A))$ is a flat $k[u] / (u^n)$-module for all $n \geq 1$. 
\end{defn} 

It is immediate from the definitions that the degeneration property is preserved by filtered colimits of dg-categories. It is also known that the degeneration property holds for categories of the form $\A=\Perf(R)$, where $R$ is a smooth and proper dg-algebra \cites{kaledin2008ncHodge,MR3702384}. In particular, this holds when $\A = \Perf(\X)$, where $\X$ is a smooth and proper Deligne-Mumford stack over $k$ \cite{hallrydh2014DM}, although a more direct argument in this case follows from \autoref{prop:Toen_de_rham} below.  If $\A$ satisfies the degeneration property, then $H_\ast(\hoch^{-}(\A))$ is a flat $k\ps{u}$-module. (See \cite{kontsevich2009notes}*{Corollary 9.1.3}).

The degeneration property owes its name to its relationship with the noncommutative Hodge-de Rham spectral sequence. This is the spectral sequence associated to the filtration\footnote{Note that the filtration is not a filtration of $k\ls{u}$-modules, as $u \cdot F^p \subset F^{p+1}$. As explained to us by D. Kaledin, this can be understood by thinking of $u$ as the Tate motive. In other words when $k \subset \bC$, rather than regarding $k\ps{u}$ simply as a complex (where $u$ has cohomological degree 2), we regard it as $H^\ast(\mathbb{P}^\infty)$ with its Hodge structure, which places $u$ in $F^1 \bC\ps{u}$.} of the complex
$$F^p \hoch^{per}(\cA) = u^p \cdot \hoch^-(\cA) \subset \hoch^{per}(\cA).$$
The $E_1$ page of the spectral sequence is $\op{gr} \hoch^{per}(\cA) \simeq \hoch(\cA)\otimes k\ls{u}$. The degeneration property implies that this spectral sequence degenerates on the first page, so the associated graded of the resulting filtration on $H_\ast(\hoch^{per}(\cA))$ is isomorphic to $H_\ast(\hoch(\cA)) \otimes k\ls{u}$. Under the assumption that $\A$ is suitably bounded, we can say something more precise: 

\begin{lem} \label{lem:cyclic_complex_colimits}
Let $\A$ be a dg-category such that $H_\ast(\hoch(\A))$ is homologically bounded above and which satisfies the degeneration property. Then there exists an (non-canonical) isomorphism $H_\ast(\hoch^{-}(\A)) \simeq H_\ast(\hoch(\A)) \otimes k\ps{u}$.
\end{lem}
\begin{proof}
This follows from the remark before Theorem 4.14 of \cite{KKP}. 
\end{proof}
The hypothesis of \autoref{lem:cyclic_complex_colimits}, that $H_\ast(\hoch(\A))$ is homologically bounded above, will apply to $\Perf(\X)$ for all smooth $k$-stacks $\X$ of finite cohomological dimension such that $\QC(\X)$ is compactly generated.

\begin{rem}
The $\Lambda$-module $\hoch(\cA)$ is functorial in $\cA$. When $\cA$ is a symmetric monoidal $k$-linear $\infty$-category, exterior tensor product followed by the symmetric monoidal product gives a natural map
\[
\hoch^-(\cA) \otimes \hoch^-(\cA) \to \hoch^-(\cA \otimes \cA) \to \hoch^-(\cA)
\]
and likewise for $\hoch^{per}(\cA)$. On the level of homology, this gives $H_\ast \hoch^{per}(\cA)$ the structure of a commutative $k\ls{u}$-algebra, and for any symmetric monoidal functor between symmetric monoidal $\infty$-categories $\cA \to \cB$, the resulting map $H_\ast \hoch^{per}(\cA) \to H_\ast \hoch^{per}(\cB)$ is a map of commutative $k\ls{u}$-algebras.
\end{rem}


\subsection{The degeneration property for quotient stacks}
\label{sect:degeneration}

In \cite{teleman2000quantization}*{Theorem 7.3}, Teleman establishes the degeneration of a commutative Hodge-de Rham sequence, which converges to the equivariant Betti-cohomology $H^\ast_G(X)$, for a smooth quotient stack $X/G$ with a complete KN stratification. The argument in \cite{teleman2000quantization} makes use of the KN stratification and has a similar flavor to the proof of \autoref{thm:motivic_main}. However the proof in the commutative case is substantially simpler. In the noncommutative situation, we are not aware of an argument to reduce the proof of degeneration to the case of the quotient of a smooth projective scheme by the action of a reductive group, as was done in \cite{teleman2000quantization}.

Using the motivic statement of \autoref{thm:motivic_main}, we can immediately deduce noncommutative HdR degeneration. The main observation is the following
\begin{lem} \label{lem:degeneration_sums}
The degeneration property is closed under direct summands and arbitrary direct sums in $\Madd_k$.
\end{lem}
\begin{proof}
The Hochschild complex $C_\bullet(-)$ is an additive invariant of dg-categories, hence factors through $\cU_k$ uniquely up to contractible choices. The claim follows from the fact that the operation $D(\Lambda) \to D(k)$ taking mapping $(M,d,B) \mapsto (M\otimes k[u]/(u^n),d+uB)$ commutes with filtered colimits and in particular infinite direct sums, and the fact that an infinite direct sum of $k[u]/(u^n)$ modules is flat if and only if every summand is flat.
\end{proof}

\begin{cor} \label{cor:degeneration_complete}
Let $G$ be a reductive group and let $X$ be a smooth $G$-quasiprojective scheme which admits a complete KN stratification. Then $\Perf(X/G)$ has the degeneration property.
\end{cor}
\begin{proof}
Combine \autoref{lem:degeneration_sums} with the conclusion of \autoref{thm:motivic_main}.
\end{proof}

\begin{ex} \label{ex:bga}
As a counterexample, consider $\Perf(B \mathbb{G}_a)$. This category is Morita equivalent to the category $\operatorname{Perf}(k[\epsilon]/(\epsilon^2))$ where $\epsilon$ has degree $-1$. By the (graded-commutative) HKR theorem, proposition 5.4.6 of \cite{loday1992}, we have that $H_*\hoch(\Perf(k[\epsilon]/(\epsilon^2)) \cong k[\epsilon]/(\epsilon^2)\otimes \operatorname{Sym}^*(d\epsilon)$, where $d\epsilon$ has degree $0$. By theorem 5.4.7 of the same book, the Connes operator goes to the de Rham differential which sends $\epsilon \to d\epsilon$ and so the spectral sequence does not degenerate. 
\end{ex}

We also observe, somewhat surprisingly, that the derived category of coherent sheaves on certain \emph{singular} quotient stacks also has the degeneration property. We will consider the following goeometric set up
\begin{itemize}
\item $X /G = X^{ss}/G \cup \bigcup_i S_i / G$ is a complete KN stratification (\autoref{def:semicomplete_KN}) of a smooth quotient stack,
\item $V$ is a $G$-equivariant locally free sheaf on $X$ such that $V|_{Z_i}$ has $\lambda_i$-weights $\leq 0$ for all $i$, and
\item $\sigma \in \Gamma(X,V)^G$ is an invariant section.
\end{itemize}
Note that the quantization-commutes-with-reduction theorem \cite{teleman2000quantization} implies that if the $\lambda_i$-weights of $V|_{Z_i}$ are strictly negative, then $\Gamma(X,V)^G \simeq \Gamma(X^{ss}(L),V)^G$ (this is referred to as \emph{adapted} in \cite{teleman2000quantization}). Using the methods of \cite{halpern2015derived} one can show that $\dim \Gamma(X,V)^G < \infty$ even when the $\lambda_i$ weight of $V|_{Z_i}$ vanishes for some $i$.

\begin{amplif} \label{amplif:zero_locus}
In the set up above, if
\begin{enumerate}
\item $\sigma$ is regular on $X^{ss}$ with smooth vanishing locus, and
\item for all $i$ the restriction of $\sigma$ to $(V|_{Z_i})^{\lambda = 0}$, the summand of $V|_{Z_i}$ which is fixed by $\lambda(\Gm)$, is regular with smooth vanishing locus,
\end{enumerate}
then there is a smooth and proper quasi-projective DM stack $\Y$ such that
$\cU_k(\DCoh(X_\sigma / G))$ is a retract of $\cU_k(\Perf(\Y))^{\oplus \bN}$.
\end{amplif}

\begin{proof}
We apply the structure theorem for the derived zero locus $X_\sigma$ in \cite{halpern2015remarks}*{Theorem 3.2}, whose derived category is just the derived category of the sheaf of cdga's over $X/G$ given by the Koszul algebra
$$\cA = (\op{Sym}(V^\dual[1]),d\phi = \phi(s)).$$
The structure theorem constructs an infinite semiorthogonal decomposition which generalizes the main structure theorem of \cite{halpern2015derived}. One factor is isomorphic to $\DCoh(X^{ss}_\sigma/G)$, and the remaining factors are isomorphic to $\DCoh(Z'_i/L_i)^w$, where $Z'_i$ denotes the derived zero locus of $\sigma$ restricted to $(V|_{Z_i})^{\lambda = 0}$, and the superscript $w$ denotes the full subcategory of $\DCoh(Z'_i/L_i)$ consisting of complexes whose homology is concentrated in weight $w$.

In order to apply this theorem, we must check that after restricting the cotangent complex $\bL_{X_\sigma / G}$ to $Z'_i/L_i$ and looking at the summand with $\lambda$-weights $<0$, there is no fiber homology in homological degree $1$. Because $X_\sigma$ is a derived zero section, we have
\[
\left(\bL_{X_\sigma/G}|_{Z_i'} \right)^{\lambda <0} \simeq \left[ (V^\dual |_{Z_i})^{\lambda <0} \to (\Omega_X|_{Z_i'})^{\lambda <0} \to \cO_{Z_i'} \otimes (\fg^\dual)^{\lambda <0} \right].
\]
So the weight hypotheses on $V|_{Z_i}$ imply that this is a two term complex of locally free sheaves in homological degrees $0$ and $-1$, and hence has no fiber homology in homological degree $1$.

Given the structure theorem for $\DCoh(X_\sigma/G)$, the proof of \autoref{lem:KN_stratification} now applies verbatim to give a finite direct sum decomposition
$$\cU_k(\DCoh(X_\sigma/G)) =\cU_k(\DCoh(X^{ss}_\sigma / G)) \oplus \bigoplus_i \cU_k(\DCoh(Z_i'/L_i)).$$
Under the hypotheses of the amplification, each factor in this direct sum decomposition is $\DCoh$ of a smooth quotient stack satisfying the hypotheses of \autoref{thm:motivic_main}, and the result follows.
\end{proof}

\begin{rem}
Note that when $V$ is strictly adapted to the KN stratification, then the condition (2) in the previous amplification is vacuous.
\end{rem}

\begin{cor} \label{cor:zero_locus_degeneration}
In the set up of \autoref{amplif:zero_locus}, the category $\DCoh(X_\sigma / G)$ has the degeneration property.
\end{cor}

We will see another approach to establishing \autoref{cor:zero_locus_degeneration} using graded LG-models in \autoref{ex:graded_lg_model} and \autoref{prop:graded_lg_degeneration}.


\section{Hodge structures on equivariant \texorpdfstring{$K$}{K}-theory}

\newcommand{\Smooth}{\mathbf{Sm}}

In this section we consider the action of a reductive group $G$ on a smooth quasi-projective $\bC$-scheme $X$. Our goal is to identify the periodic cyclic homology $\hoch^{per}(\DCoh(X/G))$ with the complexification of the Atiyah-Segal equivariant topological $K$-theory $K_{M}(X^{an})$ with respect to a maximal compact subgroup $M \subset G$ (see Section \ref{sect:top_k_theory}). Our final result, \autoref{thm:Hodge_main}, will allow us to define a pure Hodge structure of weight $n$ on $K^n_{M}(X^{an})$ in the case where $X$ admits a complete KN stratification.

Rather than construct a direct isomorphism, we study an intermediate object, the topological $K$-theory of the dg-category $K^{top}(\DCoh(X/G))$, as defined in \cite{blanc2012topological}, which admits natural comparison isomorphisms with each of these theories. In Blanc's construction, $K^{top}(\cC)$ is constructed from the geometric realization of the presheaf of spectra on the category, $\Aff$, of affine $\bC$-schemes of finite type,
$$\underline{K}(\cC) : A \mapsto K(A \otimes_\bC \cC).$$
Here $K(-)$ denotes the non-connective algebraic $K$-theory of a dg-category of \cite{blanc2012topological}*{Definition 2.10}. The geometric realization of a presheaf, $|\bullet|$, is defined to be the left Kan extension of the functor $A \mapsto \Sigma^\infty (\op{Spec} A)^{an}_+$, regarded as functor with values in spectra, along the Yoneda embedding of the category of finite type $\bC$-schemes into presheaves of spectra, $\Aff \to \Sp(\Aff)$. The geometric realization functor $|\bullet| : \Sp(\Aff) \to \Sp$ admits a right adjoint, which assigns $M \in \Sp$ to the presheaf of spectra $H_{B}(M):=\Hom_{\Sp}(\Sigma^\infty (\bullet)^{an}_+,M)$. The semi-topological $K$-theory is the geometric realization
$$K^{st}(\cC) := |\underline{K}(\cC)|,$$
regarded as a $K^{st}(\bC)$-module spectrum. By \cite{blanc2012topological}*{Theorem 4.5}, we have an isomorphism $K^{st}(\bC) \simeq bu$, where the latter denotes the connective topological $K$-theory spectrum. Choosing a generator $\beta \in \pi_2 (bu)$, one then defines the topological $K$-theory of a dg-category to be
$$K^{top}(\cC) := K^{st}(\cC) [\beta^{-1}] = |\underline{K}(\cC)| \otimes_{bu} bu[\beta^{-1}]$$

\begin{rem}
In the definition of the presheaf $\underline{K}(\cC)$, $A \otimes_\bC \cC$ denotes the derived tensor product of non-idempotent-complete dg-categories, as in \cite{keller}*{Section 4.3}. In \autoref{sect:singularities} we will also consider a symmetric monoidal structure on small stable idempotent complete dg-categories, which we will denote $\cA \itimes \cB$ to avoid confusion. In addition, following the convention of \cite{blanc2012topological}*{Definition 2.7}, throughout this paper we define the algebraic $K$-theory of a small dg-category $\cC$ to be the algebraic $K$-theory of the Waldhausen category $\Perf(\cC)$, so that it is automatically Morita invariant. We warn the reader that as a consequence, if $\cC$ is a stable dg-category, then $K(\cC)$ denotes the usual algebraic $K$-theory of the idempotent completion of $\cC$, not the $K$-theory of $\cC$ itself.
\end{rem}

We will also use the construction of a Chern character map $\op{Ch} : K^{top}(\cC) \to \hoch^{per}(\cC)$.\footnote{In order to be consistent with the rest of the paper, we use the notation $\hoch^{per}$ for the periodic cyclic homology complex of a dg-category, rather than the notation $HP$ used in \cite{blanc2012topological}. In addition, we use the notation $\otimes$ rather than $\wedge$ for the smash product of spectra and module spectra. For example $K^{top}(\cC) \otimes \bC$ is the $\bC$-module spectrum, which we canonically identify with a complex of $\bC$-modules, which is denoted $K^{top}(\cC) \wedge H\bC$ in \cite{blanc2012topological}.} First, one obtains a map of presheaves $\underline{K}(\cC) \to \underline{\hoch^{per}}(\cC)$ from the usual Chern character in algebraic $K$-theory, where $\underline{\hoch^{per}}(\cC)$ denotes the presheaf $ A \mapsto \hoch^{per}(A \otimes_\bC \cC)$. Using a version of the Kunneth formula for periodic cyclic homology, one obtains an equivalence $|\underline{\hoch^{per}}(\cC)| \simeq \hoch^{per}(\cC) \otimes_{\bC[u^\pm]} |\underline{\hoch^{per}}(\bC)|$. Then one can construct an isomorphism of presheaves $\underline{\hoch^{per}}(\bC) \simeq H_B(\bC[u^\pm])$, which leads to a map $|\underline{\hoch^{per}}(\bC)| \to \bC[u^{\pm}]$. Combining these provides a map
$$K^{st}(\cC) \to \hoch^{per}(\cC) \otimes_{\bC[u^\pm]} |\underline{\hoch^{per}}(\bC)| \to \hoch^{per}(\cC)$$
which give the Chern character after inverting $\beta$. The main result we use is \cite{blanc2012topological}*{Proposition 4.32}, which states that for a finite type $\bC$-scheme, $X$, the Chern character induces an equivalence $K^{top}(\Perf(X)) \otimes \bC \to \hoch^{per}(\Perf(X))$. Furthermore, there is a natural topologization map which is an equivalence $K^{top}(\Perf(X)) \to K(X^{an})$, and under this equivalence $\op{Ch}$ can be identified with a twisted form of the usual Chern character for $X^{an}$ under a canonical isomorphism $\hoch^{per}(X) \to H_{Betti}(X;\bQ) \otimes_\bQ \bC[u^\pm]$. More precisely, Blanc's Chern character provides an equivalence
$$K^{top}(\Perf(X)) \otimes \bQ \simeq H_{Betti}(X;\bQ) \otimes \bQ\ls{\frac{u}{2 \pi i}} \subset H_{Betti}(X^{an};\bQ) \otimes \bC\ls{u},$$
which we can alternatively express as an isomorphism
$$\pi_n (K^{top}(\Perf(X))) \otimes \bQ \simeq \bigoplus_p H^{2p-n}(X^{an};\bQ \langle p \rangle ),$$
where $\bQ \langle p \rangle \subset \bC$ denotes the subgroup $(2\pi i)^{p} \bQ$.

\subsection{Equivariant \texorpdfstring{$K$}{K}-theory: Atiyah-Segal versus Blanc}

In this section we consider a reductive group $G$ with maximal compact subgroup $M \subset G$, and a $G$-quasi-projective scheme $X$, which need not be smooth. The goal of this section will be to construct a comparison isomorphism between $K^{top}(\DCoh(X/G))$ and topological $M$-equivariant $K$-homology of $X^{an}$ with locally compact supports, whose construction and properties we recall below.

\subsubsection{Equivariant Borel-Moore homology theories from invariants of dg-categories}

Let us consider an additive invariant of $dg$-categories $E(-)$. Then we have a Borel-Moore homology theory on the category of $G$-schemes by defining $E_{BM}(X) := E(\DCoh(X/G))$ (see \autoref{rem:borel_moore_homology} below). This assignment is covariantly functorial with respect to proper $G$-equivariant maps $p : X' \to X$ by applying $E(-)$ to the pushforward functor $p_\ast : \DCoh(X'/G) \to \DCoh(X/G)$, and $E(X)$ is contravariantly functorial with respect to flat $G$-equivariant maps $f : X' \to X$ by applying $E(-)$ to the pullback functor $f^\ast : \DCoh(X/G) \to \DCoh(X'/G)$. Note also that if $X$ is a smooth $G$-space, then we have a canonical equivalence $E(\Perf(X/G)) \to E(\DCoh(X/G))$, so we can canonically identify $E_{BM}(X)$ with the ``cohomology'' theory $E(X) := E(\Perf(X/G))$, which is a form of Poincar{\'e} duality.

More precisely, let $\Space_G$ denote the category of algebraic spaces with a $G$-action and $G$-equivariant maps, and let $\Space_G^{?}$, where $?=f,p,\simeq$, denote the subcategories consisting of flat maps, proper maps, and isomorphisms respectively. We regard the additive invariant $E$ as giving two strict functors into spectra $E^p_{BM} : \Space_G^{p} \to \Sp$ and $E^f_{BM}: (\Space_G^f)^{op} \to \Sp$ along with an isomorphism of the restriction of these functors to the subcategories $\Space_G^\simeq \subset \Space_G^p$ and $\Space_G^\simeq \subset (\Space_G^f)^{op}$. Here we have used the fact that $\Space_G^\simeq$ is a groupoid to canonically identify $\Space_G^\simeq \simeq (\Space_G^\simeq)^{op}$.

Recall that an additive invariant $E(-)$ of $dg$-categories is \emph{localizing} if it takes localization sequences of $dg$-categories to exact triangles, and we say that $E(-)$ satisfies equivariant d{\'e}vissage if for any closed immersion $i : X \hookrightarrow X'$ of $G$-spaces the induced map $i_\ast : E(\DCoh(X/G)) \to E(\DCoh_{X}(X'/G))$ is an equivalence, where $\DCoh_X(X'/G)$ denotes the derived category of quasi-coherent complexes with coherent cohomology sheaves which are set theoretically supported on $X$.

\begin{lem} \label{lem:torsor_pullback}
Let $E$ be an additive invariant of $dg$-categories. Then:
\begin{enumerate}
\item The base change formula holds: if $p : X \to Y$ is proper, $f : Y' \to Y$ is flat, and $f'$ and $p'$ are base change of $f$ and $p$, then we there is an equality of compositions $p_\ast \circ f^\ast = (f')^\ast \circ (p')_\ast : E_{BM}(X) \to E_{BM}(Y')$ in the homotopy category of spectra $\op{Ho}(\Sp)$;
\item If $\cV$ is a $G$-equivariant vector bundle on $X$ and $\pi : \bP(\cV) \to X$ is the projection, then $\pi^\ast : E_{BM}(X) \to E_{BM}(\bP(\cV))$ followed by $(-)\otimes \cO_{\bP(\cV)}(k) : E_{BM}(\bP(\cV)) \to E_{BM}(\bP(\cV))$ is a split injection;
\item The previous maps, where $k$ ranges from $0,\ldots,n-1$, define a canonical equivalence $$E_{BM}(\bP(\cV)) \simeq E_{BM}(X)^{\oplus n}.$$
\end{enumerate}
Furthermore, if $E$ is localizing and satisfies equivariant d{\'e}vissage, then:
\begin{enumerate}
\setcounter{enumi}{3}
\item For any closed immersion of $G$-spaces $i : X \hookrightarrow X'$ with open complement $j :U \subset X'$, there is an exact triangle $E_{BM}(X) \xrightarrow{i_\ast} E_{BM}(X') \xrightarrow{j^\ast} E_{BM}(U) \to$;
\item If $\pi : V \to X$ is a $G$-equivariant torsor for $\cV$, then $\pi^\ast : E_{BM}(X) \to E_{BM}(V)$ is an equivalence.
\end{enumerate}
\end{lem}
\begin{proof}
The claim (1) follows immediately by applying $E(-)$ to the usual base change formula. Claims (2) and (3) follow from the fact that $E$ is an additive invariant and there is a canonical semiorthogonal decomposition $$\DCoh(\bP(\cV)/G) = \langle \pi^\ast(\DCoh(X/G)), \pi^\ast(\DCoh(X/G)) \otimes \cO_{\bP(\cV)}(1),\ldots,\pi^\ast(\DCoh(X/G)) \otimes \cO_{\bP(\cV)}(n-1)\rangle.$$

Now let us assume that $E$ satisfies d{\'e}vissage. (4) follows immediately from the d{\'e}vissage condition and the localization exact triangle $$E(\DCoh_{X}(X'/G)) \to E(\DCoh(X'/G)) \to E(\DCoh((X'-X)/G)) \to.$$
For (5), we note that Thomason's proof when $E(-)$ is algebraic $K$-theory in \cite{thomason1987algebraic}*{Theorem 4.1} generalizes to this context as well: the torsor $V$ is classified by an extension $0 \to \cV \to \cV' \to \cO_{X} \to 0$ of $G$-equivariant locally free sheaves on $X$. $V$ is isomorphic to the complement of the closed embedding $\cP(\cV) \hookrightarrow \cP(\cV')$. Using the direct sum decomposition of $E(\bP(\cV'))$ and $E(\bP(\cV))$ from (3) one can show that $E(i_\ast) : E(\bP(\cF)) \to E(\bP(\cW))$ is a split injection, and the pullback map gives an equivalence of the cofiber with $E(X)$. It follows from the localization sequence of (4) that the cofiber of $E(i_\ast)$ can is canonically identified with $E(\bP(\cV')-\bP(\cV))$ under the restriction map, so we have that pullback gives an equivalence $E(X) \simeq E(V)$.
\end{proof}

\begin{ex}
Algebraic $K$-theory of $dg$-categories as defined in \cite{blanc2012topological}*{Definition 2.10} is a localizing invariant and satisfies d{\'e}vissage, and it follows that $K^{top}(-)$ is localizing and satisfies d{\'e}vissage as well.
\end{ex}

\begin{rem} \label{rem:borel_moore_homology}
The theory $E_{BM}(X)$ gives a Borel-Moore homology theory in the sense of \cite{levine2007algebraic}*{Definition 2.1.2}, but whose source is the category of $G$-spaces and which takes values in $\op{Ho}(\Sp)$ rather than graded abelian groups. We are forced to use the homotopy category of spectra because we have only formulated the base change formula in this setting. In order to formulate a ``strict'' Borel-Moore homology theory valued in spectra, one would have to define the source category as an $\infty$-category of correspondences, as in \cite{gaitsgory2017study}*{Section V}.
\end{rem}

\begin{rem} \label{rem:borel_moore_stacks}
Although we are only interested in quotient stacks here, one can define a Borel-Moore homology theory $E_{BM}(\cX) := E(\DCoh(\cX))$ for arbitrary algebraic stacks, and the statements and proofs of the previous lemma extend to this context.
\end{rem}

\subsubsection{Atiyah-Segal equivariant $K$-theory} \label{sect:top_k_theory}

We now turn our attention to topological equivariant $K$-theory with respect to a maximal compact subgroup $M \subset G$. There is a $K$-cohomology theory for topological $M$-spaces constructed in \cite{atiyah1969equivariant}. Below we use the more systematic description of $K_{M}(X)$ in terms of equivariant stable homotopy theory as the spectrum obtained by taking level-wise $M$-equivariant mapping spaces from $X$ to the naive $M$-spectrum underlying the $M$-spectrum $bu_M$. For details on the non-equivariant and equivariant stable homotopy category, we refer the reader to \cite{lewis1986equivariant} and \cite{may1996equivariant}*{Chapters XII, XIV, and XVI}. $K_M(X^{an})$ can also be constructed as the $K$-theory spectrum associated via Quillen's $Q$-construction to the exact category of $M$-equivariant complex vector bundles on $X^{an}$.

We will also consider the Atiyah-Segal equivariant $K$-homology with locally compact supports $K_M^{c,\dual}(X^{an})$. This theory was studied in \cite{thomason1988equivariant}*{Section 5} under the notation $G^{AS}(G,X)$, and our discussion follows this reference closely. In particular, we refer the reader there for a nice discussion contextualizing $K_M^{c,\dual}(-)$ with respect to several other versions of equivariant $K$-theory. We have chosen to denote the $M$-equivariant $K$-homology with locally compact supports as $K^{c,\dual}_{M}(X^{an})$ because it is the $M$-equivariant Spanier-Whitehead dual (see \cite{may1996equivariant}*{XVI.7}) of the $M$-spectrum of equivariant $K$-theory with compact supports constructed in \cite{segal1968equivariant}*{Page 136}, which we denote $K^c_{M}(X^{an})$.

In particular, $K^c_M(X)$ is strictly covariantly functorial for open immersions \cite{segal1968equivariant}*{Proposition 2.9} and strictly contravariantly functorial for proper maps. Dually, it follows that $K^{c,\dual}_M(-)$ can be regarded as both a strict functor $\Space_G^{p} \to \Sp$ and a strict functor $(\Space_G^{o})^{op} \to \Sp$, where $\Space_G^{o} \subset \Space_G$ denotes the subcategory of open immersions, and these two functors are isomorphic on the subcategory $\Space_G^\simeq$. Poincar{\'e} duality in this context gives a canonical isomorphism $K_M(X^{an}) \simeq K_M^{c,\dual}(X^{an})$ for complex $G$-manifolds -- a priori this isomorphism depends on a choice of orientation, but complex manifolds are canonically oriented for $K$-theory (see \cite{segal1968equivariant}*{Section 3} and \cite{may1996equivariant}*{XIV}).

\begin{lem} \label{lem:torsor_pullback_segal}
The properties of \autoref{lem:torsor_pullback} also hold for $K_M^{c,\dual}(X^{an})$. Namely:
\begin{enumerate}
\item The base change formula holds with respect to a proper map $p : X \to Y$ and an open immersion $f : U \subset Y$;
\item If $\cV$ is a $G$-equivariant vector bundle on $X$ of dimension $n$, then there is a canonical equivalence $K_M^{c,\dual}(\bP(\cV)^{an}) \simeq K_M^{c,\dual}(X^{an})^{\oplus n}$.
\item For any closed immersion of $G$-spaces $i : X \hookrightarrow X'$ with open complement $j :U \subset X'$, there is an exact triangle $E_{BM}(X) \xrightarrow{i_\ast} E_{BM}(X') \xrightarrow{j^\ast} E_{BM}(U) \to$;
\item If $\pi : V \to X$ is a $G$-equivariant torsor for $\cV$, then there is a canonical equivalence $E_{BM}(X) \to E_{BM}(V)$.
\end{enumerate}
\end{lem}
\begin{proof}
The base change formula arises from the observation that if $U' := p^{-1}(U) \subset X$, then there is a commutative diagram of maps of one point compactifications
\[
\xymatrix{X^+ \ar[r] \ar[d] & (U')^+ \ar[d] \\ Y^+ \ar[r] & U^+},
\]
and it is precisely pullback along these maps which is used to define proper pullback and pushforward along an open immersion in $K_M^c(-)$. The canonical equivalence (2) is dual to the equivalence $K_M^c(\bP(\cV)) \simeq K_M^c(X)^{\oplus n}$ of \cite{segal1968equivariant}*{Proposition 3.9}. (3) is \cite{segal1968equivariant}*{Proposition 2.9}. One proves (4) using the same method as the proof of \autoref{lem:torsor_pullback}: the canonical functor $K_M^{c,\dual}(X^{an}) \to K_M^{c,\dual}(V)$ is the inclusion of the $0^{th}$ piece of the direct sum decomposition $K_M^{c,\dual}(\bP(\cV')) \simeq K_M^{c,\dual}(X^{an})^{\oplus n+1}$ followed by restriction to the open subset $V \subset \bP(\cV')$.
\end{proof}

\begin{rem}
It is a priori not clear how to define a canonical pullback map $f^\ast : K_M^{c,\dual}(Y^{an}) \to K_M^{c,\dual}(X^{an})$ for a flat map $f: X \to Y$. That is why we have only stated the projection formula for open immersions, and it's why in (2) and (4) above we can not state that the canonical isomorphisms are given by pullback. Indeed it is not immediate in (4) that the isomorphism $K_M^{c,\dual}(X^{an}) \simeq K_M^{c,\dual}(V^{an})$ is independent of the presentation of $V$ as a torsor for $\cV$. If $f$ is the restriction of a flat map $X' \to Y'$ along a closed embedding $Y \hookrightarrow Y'$, then one could provisionally define a pullback map using the localization sequence (3), but perhaps the most uniform definition of flat pullback for $K_M^{c,\dual}$ is via the isomorphism of \autoref{thm:topological_comparison} below.
\end{rem}

\subsubsection{The comparison map}

We now construct a comparison map between Borel-Moore homology theories $\rho_{G,X} : K^{top}(\DCoh(X/G)) \to K_{M}^{c,\dual}(X^{an})$ by first constructing a comparison isomorphism between the corresponding cohomology theories $K^{top}(\Perf(X/G)) \to K_M(X^{an})$. Observe that we have a natural transformation of presheaves of spectra
$$K^{alg}(\operatorname{Perf}(X/G\times T)) \to K_{M}(X^{an} \times T^{an})$$ 
which is induced by the functor of exact categories which sends an algebraic $G$-vector bundle to its underyling complex topological vector bundle equipped with the induced action of ${M}$ (This functor is symmetric monoidal, and hence induces a map of $K$-theory spectra \cite{thomason1988equivariant}*{Section 5.4}). By the \autoref{lem:equiv_homotopy} below, the presheaf $K_{M}(X^{an} \times (-)^{an})$ is weakly equivalent to $\Hom_{Sp}(\Sigma^\infty (-)_+^{an},K_{M}(X^{an}))$, where $\Hom_{\Sp}(\bullet,\bullet)$ denotes the internal function spectrum in the category of spectra. Composing this with the comparison map above gives a map of presheaves
$$K^{alg}(\operatorname{Perf}(X/G) \otimes_\bC \cO_T)) \to \Hom_{Sp}(\Sigma^{\infty} T^{an}_+, K_{M}(X^{an})),$$
where both sides are regarded simultaneously as presheaves in the $G$-space $X$ and the affine scheme $T$. Here we have used the natural Morita equivalence $\Perf(X/G)\otimes_\bC \cO_T \xrightarrow{\simeq} \Perf(X/G \times T)$. By the adjunction defining the geometric realization \cite{blanc2012topological}*{Definition 3.13} for presheaves on the category of affine schemes, this gives a map
\[
K^{st}(\Perf(X/G)) \to K_M(X^{an})
\]
of presheaves in the $G$-space $X$. This natural transformation of presheaves on $\Space_G$ will be $bu$-linear by construction. Because $K_M(X^{an})$ already satsifies Bott periodicity, this comparison map extends uniquely to the localization, giving our final comparison map
\[
\rho_{G,X}^{perf} : K^{top}(\Perf(X/G)) := K^{st}(\Perf(X/G)) \otimes_{bu} bu[\beta^{-1}] \to K_M(X^{an}).
\]

The following technical lemma was used in the construction above:

\begin{lem} \label{lem:equiv_homotopy}
For any topological space $Y$ and topological $M$-space $X$, we have a natural weak equivalence in $\Sp$,
$$\Hom_{Sp}(\Sigma^{\infty} Y_+ , K_{M}(X))  \cong K_{M}(X \times Y),$$
where on the right $Y$ is regarded as an $M$-space with trivial $M$ action.
\end{lem}

\begin{proof}
We fix a universe $U$ for forming the equivariant stable homotopy category $\Sp_M$ as in \cite{lewis1986equivariant}. The ``change of universe" functor taking an $E \in \Sp_M$ to its underlying naive $M$-spectrum admits a left adjoint, as does the functor from naive $M$-spectra to spectra which applies the $M$-fixed point functor level-wise. We will denote the composition of these to functors as $(-)^M : \Sp_M \to \Sp$, and it therefore has a left adjoint, which we denote $\iota$. By definition we have that $$ K_{M}(X):= \left(\iHom_{\Sp_M}(\Sigma^{\infty}_U X_+, bu_{M})\right)^{M} $$ where $bu_{M}$ is the $M$-spectrum representing equivariant K-theory, $\Sigma^\infty_U$ is the stabilization functor from pointed $M$-spaces to $M$-spectra, and $\iHom_{\Sp_M}$ is the internal function spectrum in the symmetric monoidal category of $M$-spectra \cite{lewis1986equivariant}*{page 72}.Thus by the (spectrally enhanced) adjunction and the definition of inner Hom in a symmetric monoidal category we have
$$\Hom_{\Sp}(\Sigma^\infty Y_+,K_M(X)) \simeq \left( \iHom_{\Sp_M}(\iota(\Sigma^\infty Y_+) \wedge \Sigma^\infty_U(X_+),bu_M) \right)^M$$ The claim now follows from the natural isomorphism $\iota(\Sigma^\infty Y_+) \simeq \Sigma^\infty_U (Y_+)$ \cite{lewis1986equivariant}*{Remark II.3.14(i)}, where $Y$ is regarded as an $M$-space with trivial $M$ action, the fact that $\Sigma^\infty_U$ maps smash products of pointed $M$-spaces to smash products of $M$-spectra \cite{lewis1986equivariant}*{Remark II.3.14(iii)}, and the fact that $Y_+ \wedge X_+ \simeq (Y \times X)_+$ for pointed $M$-spaces.
\end{proof}

Now that we have constructed a comparison map of presheaves $K^{top}(\Perf(X/G)) \to K_M(X^{an})$, we construct a comparison map for Borel-Moore homology $\rho_{G,X} : K^{top}(\DCoh(X/G)) \to K_M^{c,\dual}(X^{an})$ for a $G$-quasi-projective scheme $X$ as follows: Choose a $G$-equivariant closed immersion $i : X \hookrightarrow Z$ into a smooth $G$-scheme $Z$ with open complement $j : U \subset X$. Then by \autoref{lem:torsor_pullback} and \autoref{lem:torsor_pullback_segal}, and using the canonical equivalence $\Perf(Z/G) \simeq \DCoh(Z/G)$, there is a dotted arrow completing the map of exact triangles
\begin{equation} \label{eqn:define_comparison}
\xymatrix{K^{top}(\DCoh(X/G)) \ar[r]^{i_\ast} \ar@{-->}[d]^{\rho_{G,X}} & K^{top}(\Perf(Z/G)) \ar[r]^{j^\ast} \ar[d]^{\rho_{G,Z}^{perf}} & K^{top}(\Perf(U/G)) \ar[d] \ar[r] \ar[d]^{\rho_{G,U}^{perf}} & \\ 
K_M^{c,\dual}(X^{an}) \ar[r]^{i_\ast} & K_M^{c,\dual}(Z^{an}) \ar[r]^{j^\ast} & K_M^{c,\dual}(U^{an}) \ar[r] & }
\end{equation}

\begin{lem} \label{lem:comparison}
The map $\rho_{G,X} : K^{top}(\DCoh(X/G)) \to K_M^{c,\dual}(X^{an})$ is canonically defined up to homotopy, independent of the choice of embedding $X \hookrightarrow Z$, and agrees with the comparison map $\rho_{G,X}^{perf}$ if $X$ is smooth. $\rho_{G,X}$ commutes (up to homotopy) with proper pushforwards, restriction to open $G$-subschemes, the canonical equivalence between the Borel-Moore homology of $X$ and that of a torsor for a vector bundle over $X$, and restriction from $G$-equivariant $K$-theory to $H$-equivariant $K$-theory for a reductive subgroup $H \subset G$ such that $H$ is the complexification of $H \cap M$.
\end{lem}
\begin{proof}

The fact that the map $\rho_{G,X}$ is canonical up to homotopy follows from the fact that the rest of the diagram is strictly commutative and the formation of homotopy fibers is functorial in the category of spectra. So, the map between homotopy fibers of the two restriction maps $j^\ast$ is canonical, and the d{\'e}vissage isomorphism with the homotopy fiber of $j^\ast$ induced by the Gysin map $i_\ast$ is canonical for both the theories $K^{top}(\DCoh(-/G))$ and $K_M^{c,\dual}((-)^{an})$. The fact that $\rho_{G,X} \sim \rho_{G,X}^{perf}$ when $X$ is smooth follows from considering the identity embedding $X \hookrightarrow X$. Independence of the choice of embedding is essentially proved in \cite{thomason1988equivariant}, which is an extension to the equivariant setting of \cite{thomason1983riemann}. For the benefit of the reader, we explain the conceptual core of argument:

Define a category $\op{Emb}$ of ``virtual embeddings'' whose objects are $G$-quasi-projective schemes and whose morphisms $X \rightsquigarrow Z$ consist of a $G$-equivariant closed subscheme $V \hookrightarrow Z$ along with a $G$-equivariant map $V \to X$ which can be factored as a composition of maps which are torsors for locally free sheaves. Composition is given by pullback of closed subschemes. Then in the proofs of \cites{thomason1988equivariant,thomason1983riemann}, Thomason shows that given two maps $X \rightsquigarrow Z_1$ and $X \rightsquigarrow Z_2$, there is a linear action of $G$ on $\bA^n$ and maps $Z_1 \rightsquigarrow \bA^n$ and $Z_2 \rightsquigarrow \bA^n$ such that the two compositions $X \rightsquigarrow \bA^n$ agree.\footnote{More precisely, the proof of \cite{thomason1988equivariant}*{Proposition 5.8} shows that for any $G$-quasi-projective $X$ there is a map $X \rightsquigarrow \bA^n$ for some linear represenation of $G$. Thus it suffices to consider the case of two maps $X \rightsquigarrow \bA^{n_i}$, $i=1,2$. Next if $V \to X$ is a composition of torsors for locally free sheaves and $V \hookrightarrow \bA^{n_i}$, $i=1,2$ are two $G$-equivariant closed embeddings, then the proof of \cite{thomason1983riemann}*{Lemma 4.2} works equivariantly to constuct an equivariant embeddings $\bA^{n_i} \hookrightarrow \bA^{n_1} \times \bA^{n_2}$ such that the two induced embeddings $V \hookrightarrow \bA^{n_1} \times \bA^{n_2}$ agree. Thus it suffices to show that for any two maps $X \rightsquigarrow \bA^{n_i}$ corresponding to two fibrations $V_i \to X$, one can compose with maps $\bA^{n_i} \rightsquigarrow \bA^{n'_i}$ such that if $V_i' \subset \bA^{n_i'}$ corresponds to the compositions $X \rightsquigarrow \bA^{n_i'}$, then $V_1' \simeq V_2'$ over $X$. This follows from the proof of \cite{thomason1983riemann}*{Proposition 4.7}, which also works equivariantly.} In particular, the under-category $\op{Emb}_{X/}$ is filtered, as is the category $\op{Emb}_{X/}^{sm}$ of virtual embeddings $X \rightsquigarrow Z$, where $Z$ is smooth.

Now let $X \rightsquigarrow Z$ be a virtual embedding into a smooth $Z$, which corresponds to $(V \hookrightarrow Z, \pi : V \to X)$. For ease of notation let $E^1(X) = K^{top}(\Perf(X/G))$ and let $E^2(X)=K_M(X^{an})$ for any smooth $G$-quasi-projective scheme $X$. Then we define
\[
E^i(X \rightsquigarrow Z) := E^i_{V}(X) = \op{fib}(E^i(Z) \to E^i(Z-V)), i=1,2.
\]
Note that by functoriality of the homotopy fiber the comparison maps on cohomology theories $\rho_{G,Z}^{perf}$ and $\rho_{G,Z-V}^{perf}$ induce a map
\[
\rho(X \rightsquigarrow Z) : E^1(X \rightsquigarrow Z) \to E^2(X \rightsquigarrow Z).
\]
We claim that the assignment $E^i(X \rightsquigarrow Z)$ can be extended to a functor $E^i : \op{Emb}_{X/}^{sm} \to \op{Ho}(\Sp)$ in which all arrows in $\op{Emb}_{X/}^{sm}$ map to isomorphisms, and the comparison map $\rho(X \rightsquigarrow Z)$ defines a natural transformation of functors $E^1 \to E^2$. Indeed given a composition $X \rightsquigarrow Z_1 \to Z_2$, which we regard as a map $Z_1 \rightsquigarrow Z_2$ in $\op{Emb}_{X/}^{sm}$, consider the following diagram, where the central square is cartesian:
\[
\xymatrix{ & & V'_2 \ar[dl] \ar@{^{(}->}[dr] \ar@/^15pt/[ddrr]^i \ar@/_15pt/[ddll]_\pi & & \\ & V_1 \ar@{^{(}->}[dr] \ar[dl] & & V_2 \ar@{^{(}->}[dr] \ar[dl] & \\ X \ar@{~>}[rr] & & Z_1 \ar@{~>}[rr] & & Z_2 }.
\]
Then we have canonical isomorphisms $E^i_{V_1}(Z_1) \xrightarrow{\pi^\ast} E^i_{V'_2}(V_2) \xrightarrow{i_\ast} E^i_{V_2'}(Z_2)$ for $i=1,2$, where the construction of the latter map uses Poincar{\'e} duality (all of the $G$-schemes involved are smooth) to convert cohomology to homology with locally compact supports (i.e. Borel-Moore homology), and $i_\ast$ is an isomorphism by d{\'e}vissage for the Borel-Moore theory. We leave the somewhat involved diagram chase to the reader to show that the isomorphisms $E^{i}(X \rightsquigarrow Z_1) \to E^i(X \rightsquigarrow Z_2)$ commute with the comparison maps $\rho(X \rightsquigarrow Z_1)$ and $\rho(X \rightsquigarrow Z_2)$ up to homotopy.

Our comparison map $\rho_{G,X}$ arises from choosing an embedding $X \rightsquigarrow Z$ where $V=X$ and combining $\rho(X \rightsquigarrow Z)$ with the canonical d{\'e}vissage isomorphisms $K^{top}(\DCoh(X/G)) \simeq E^1(X \rightsquigarrow Z)$ and $K_M^{c,\dual}(X^{an}) \simeq E^2(X \rightsquigarrow Z)$. We observe that using part (5) of \autoref{lem:torsor_pullback} and part (4) of \autoref{lem:torsor_pullback_segal}, for any $X \rightsquigarrow Z$ we can compose the canonical equivalence $K^{top}(\DCoh(X/G)) \to K^{top}(\DCoh(V/G))$ followed with the canonical d{\'e}vissage equivalence $K^{top}(\DCoh(V/G)) \to E^1_{V}(Z)$, and likewise for $K_M^{c,\dual}(-)$. Given a map $V \to Z$ which is a composition of torsors for vector bundles, and given a closed $G$-subscheme $X \hookrightarrow Z$, the canonical equivalences $K^{top}(\DCoh(Z/G)) \to K^{top}(\DCoh(V/G))$ and $K^{top}(\DCoh(X/G)) \to K^{top}(\DCoh(V|_X / G))$ commute with the pushforward maps, and likewise for $K_M^{c,\dual}(-)$. It follows that for any composition $X \rightsquigarrow Z_1 \rightsquigarrow Z_2$, the composition of the canonical maps $K^{top}(\DCoh(X/G)) \to E^1(X \rightsquigarrow Z_1) \to E^1(X \rightsquigarrow Z_2)$ is homotopic to the canonical map $K^{top}(\DCoh(X/G)) \to E^1(X \rightsquigarrow Z_2)$, and the same is true for $K_M^{c,dual}(-)$. Combining these facts, we see that the comparison map $\rho_{G,X}$ defined via the canonical equivalences for any smooth virtual embedding $X \rightsquigarrow Z$
\[
K^{top}(\DCoh(X/G)) \simeq E^1(X \rightsquigarrow Z) \xrightarrow{\rho(X \rightsquigarrow Z)} E^2(X\rightsquigarrow Z) \simeq K_M^{c,\dual}(X^{an})
\]
is independed of the choice of $X \rightsquigarrow Z$, because the category $\op{Emb}_{X/}^{sm}$ is filtered.

Now that we have shown that the comparison map $\rho_{G,X} : K^{top}(\DCoh(X/G)) \to K_M^{c,\dual}(X^{an})$ is independent, up to homotopy, of the choice of embedding $X \rightsquigarrow Z$, it is fairly straightforward to show that $\rho_{G,X}$ commutes with restriction to an open subset and restriction of requivariance from $G$ to a reductive subgroup $H \subset G$ which is the complexification of $H \cap M$, because the corresponding claims hold for the comparison map $\rho_{G,X}^{perf}$ for smooth $G$-schemes. The fact that $\rho_{G,X}$ commutes with the equivalences $K^{top}(\DCoh(X/G)) \simeq K^{top}(\DCoh(V/G))$ and $K_M^{c,\dual}(X^{an}) \simeq K_M^{c,\dual}(V^{an})$ for a map $V \to X$ which is a torsor for a vector bundle follows from the more general description of $\rho_{G,X}$ above, because any closed embedding into a smooth $G$-scheme $V \hookrightarrow Z$ defines both a virtual embedding $X \rightsquigarrow Z$ and $V \rightsquigarrow Z$, and the canonical equivalence $K^{top}(\DCoh(X/G)) \simeq E^1(X \rightsquigarrow Z)$ is the composition of the equivalences $K^{top}(\DCoh(X/G)) \simeq K^{top}(\DCoh(V/G)) \simeq E^1(V\rightsquigarrow Z)$ and likewise for $K_M^{c,\dual}(-)$.

Finally, one shows that $\rho_{G,X}$ commutes pushforward along a proper map $p : X \to Y$ by choosing a closed immersion $Y \hookrightarrow Z$ into a smooth $G$-scheme $Z$ and a projective bundle $\bP(\cV)$ over $Z$ such that $p$ factors through a closed immersion $X \hookrightarrow Y$. One shows that $\rho$ commutes with pushforward along the projection $\bP(\cV|_Y) \to Y$ using d{\'e}vissage and the corresponding fact for $\rho_{G,Z}^{perf}$ for smooth $G$-schemes. It therefore suffices to show that $\rho$ commutes with closed immersions. For any closed immersion of $G$-schemes $X_1 \hookrightarrow X_2$, we can choose a closed immersion into a smooth $G$-scheme $X_2 \hookrightarrow Z$. The fact that $\rho$ commutes with pushforward along $X_1 \hookrightarrow X_2$ follows from the fact that $i_\ast$ is part of a canonical map of exact triangles from $E^i(X_1) \to E^i(Z) \to E^i(Z-X_1) \to$ to $E^i(X_2) \to E^i(Z) \to E^i(Z-X_2)$ for both theories $i=1,2$.

\end{proof}

Note that in \autoref{lem:comparison} we have passed from strict presheaves on the category of algebraic $G$-spaces with values in spectra, to presheaves on the category of \emph{quasi-projective} $G$-schemes with values in the \emph{homotopy category} $\op{Ho}(\Sp)$. Presumably neither of these relaxations are necessary, but they simplify our discussion.

\begin{thm} \label{thm:topological_comparison}
The natural map $\rho_{G,X} : K^{top}(D^b \Coh(X/G)) \xrightarrow{\simeq} K^{c,\dual}_{M}(X^{an})$ construced in \eqref{eqn:define_comparison} and \autoref{lem:comparison} is a weak equivalence.
\end{thm}

\begin{proof} [Proof of \autoref{thm:topological_comparison}]

By the commutativity of the diagram \eqref{eqn:define_comparison} and the fact that $\rho_{G,X}^{perf} \sim \rho_{G,X}$ when $X$ is smooth, it suffices to show that $\rho_{G,X}$ is a weak equivalence when $X$ is smooth $G$-quasi-projective scheme. Let $T \subset G$ be a maximal torus which is the complexification of a compact maximal torus $T_c = T\cap M$, and let $B \subset G$ be a Borel subgroup containing $T$. For both $G$-equivariant cohomology theories, pullback along the map $G\times_B X \to X$ is an injection which is canonically split by the pushforward map. The map which forgets from $G$-equivariance to $T$-equivariance, then restricts along the $T$-equivariant map $\{1\} \times X \to G \times_B X$ induces an equivalences $K_M(G \times_B X) \to K_{T_c}(X)$ because topologically $G \times_B X = M \times_{T_c} X$ and the same maps induces an equivalence $K^{top}(\Perf(G \times_B X/G)) \to K^{top}(\Perf(X/T))$ because $X/T \to G\times_B X / G \simeq X/B$ is a composition of torsors for line bundles on $X/B$. The comparison map $\rho_{G,X}^{perf}$ commutes with the operations of pullback, proper pushforward, and restriction to subgroups, so $\rho$ is compatible with the splitting of the inclusions $K^{top}(\Perf(X/G)) \hookrightarrow K^{top}(\Perf(X/T))$ and $K_{M}(X^{an}) \hookrightarrow K_{T_c}(X^{an})$, so it suffices to prove the claim when $G=T$ is a torus.

We can stratify $X/T$ by smooth $T$-schemes of the form $U \times (T/T')$, where $T' \subset T$ is an algebraic subgroup and $T$ acts trivially on $U$. Using the localization sequences for closed immersions in \autoref{lem:torsor_pullback} and \autoref{lem:torsor_pullback_segal} and the compatibility of $\rho_{G,X}$ with pushforward along closed immersions and restriction to open subsets, \autoref{lem:comparison}, it suffices to prove the claim for schemes of this form. The fact that $\Perf(U \times (T/T') / T) \simeq \Perf(U \times BT') \simeq \bigoplus_{\chi} \Perf(U)$, where $\chi$ ranges over the group of characters of the diagonalizable group $T'$, implies that $K^{top}(\Perf(U \times (T/T') / T)) \simeq \bigoplus_{\chi} K^{top}(\Perf(U))$. There is an analogous decomposition of $K_{T_c}(U \times T/T')$, and $\rho_{T,U \times T/T'}$ respects this direct sum decomposition because the summands are the essential image of pullback along the map $U \times T/T' \to U$ followed by tensoring with the various characters of $T'$. We note that when the group is trivial, our comparison map agrees with the one constructed in \cite{blanc2012topological}*{Proposition 4.32}, therefore it is an equivalence, and the claim follows.
\end{proof}

\begin{rem}
If $G$ is not necessarily reductive, then one can choose a decomposition $G = U \rtimes H$, where $H$ is reductive and $U$ is a connected unipotent group. As in the first step in the proof of \cite{thomason1988equivariant}*{Theorem 5.9}, one shows that the map of stacks $X/H \to X/G$ can be factored as a sequence of torsors for vector bundles, so the canonical restriction map $K^{top}(\Perf(X/G)) \to K^{top}(\Perf(X/H))$ is an equivalence by \autoref{lem:torsor_pullback}. Combining this with the previous theorem shows that for a maximal compact subgroup $M \subset H \subset G$, the topologization functor is an equivalence $K^{top}(\Perf(X/G)) \to K_M(X^{an})$ as presheaves of spectra on $\Smooth_G$, and we have a comparison isomorphism $\rho_{G,X} : K^{top}(\DCoh(X/G)) \to K^{c \dual}_M(X^{an})$.
\end{rem}

\subsection{The case of smooth Deligne-Mumford stacks}
\label{sect:DM_invariants}

Here we provide an explicit computation of the periodic cyclic homology of $\Perf(\X)$ for a smooth Deligne-Mumford stack of finite type over $\bC$ and study its noncommutative Hodge theory when it is proper. The results of this section are likely known to experts.

Given a smooth scheme $U$, we can consider its de Rham complex, $0 \to \cO_U \to \Omega^1_U \to \cdots$, a complex of vector spaces. We can regard this as a $\Lambda$-module $\Omega_\bullet(U)$ by defining $\Omega_p(U):= \Omega^{-p}_U$ and letting $B$ act via the de Rham differential. Even though the $\Lambda$-module structure is not $\cO_U$-linear, it still defines a sheaf of $\Lambda$-modules on the small site $\X_{et}$ for any smooth DM stack $\X$. We define the de Rham cohomology of a smooth Deligne-Mumford stack $\X$ to be the $\Lambda$-module
$$H_{dR}(\X) := \RGamma(\X_{et}, \Omega_\bullet).$$

There are several other ways to present $H_{dR}(\X)$. First note that we can equivalently restrict to the sub-site of \'{e}tale maps $U \to \X$ for which $U$ is affine, which we denote $\X^{\tinyop{aff}}_{et}$, because it has an equivalent topos of sheaves, i.e. the canonical map is an equivalence
$$H_{dR}(\X) \xrightarrow{\simeq} \RGamma(\X^{\tinyop{aff}}_{et},\Omega_\bullet).$$
We can consider the sheaf of $\Lambda$-modules on $\X^{\tinyop{aff}}_{et}$ given by $U \mapsto \hoch(\cO_U)$, the Hochschild complex of coordinate algebra. This admits a canonical map to the presheaf of $\Lambda$-modules given by $U \mapsto \hoch(\Perf(U))$. Likewise, for any smooth affine scheme the map $\operatorname{HKR}$ is a map of $\Lambda$-modules $\hoch(\cO_U) \to \Omega_\bullet(U)$ and compatible with \'{e}tale base change, so they induce maps of presheaves on $\X_{et}^{\tinyop{aff}}$
\begin{lem} \label{lem:local_HKR}
The canonical maps
\[
\RGamma(\X_{et}^{\tinyop{aff}},\Omega_\bullet) \gets \RGamma(\X_{et}^{\tinyop{aff}},\hoch(\cO_{-})) \to \RGamma(\X_{et}^{\tinyop{aff}},\hoch(\Perf(-)))
\]
are all equivalences of $\Lambda$-modules.
\end{lem}
\begin{proof}
These maps are all equivalences for affine $U$ at the level of underlying complexes. The result follows formally from the fact that a map of $\Lambda$-modules is an equivalence if and only if the underlying map of complexes is an equivalence, and the forgetful functor taking a $\Lambda$-module to its underlying complex commutes with limits, hence commutes with $\RGamma$.
\end{proof}

The following is due to Toen, and essentially follows the argument of \cite{toen1999theoremes} in the case of algebraic $K$-theory. We will need to use both the derived inertia stack $I_\X$ and its underlying classical stack $I^{cl}_\X \subset I_\X$.

\begin{prop}[Toen, unpublished] \label{prop:Toen_de_rham}
Let $\X$ be a smooth Deligne-Mumford stack, and let $I^{cl}_\X$ denote its classical inertia stack. There is a natural isomorphism of $\Lambda$-modules $\hoch (\Perf(\X)) \to H_{dR}(I^{cl}_\X)$.
\end{prop}

The idea of the proof is to show that the formation of both complexes is local in the etale topology over the coarse moduli space of $\X$, so one can reduce to the case of a global quotient. Thus a key observation is that the formation of the derived inertia stack $I_\X$ is \'{e}tale local.

\begin{lem}
Let $\X \to X$ be a map from a stack to a separated algebraic space, and let $\U \to U$ be the base change along an etale map $U \to X$. Then $I_\U \simeq I_\X \times_X U$.
\end{lem}
\begin{proof}
This can be seen for the derived inertia stack from a functor-of-points definition of $I_\U$. We let $\U(T)$ denote the $\infty$-groupoid of maps from $T$ to $\U$ for a derived affine scheme $T$.
\begin{align*}
I_{\U}(T) &= \U(T) \times_{\U(T) \times \U(T)} \U(T) \\
&\sim \Map(S^1,\U(T)) \\
&\sim \Map(S^1,\X(T)) \times_{\Map(S^1,X(T))} \Map(S^1,U(T))
\end{align*}
So in order to show that $I_\U(T) \simeq I_\X(T) \times_{X(T)} U(T)$, it will suffice to show that $I_U \simeq I_X \times_X U$ in the derived sense. Consider the following diagram, in which each square is Cartesian and the vertical arrows are closed immersions:
$$\xymatrix{U \ar[r] \ar@/^10pt/[rr]^\simeq & U \times_X U \ar[r] \ar[d] & \Gamma \ar[r] \ar[d] & X \ar[d] \\ & U \times U \ar[r] & U \times X \ar[r] & X \times X }$$
Here $\Gamma$ denotes the graph of the morphism $U \to X$. Then by definition $I_X$ is the derived self intersection of the closed subspace $X \to X \times X$, so in order to prove the claim it will suffice to show that $\Gamma \times_{U \times X} \Gamma$ is isomorphic to $I_\cU$ as a derived scheme over $U$.

The map $U \to \Gamma$ is an isomorphism on underlying classical algebraic spaces, and it follows from the fact that $U \to U\times_X U$ is an etale closed immersion of closed substacks of $U \times U$ that the induced map $I_U \to \Gamma \times_{U \times X} \Gamma$ induces an isomorphism on cotangent complexes as well, hence it is an isomorphism.
\end{proof}

\begin{proof} [Proof of \autoref{prop:Toen_de_rham}]
The pullback functor along the projection $I^{cl}_\X \to \X$ induces a map $\hoch(\Perf(\X)) \to \hoch(\Perf(I^{cl}_\X))$. For any \'{e}tale $U / I^{cl}_\X$, the pullback functor induces a natural map $\hoch(\Perf(I^{cl}_\X)) \to \hoch(\Perf(U))$. Thus we get a map of presheaves of $\Lambda$-modules
$$\hoch(\Perf(\X)) \to \RGamma((I^{cl}_\X)_{et}^{\tinyop{aff}},\hoch(\Perf(-))) \simeq H_{dR}(I^{cl}_\X).$$
Note that if $p : \X \to X$ is the coarse moduli space of $\X$, then the map constructed above is functorial with respect to pullback along maps $U \to X$.

We claim that $\hoch(\Perf(\X))$, regarded as a presheaf over $X$, has \'{e}tale descent. Indeed, consider any \'{e}tale map $U \to X$, and let $\U = \X \times_X U$. Because the derived category of $\U$ is compactly generated \cite{hallrydh2014DM}, we can identify
$$\hoch(\Perf(\U)) \simeq \RGamma(\U,\Delta^\ast \Delta_\ast (\cO_\U)) \simeq \RGamma(U,(p_U)_\ast \cO_{I_\U}),$$
where $\Delta : \U \to \U \times \U$ is the diagonal, $p_U : \U \to U$ is the base change of $p$, and $\cO_{I_\U}$ is the structure sheaf of the derived inertia stack, regarded as a finite algebra over $\cO_\U$. In the previous lemma, we saw that the formation of $I_\U$ commutes with \'{e}tale base change, so this combined with the projection formula implies that $\RGamma(U,p_\ast \cO_{I_\U}) \simeq \RGamma(U,p_\ast(\cO_{I_\X})|_U)$, functorially in $U$. The presheaf $U/X \mapsto \RGamma(U,p_\ast(\cO_{I_\X})|_U)$ has \'{e}tale descent, so $U \mapsto \hoch(\Perf(\U))$ does as well.

Thus in order to show that $\hoch(\Perf(\X)) \to H_{dR}(I^{cl}_\X)$ is an equivalence, it suffices to verify this after base change to an \'{e}tale cover of $X$. We can find such a $U \to X$ such that $\U = \X \times_X U$ is a global quotient of a scheme by a finite group action. In that case, the result is shown in \cite{baranovsky2003orbifold}*{Proposition 4}.
\end{proof}

Finally after applying the Tate construction, i.e. passing to periodic cyclic homology, we can compare this to the cohomology of $|\X^{an}|$, the geometric realization of the underlying topological stack (in the analytic topology) associated to $\X$ \cite{noohi2012homotopy}, as well as the cohomology of a coarse moduli space $\X \to X$.

\begin{lem} \label{lem:cdim}
Let $\X$ be a Noetherian separated DM stack of finite type over a Noetherian base scheme. Assume that $\X$ has finite dimension. Then $\X$ has finite \'{e}tale cohomological dimension with $\bQ$-linear coefficients, and the functor $\RGamma(\X_{et},-)$ commutes with filtered colimits.
\end{lem}

\begin{proof}
We first claim that the pushforward along the projection to the coarse moduli space $p: \X \to X$ is exact. Indeed this can be checked \'{e}tale locally on $X$, and so we may assume that $\X$ is a global quotient $U /G$, where $G$ is a finite group. One can factor $p$ as $U / G \to X \times BG \to X$ -- pushforward along the first is exact by \cite{stacksProject}*{Tag 03QP}, and the second is exact because we are using characteristic $0$ coefficients.

It now suffices to prove the claim when $\X = X$ is a Noetherian separated algebraic space of finite type over a Noetherian base scheme. In this case, we can apply the induction principle of \cite{stacksProject}*{Tag 08GP} and the fact that \'{e}tale cohomology takes elementary excision squares to homotopy cartesian squares to reduce to the case of affine schemes. In this case, the result follows from the fact that derived global sections of characteristic $0$ sheaves on a Noetherian scheme in the \'{e}tale topology argees with that in the Nisnevich topology, and the Nisnevich topology has cohomological dimension $\leq d$.

Finally, the implication that finite cohomological dimension implies commutation with filtered colimits in the unbounded derived category is \cite{cisinski2013etale}*{Lemma 1.1.7}.

\end{proof}

\begin{lem} \label{lem:deRham}
There are natural isomorphisms $$H_{dR}(\X)^{\Tate} \simeq C^\ast_{sing}(|\X^{an}|;\bQ)\otimes_\bQ \bC\ls{u} \simeq C^\ast_{sing}(X;\bQ) \otimes_\bQ \bC\ls{u}$$
\end{lem}
\begin{proof}
The de Rham isomorphism gives a canonical isomorphism of pre-sheaves of $\bC\ls{u}$-modules on $\X^{\tinyop{aff}}_{et}$ between $U \mapsto (\Omega_{\bullet}(U))^{\Tate}$ and $U \mapsto C^\ast_{sing}(U^{an};\bC)\ls{u}$, so we have a canonical isomorphism\footnote{All of the singular complexes we will encounter have finite dimensional total cohomology, so $M\ls{u} \simeq M \otimes_\bC \bC\ls{u}$.}
\[
C_{sing}^\ast (|\Y^{an}|;\bC) \otimes \bC\ls{u} \simeq \RGamma(\X_{et}^{\tinyop{aff}},\Omega_\bullet(-)^{\Tate}).
\]
It therefore suffices to show that the Tate construction commutes with taking derived global sections for the sheaf of $\Lambda$-modules $\Omega_\bullet$. For this we observe that the functor $M \mapsto M^{\Sone}$ commutes with homotopy limits, and hence with derived global sections, and $M^{\Tate}$ is the filtered colimit of $M^{\Sone} \to M^{\Sone}[2] \to M^{\Sone}[4] \to \cdots$, so its formation commutes with $\RGamma$ by the previous lemma.

Finally, one can check that the pullback map $C^\ast_{sing}(Y^{an};\bQ) \to C^\ast_{sing}(|\Y^{an}|;\bQ)$ is an equivalence locally in the analytic topology on $Y^{an}$. Locally $\Y^{an}$ is isomorphic to a global quotient of a scheme by a finite group, for which the fact is well-known.
\end{proof}

\subsection{Equivariant \texorpdfstring{$K$}{K}-theory and periodic cyclic homology}

For a dg-category, $\cC$, it is natural to ask if the Chern character induces an equivalence $K^{top}(\cC) \otimes \bC \to \hoch^{per}(\cC)$. This is referred to as the lattice conjecture in \cite{blanc2012topological}, where it is conjectured to hold for all smooth and proper dg-categories. Here we observe some situations in which the lattice conjecture holds, even for categories which are not smooth and proper.

\begin{thm}[Lattice conjecture for smooth quotient stacks] \label{thm:lattice_conjecture}
Let $G$ be an algebraic group acting on a smooth quasi-projective scheme $X$. If $X/G$ admits a semi-complete KN stratification (\autoref{def:semicomplete_KN}), then the Chern character induces an equivalence $K^{top}(\Perf(X/G)) \otimes \bC \to \hoch^{per}(\Perf(X/G))$.
\end{thm}

\begin{lem} \label{lem:localization_DM}
Let $\X$ be a smooth Deligne-Mumford stack, and let $i : \Z \hookrightarrow \X$ be a smooth closed substack. Then the pushforward functor fits into a fiber sequence
$$\hoch^{per}(\Perf(\Z)) \xrightarrow{i_\ast} \hoch^{per}(\Perf(\X)) \xrightarrow{j^\ast} \hoch^{per}(\Perf(\X-\Z)).$$
\end{lem}

\begin{proof}
This follows from \autoref{prop:Toen_de_rham}, combined with \autoref{lem:deRham} and the usual Gysin sequence for the regular embedding of inertia stacks $I_\Z \hookrightarrow I_\X$.
\end{proof}

\begin{proof} [Proof of \autoref{thm:lattice_conjecture}]

Because $K^{top}(-) \otimes \bC$ and $\hoch^{per}(-)$ are both additive invariants, proving that the natural transformation
$$K^{top}(\Perf(X/G)) \otimes \bC \to \hoch^{per}(\Perf(X/G))$$
is an equivalence for smooth projective-over-affine $X$ and reductive $G$ reduces to the case where $X/G$ is Deligne-Mumford by \autoref{thm:motivic_main}.

Note that the only point in the proof of \autoref{lem:torsor_pullback} which does not immediately apply to an arbitrary additive invariant is the localization sequence for a closed immersion. Therefore \autoref{lem:localization_DM} implies that \autoref{lem:torsor_pullback} applies to the presheaf $\hoch^{per}(\Perf(-))$, because the only stacks that appear in the proof are DM.

We can now imitate the proof of \autoref{thm:topological_comparison}: $\Perf(X/G)$ is a retract of $\Perf(X/B)$, and $\Perf(X/B) \to \Perf(X/T)$ induces an equivalence for both invariants $K^{top}(\bullet)$ and $\hoch^{per}(\bullet)$, by \autoref{lem:torsor_pullback}. Thus it suffices to consider smooth DM stacks of the form $X/T$. Any such stack admits a stratification by smooth stacks of the form $U \times B\Gamma$ for some finite group $\Gamma$, and by \autoref{lem:localization_DM} it suffices to prove the theorem for such stacks. Thus $\op{Ch} \otimes \bC$ is an equivalence because it is an equivalence for smooth schemes and $\Perf(U \times B\Gamma) \simeq \bigoplus_\chi \Perf(U)$, the sum ranging over characters of $\Gamma$.
\end{proof}

\begin{cor}[Lattice conjecture for smooth DM stacks] \label{cor:lattice_conjecture_DM}
Let $\cX$ be a smooth Deligne-Mumford stack. Then the Chern character induces an equivalence $K^{top}(\Perf(\cX)) \otimes \bC \to \hoch^{per}(\Perf(\cX))$.
\end{cor}
\begin{proof}
We have established a localization sequence for closed immersions of smooth DM stacks for $K^{top}(\Perf(-))$ in \autoref{lem:torsor_pullback}. We do not know if the localizing invariant of dg-categories $HP(-)$ satisfies d{\'e}vissage in the sense of \autoref{lem:torsor_pullback} for arbitrary closed immersions of stacks, but the localization sequence for closed immersions of smooth DM stacks follows from the comparison results \autoref{prop:Toen_de_rham} and \autoref{lem:deRham} and the corresponding fact for singular cohomology.

We therefore have a map of fiber sequences for any closed immersion of smooth DM stacks over $\bC$,
$$\xymatrix@R=10pt{K^{top}(\DCoh(\cZ)) \otimes \bC \ar[r] \ar[d]^{\op{Ch}} & K^{top}(\DCoh(\cX)) \otimes \bC \ar[d]^{\op{Ch}} \ar[r] & K^{top}(\DCoh(\cU)) \otimes \bC \ar[d]^{\op{Ch}} \\ HP(\DCoh(\cZ)) \ar[r] & HP(\DCoh(\cX)) \ar[r] & HP(\DCoh(\cU)) }.$$
From \cite{laumon2000champs}*{Corollaire 6.1.1} every smooth DM stack of finte type admits a stratification by locally closed substacks which are quotients of a smooth affine scheme by a finite group. The corollary follows by applying \autoref{thm:lattice_conjecture} and the fiber sequence above inductively to this stratification. 
\end{proof}

\subsection{Hodge structure on equivariant \texorpdfstring{$K$}{K}-theory}

We can now prove the final result of this paper, the construction of a pure Hodge structure on the equivariant $K$-theory. What we mean by a pure Hodge structure on a spectrum $E$ in this case is simply a Hodge structure on the homotopy groups of that spectrum $\pi_\ast(E)$, i.e. for each $n$ a weight $n$ Hodge structure on $\pi_n(E)$ is a descending filtration of $\pi_n(E) \otimes \bC$ such that
$$\pi_n(E) \otimes \bC = F^p \pi_n(E)\otimes \bC \oplus \overline{F^{n+1-p} \pi_\ast(E) \otimes \bC}, \forall p$$

\begin{thm} \label{thm:Hodge_main}
Let $X$ be a smooth quasi-projective $\bC$-scheme, let $M$ be a compact Lie group whose complexification $G$ acts on $X$. Then if $X/G$ admits a complete KN stratification, the Chern character isomorphism
\[
K_{M}(X^{an}) \otimes \bC \to \hoch^{per}(\Perf(X/G))
\]
combined with the noncommutative Hodge-de Rham sequence induces a pure Hodge structure of weight $n$ on $K^n_{M}(X^{an})$ with a canonical isomorphism
\[
\op{gr}^p_{Hodge} (K^n_{M}(X^{an}) \otimes \bC) \simeq H^{n-2 p} \RGamma(I_\X,\cO_{I_\X}),
\]
where $I_\X$ denotes the derived inertia stack of $\X := X/G$. The Hodge filtration on $K^n_M(X^{an})$ is compatible with pullback maps, and in particular it is a filtration of $\Rep(M)$-modules.
\end{thm}

\begin{rem}
As we will see in the proof, this claim also holds for arbitrary smooth and proper DM stacks over $\bC$, without requiring that $\X$ is a global quotient.
\end{rem}

\begin{proof}
The degeneration property follows from \autoref{cor:degeneration_complete}, and we have a Chern character isomorphism from \autoref{thm:topological_comparison} combined with \autoref{thm:lattice_conjecture}, so all we have to do is check that the filtration on $\hoch^{per}(\Perf(X/G))$ coming from the HdR spectral sequence combined with the rational structure coming from the Chern character defines a weight $n$ pure Hodge structure on $\pi_{-n}(K^{top}(\Perf(X/G)) \otimes \bQ$. This claim is closed under arbitrary direct sums and summands in $\Madd_k$, so by \autoref{thm:motivic_main} it suffices to prove this claim for smooth and proper DM stacks which are global quotients of a $G$-quasi-projective scheme by a reductive group $G$.

For a smooth and proper DM stack, \autoref{prop:Toen_de_rham} gives an isomorphism $H_{dR}(I^{cl}_\X) \simeq \hoch(\Perf(\X))$ of $\Lambda$-modules. Note that $I^{cl}_\X$ is itself a smooth and proper DM stack, and for any smooth DM stack $\Y$ the complex $H_{dR}(\Y)^{\Tate}$ is canonically equivalent to the usual de Rham complex of \cite{satriano2012derham} tensored with $\bC\ls{u}$,
$$\RGamma\left(\Y, \left[0 \to \cO_\Y \to \Omega_\Y^1 \to \cdots \right] \right) \otimes_\bC \bC\ls{u}.$$
However, the usual Hodge filtration differs slightly from the noncommutative one. We have a canonical isomorphism
\begin{align*}
H^n (H_{dR}(\Y)^{\Tate}) &\simeq \bigoplus_{l \equiv n \mod 2} H^l(\Y;\bC) \\
\op{gr}^p F^\bullet_{nc} H^n(H_{dR}(\Y)^{\Tate}) &\simeq \bigoplus_i \RGamma(\Y,\Omega_\Y^i[i-2p]))
\end{align*}
Because the cyclic complex $H_{dR}(\Y)$ has the degeneration property, we may commute taking cohomology $H^n$ and taking associated graded $\op{gr}^p$, so we have
$$\op{gr}^p F^\bullet_{nc} H^n(H_{dR}(\Y)^{\Tate}) \simeq \bigoplus_i H^{n+i-2p}(\Y,\Omega^i_\Y)$$

Thus on each direct summand $H^l(\Y;\bC)$ of $H^n(H_{dR}(\Y)^{\Tate})$, the subquotient $H^{l-p'}(\Y,\Omega_\Y^{p'})$ shows up in $F^p_{nc}$ if and only if $l-p'  = n+p' -2 p''$ for some $p'' \geq p$. In other words the subquotients appearing are those for which $p' \geq p + \frac{l-n}{2}$. It follows that under the direct sum decomposition above we have
$$F^p_{nc} H^n(H_{dR}(\Y)^{\Tate}) \simeq \bigoplus_{l \equiv n \mod 2} F^{p+ \frac{l-n}{2}}_{classical} H^l(\Y;\bC)$$
Thus under the isomorphism $H^n(H_{dR}(\Y)^{\Tate}) \simeq \bigoplus_{l \equiv n \mod 2} H^l(\Y;\bQ) \otimes \bC$ of \autoref{lem:deRham}, the noncommutative Hodge filtration corresponds to the Hodge filtration on $\bigoplus_{l \equiv n \mod 2} H^l(\Y;\bQ)\langle \frac{l-n}{2} \rangle.$ We claim that this rational structure on $H^n(H_{dR}(I^{cl}_\X)^{\Tate})$ agrees with the one induced by the equivalence $K^{top}(\Perf(\X)) \otimes \bC \simeq H_{dR}(I_\X^{cl})^{\Tate}$ of \autoref{thm:lattice_conjecture} and \autoref{lem:deRham}, so that we have an isomorphism of Hodge structures
\begin{equation} \label{eqn:DM_Hodge_str}
\pi_{-n}(K^{top}(\Perf(\X))) \otimes \bQ \simeq \bigoplus_{l \equiv n \mod 2} H^l(I_\X^{cl};\bQ) \langle \frac{l-n}{2} \rangle.
\end{equation}
The Hodge structure on the $l^{th}$ rational cohomology of the de Rham complex of a smooth DM stack has weight $l$ (see \cite{steenbrink1977mixed}), so it would follow that $\pi_{-n}(K^{top}(\Perf(\X))) \otimes \bQ$ has a Hodge structure of weight $n$.

For the claim about the rational structure of $H_{dR}(I_\X^{cl})$, note that the isomorphism $H_{dR}(I^{cl}_\X)^{\Tate} \simeq K^{top}(\Perf(\X)) \otimes \bC$ results from applying the derived global sections functor to isomorphic sheaves on the \'{e}tale site of $I^{cl}_\X$,
\begin{align*}
\RGamma((I^{cl}_\X)^{\tinyop{aff}}_{et}, K^{top}(\Perf(-)) \otimes \bC) &\simeq \RGamma((I^{cl}_\X)^{\tinyop{aff}}_{et}, \hoch^{per}(\Perf(-))) \\
&\simeq \RGamma((I^{cl}_\X)^{\tinyop{aff}}_{et},C^\ast_{sing}((-)^{an};\bC)\ls{u}).
\end{align*}
But according to \cite{blanc2012topological}*{Proposition 4.32}, the noncommutative Chern character for smooth $\bC$-schemes factors through the twisted Chern character under the natural equivalence $\hoch^{per}(\Perf(X)) \simeq H_{dR}(X)^{\Tate} \simeq C^\ast_{sing}(X;\bC\ls{u})$. It follows that the isomorphism above is the complexification of a map of presheaves of $\bQ$-complexes on $(I^{cl}_\X)^{\tinyop{aff}}_{et}$
$$K^{top}(\Perf(-)) \otimes \bQ \to C^\ast_{sing}((-)^{an};\bQ) \otimes_\bQ \bQ\ls{\frac{u}{2\pi i}},$$
which is also a level-wise weak equivalence. Thus the rational structure on $H^n(H_{dR}(I_\X^{cl})^{\Tate})$ agrees with that of the Hodge structure of \autoref{eqn:DM_Hodge_str}.

\end{proof}

\begin{rem}
For any of the quotient stacks appearing in \autoref{amplif:zero_locus}, the theorem above still holds for $\DCoh(\X)$ with the same proof, with the exception of the explicit computation of $\op{gr}^p H^n(C^{per}_\bullet(\DCoh(\X)))$ when $\X$ is not smooth. In particular we have:
\begin{itemize}
\item A canonical isomorphism $K_M(X^{an}) \otimes \bC \to C^{per}_\bullet(\DCoh(X/G))$ which factors through the complexification of the Chern character $K^{top}(\DCoh(X/G)) \to C^{per}_\bullet(\DCoh(X/G))$;
\item The degeneration property for $\DCoh(X/G)$; and
\item A pure Hodge structure of weight $n$ on $\pi_{-n} K^{top}(\DCoh(\X))$ coming from the degeneration of the noncommutative Hodge-de Rham sequence which is $\Rep(G)$-linear.
\end{itemize}
\end{rem}


\section{Hodge-de Rham degeneration for singularity categories}
\label{sect:singularities}

In this section we extend our methods to establish the degeneration property for certain ``dg-categories of singularities'' $\MF(X/G,W)$ associated to an equivariant Landau-Ginzburg model, i.e., a smooth $G$-variety $X$ and a $G$-invariant function $W : X \to \bA^1$. The notation $\MF$ is more frequently used for categories of matrix factorizations, which are equivalent to singularity categories for LG models on smooth schemes \cite{MR2910782}. We have chosen to use $\MF$ to denote singularity categories for consistency with \cite{Preygel}.

The categories $\MF(X/G,W)$ will be $\bZ/2\bZ$-graded, and we will need a suitable $\infty$-categorical model to work with these. The $\infty$-category of (essentially) small idempotent complete stable $\infty$-categories, $\op{Cat}_\infty^{perf}$, admits a symmetric monoidal structure \cite{blumberg2013universal}*{Section 3.1}. For $\cA, \cB \in \op{Cat}_\infty^{perf}$, $\cA \itimes \cB$ is the category of compact objects in $\op{Ind}(\cA) \widehat{\otimes} \op{Ind}(\cB)$, which is idempotent complete.

For any $E_\infty$-algebra $R$, $\Perf(R)$ is canonically a commutative algebra object in $\op{Cat}_{\infty}^{perf}$. We let
\[
\op{LinCat}_R^{sm} = (\Perf(R)^\otimes)\Mod(\op{Cat}_\infty^{perf})
\]
denote the $\infty$-category of $\Perf(R)^\otimes$-module objects. $\op{LinCat}_R^{sm}$ is equivalent, via the ind-completion functor, to the $\infty$-category of $(R\Mod)^\otimes$-module objects in the $\infty$-category of compactly generated presentable stable infinity categories with functors which preserve colimits and compact objects, as in \cite{DAG-VII}*{Definition 6.2}. Because $(R\Mod)^\otimes$ is a commutative algebra object, $\op{LinCat}_R^{sm}$ has a canonical symmetric monoidal structure \cite{lurie2016higher}*{Thm.~4.5.2.1}. In addition, if $R \to R'$ is a map of $E_\infty$-algebras, the tensor product induces a map of commutative algebra objects $\Perf(R)^\otimes \to \Perf(R')^\otimes$, and this induces a symmetric monoidal pullback functor \cite{lurie2016higher}*{Thm.~4.5.3.1} which we denote
\[
(-)\itimes_R R' \colon \op{LinCat}_R^{sm} \to \op{LinCat}_{R'}^{sm}
\]
These constructions work just as well with $\Perf(R)^\otimes$ and $\Perf(R')^\otimes$ replaced by any other essentially small stable idempotent complete symmetric monoidal $\infty$-categories.

Any $\infty$-category $\cC$ in $\op{LinCat}_R^{sm}$ is canonically enriched over $R\Mod$, i.e. $\RHom(E,F) \in R\Mod$ for any $E,F \in \cC$, via an inner-Hom construction. In fact, by \cite{cohn2013differential}, $\op{LinCat}_R^{sm}$ is equivalent to the $\infty$-category of categories enriched in $R$-module spectra. Regarding $k\ls{\beta}$, where $\beta$ is variable of homological degree $-2$, as an $E_\infty$-algebra via the forgetful functor from dg-algebras to $E_\infty$-algebras, this identifies $\op{LinCat}_{k\ls{\beta}}^{sm}$ with the $\infty$-category of categories enriched over dg-$k\ls{\beta}$-modules, or $\bZ/2\bZ$-graded dg-categories. This justifies using $\op{LinCat}_{k\ls{\beta}}^{sm}$ as our model for $\bZ/2\bZ$-graded dg-categories.


The main results of this section, \autoref{thm:motivic_main_lg} and \autoref{prop:MFdeg}, will establish a $k\ls{\beta}$-linear version of the degeneration property for some categories of singularities on quotient stacks.

\subsection{Preliminaries on categories of singularities on stacks}


There have been several concrete approaches to developing the general theory of singularity categories \cites{MR3366002,lin2013global,MR2910782}. We will mostly use the perspective of \cite{Preygel}, extended more recently in \cite{BRT}, which is more abstract but has the advantage of allowing one to deduce results about singularity categories directly from the analogous statement for derived categories of coherent sheaves (see \autoref{lem:main_SOD_lemma}). We will summarize the main definitions and lemmas we will use from \cite{Preygel}. As elsewhere in the paper $k$ denotes a field of characteristic zero.

\begin{defn} A Landau-Ginzburg (LG) model is a pair $(\X,W)$, where $\X$ is a smooth finite type $k$-stack such that the automorphism groups of its geometric points are affine and $W$ is a morphism $$W: \X \to \mathbb{A}^1.$$ \end{defn} 
In particular, $\X$ is a QCA stack over $k$ in the sense of \cite{drinfeld2013some}. Our primary examples of interest will be quotient stacks $\X:=X/G$ over $k$. By generic smoothness, if $W$ is non-constant on every component of $\X$, then $W$ has only finitely many critical values in $\bA^1$. Throughout this paper we let $\Crit(W)$ denote the component of the vanishing locus of $dW \in \Gamma(\X,\Omega^1_\X)$ which lies set theoretically in $\X_0:= \X \times_{\mathbb{A}^1} \{0\}$.

%
%

We now equip the bounded derived category of coherent sheaves on the zero fiber, $\DCoh(\X_0)$, with a $k\ps{\beta}$-linear structure, where $\beta$ is a variable of homological degree -2. This arises from a homotopical $\Sone$-action on the category $\DCoh_{\X_0}(\X)$, in the terminology of \cite{Preygel}, which concretely refers to a natural action of $H_\ast(S^1;k) \simeq \Lambda$ on the Hom-complexes of the category. The formal variable $\beta$ arises via the same construction which leads to the formal variable $u$ acting on $\hoch^-(\cA)$, but we use different variable names to avoid confusion between these two $S^1$-actions, especially when we discuss the $k\ls{\beta}$-linear negative cyclic homology below.

$\op{Spec}(\Lambda)$ admits the structure of a derived group scheme, so the $\infty$-category $\DCoh(\Lambda)$, as well as its ind-completion $\IC(\Lambda)$, admits a symmetric monoidal structure given by the convolution product ``$\circ$'': Given $F,G \in \DCoh(\Lambda)$, $F \circ G := m_\ast (F \boxtimes G)$, where $m : \op{Spec}(\Lambda) \times \op{Spec}(\Lambda) \to \op{Spec}(\Lambda)$ is the group multiplication. So the underlying complex of $F \circ G$ is $F \otimes_k G$, with $\Lambda$-module structure given by letting $B$ act on homogeneous elements by
\[
B_{F \otimes G}(m \otimes n) := B_F(m) \otimes n + (-1)^{|m|} m \otimes B_G(n).
\]
The following is an enrichment of standard Koszul duality results.

\begin{lem}\label{lem:Tate}
The functor
\begin{align*}
\DCoh(\Lambda) &\to \Perf(k\ps{\beta}) \\ 
V &\mapsto V^{\Sone}:=\RHom_{\Lambda}(k,V)
\end{align*}
extends to a symmetric monoidal equivalence, leading to a symmetric monoidal equivalence
$$\IC(\Lambda)^{\otimes} \cong (k\ps{\beta} \Mod)^{\otimes} $$
\end{lem}
\begin{proof}
This is \cite{Preygel}*{Proposition 3.1.4}, see also \cite{BRT}*{Remark 2.38, Lemma 2.39}.
\end{proof}
The proof of \autoref{lem:Tate} relies on an elementary but important observation. Let $(V,d)$ be a complex with a $\Lambda$-action. There is a quasi-isomorphism of complexes
$$V^{\Sone} \cong  (V\ps{\beta}, d+\beta B).$$
Formal completion does not commute with the formation of tensor products of complexes, but the formation of the complex $(V[\beta], d + \beta B)$ does commute with forming tensor products of complexes. So the crux of the proof of \autoref{lem:Tate} is the following:

\begin{lem} \label{lem:hom_bounded_above} 
The natural inclusion of complexes
$$ (V[\beta],d+ \beta B) \to (V\ps{\beta}, d + \beta B)$$
is a quasi-isomorphism whenever $(V,d)$ is \emph{homologically-bounded above}.
\end{lem}
\begin{proof}
By definition $(V\ps{\beta},d + \beta B)$ is the inverse limit of the complexes $(V[\beta]/(\beta^n), d + \beta B)$. If $V$ is homologically bounded above, then for any $m$ the canonical map $(V[\beta],d+\beta B) \to (V[\beta]/(\beta^n),d+\beta B)$ induces an isomorphism in homology of degree $\geq m$ for $n \gg 0$, which implies the claim.
\end{proof}

As described in \cite{Preygel}*{Construction 3.1.5}, \cite{BRT}*{Remark 2.38}, the stack $\X_0$ admits an action by the derived group scheme $\op{Spec}(\Lambda)$ which defines the upper horizontal arrow in the cartesian square: 
\begin{align} \label{eq:action}
\xymatrix{
\X_0 \times \op{Spec}(\Lambda)  \ar[r] \ar[d]_{p_1} & \X_0 \ar[d]^i\\
\X_0 \ar[r]^i & \X
}
\end{align}
The action \eqref{eq:action} can be described (\cite{Preygel}*{Remark 3.1.7}) concretely by a cosimplicial, commutative dg-$\mathcal{O}_{\X}$ algebra whose cosimplicial degree $n$ piece is given by $\mathcal{A}_n := \langle \mathcal{O}_{\X}[B_{\X}, B_1,\cdots, B_n], dB_i=0, dB_{\X}=W \rangle$, where each of the variables $B_i$ have degree one.


\begin{defn}\label{defn:MFcat} We define $\PreMF(\X,W) :=\DCoh(\X_0)$ with the additional $k\ps{\beta}$-linear structure induced by \autoref{lem:Tate} and the $\DCoh(\Lambda)^\otimes$-module structure on $\DCoh(\X_0)$ induced by the action \eqref{eq:action}.
\end{defn}

It is useful to note that the $k\ps{\beta}$-linear structure has a concrete dg-model which is described in \cite{Preygel} and which we now recall. Observe that $\mathcal{O}_{\X_0} \cong \mathcal{A}:=(\mathcal{O}_{\X}[B_{\X}],dB_{\X}=W)$, where $B_{\X}$ is a variable of degree one. Pushforward defines a canonical equivalence: $$\DCoh(\X_0) \cong \mathcal{A}\Mod(\Perf(\X)^\otimes),$$ where the right hand side denotes the category of coherent $\mathcal{A}$-modules. There are natural adjoint functors:\begin{align*} i_*:  \mathcal{A}\Mod(\Perf(\X)^\otimes) \to \DCoh(\X) \\ i^*: \DCoh(\X) \to \mathcal{A}\Mod(\Perf(\X)^\otimes)\end{align*} given by forgetting the $\cA$-module structure and tensoring with $\mathcal{A}$, respectively.
Now given two dg-$\mathcal{A}$-modules $M,N$ whose underlying complex of $\cO_\X$-modules has bounded coherent homology, the Hom-complex $\Hom_{\X}(i_*M,i_*N)$ inherits a $\Lambda$-module structure given by \begin{align} \label{eq: S1actin} B :  \phi \mapsto B_{\X} \circ \phi - (-1)^{|\phi|} \phi \circ  B_{\X} \end{align} The following proposition shows that this $\Lambda$-module structure is enough to recover Homs as $\mathcal{A}$-modules:


\begin{lem}[\cite{Preygel}*{Section 3.3}; see also \cite{Preygel}*{Proposition 3.2.1}] \label{lem:circle_actions}
Given objects $M, N \in  \mathcal{A}\Mod(\Perf(\X)^\otimes)$, we have an equivalence $\Hom_{\X}(i_*M,i_*N)^{\Sone} \cong \Hom_{\mathcal{A}}(M,N)$, where the $S^1$-action is given by \eqref{eq: S1actin}. \end{lem}

\begin{proof}
For any object $M \in \mathcal{A}\Mod(\Perf(\X)^\otimes)$, we can construct a complex of $B_\X$-modules, $M[B]$, where as usual $B$ is a variable of degree 1 with $B^2=0$ and the action of $B_\X$ is given by $(B_\X)_M + B$, where left multiplication by $B$ on $M[B]$ anti-commutes with the action of $(B_\X)_M$. There is an adjunction:\footnote{More conceptually, we have that $M[B]$ is isomorphic to $i^*i_*M = \mathcal{A} \otimes_{\cO_\X} M$. Both are isomorphic to $M \oplus M[1]$ as complexes of $\cO_\X$-modules, and the isomorphism $\cA \otimes_{\cO_\X} M \to M[B]$ which intertwines the action of $B_\X$ is $(m,m') \mapsto (m+B_\X m', m')$.}  \begin{equation} \label{eq:adjuncting} \begin{aligned} \Hom_{\X}(i_*M,i_*N) &\cong \Hom_\mathcal{A}(M[B],N) \\ \phi(m) &\mapsto \tilde{\phi}(m+Bm'):=\phi(m-B_\X \cdot m') + (-1)^{|\phi|} B_\X \phi(m') \\ \tilde{\psi}(m) := \psi(m+B\cdot 0) &\mapsfrom \psi(m+Bm')\end{aligned} \end{equation}
The $\Lambda = k[B]$-module structure, which on the left hand side is given by $B : \phi \mapsto B_\X \circ \phi - (-1)^{|\phi|} \phi \circ B_\X$, corresponds under this isomorphism to the $\Lambda$-module structure given by $B : \phi(-) \mapsto - \phi(B \cdot (-) )$.

We then extend $M[B]$ to a resolution of $M$ as an $\cA$-module by forming the complex
\begin{align}\label{eq:KoszulTate} \xymatrix{ (M[B]\otimes_k (k\ls{\beta} / \beta k\ps{\beta}); d = d_M - B\beta \cdot (-))  \ar[rr]^-{\cong}_-{B,\beta^{-n} \mapsto 0} & &  M } \end{align}
where $\beta$ is a variable of homological degree $-2$, which we refer to as the Koszul-Tate resolution of $M$. Using Equations \eqref{eq:adjuncting} and \eqref{eq:KoszulTate}, we have that $\Hom_{\mathcal{A}}(M,N)$ can be computed as
$$ \Hom_{\mathcal{A}}(M[B]\otimes_k (k\ls{\beta} / \beta k\ps{\beta}), N) \cong \Hom_{\mathcal{A}}(M[B],N)\ps{\beta} \cong \Hom_{\X}(i_*M,i_*N)^{\Sone},$$
where it is evident that the differential on the first term agrees with the differential used to compute the invariants for the $\Sone$-action defined in \eqref{eq: S1actin}.
\end{proof}

The natural $k\ps{\beta}$-linear structure on the complex $\Hom_{\X}(i_*M,i_*N)^{\Sone}$ provides an explicit model for the $k\ps{\beta}$-linear structure from \autoref{defn:MFcat}. Note that in \eqref{eq:KoszulTate}, we constructed a canonical quasi-isomorphism of $\cA$-modules $M \cong M[B] \otimes_k (k\ls{\beta} / \beta k\ps{\beta})$, where the latter has an explicit action by $k\ps{\beta}$. Under the resulting quasi-isomorphism \[\Hom_\cA(M,N) \simeq \Hom_\cA(M[B] \otimes_k (k\ls{\beta} / \beta k\ps{\beta}),N),\]the action of $\beta$ corresponds to $\phi(-) \mapsto \phi(\beta \cdot (-))$. The fact that composition is $k\ps{\beta}$-linear follows from the elementary calculation that composition $\Hom_\X(i_*M,i_*N) \otimes_k \Hom_\X(i_*N,i_*P) \to \Hom_\X(i_*M,i_*P)$ is $\Lambda$-linear, and \autoref{lem:Tate}.


With all of this in place, we turn to defining our main objects of interest, the categories of singularities: 

\begin{defn}
We define the category $\MF(\X,W)$ to be $$\MF(\X,W) :=\PreMF(\X,W)\itimes_{k\ps{\beta}} k\ls{\beta}.$$
\end{defn}

This definition is justified by the following lemma:

\begin{lem}[\cite{Preygel}*{Proposition 3.4.1}] $\MF(\X,W)$ is a dg-enhancement of the idempotent completion of the triangulated category $$H^0(\DCoh(\X_0))/H^0(\Perf(\X_0)).$$   \end{lem}
\begin{proof} 

Let $M \in \DCoh(\X_0)$. The lemma can be easily reduced to the following assertion: $M \in \Perf(\X_0)$ iff $\beta^n=0 \subset \Hom(M,M)$ for large enough $n$. To prove this assertion, recall that $M$ is perfect iff it is compact in $\QC(\X_0)$, because $\cX_0$ is QCA \cite{drinfeld2013some}*{Cor.~1.4.3}.

On one hand $\beta^n \in \Hom_\cA(M,M)$ corresponds under the quasi-isomorphism \eqref{eq:KoszulTate} to simply multiplying by $\beta^n$. We observe that the kernel of the surjective map $$\ker \left(\beta^n : M[B] \otimes (k\ls{\beta}/\beta k \ps{\beta}) \to M[B] \otimes (k\ls{\beta}/\beta k \ps{\beta}) \right) = M[B] \cdot \beta^{-n-1} \oplus \cdots \oplus M[B] \cdot \beta^0,$$ is a compact $\cA$-module because the associated graded of the $\beta$-adic filtration is a direct sum of finitely many copies of the compact $\cA$-module $M[B] \simeq \cA \otimes_{\cO_X} M$. Giving a null-homotopy of $\beta^n$, which is equivalent to giving a null homotopy of the composition $\beta^n \circ \id_{M}$, is equivalent to giving a factorization of $\id_{M[B] \otimes (k\ls{\beta}/\beta k \ps{\beta})}$ through the subcomplex $\ker(\beta^n)$. In this case $M$ is a homotopy retract of the compact $\cA$-module $\ker(\beta^n)$ and is thus compact.

Conversely, if $M[B] \otimes (k\ls{\beta}/\beta k \ps{\beta}) = \bigcup_n \ker(\beta^n)$ is a compact $\cA$-module, the identity morphism factors through some $\ker(\beta^n)$ for some $n$. This proves the lemma. 
\end{proof} 


There is another point of view on the $k\ps{\beta}$ linear structure which will be useful below. Let $(\X,W)$ be an LG-model. According to \cite{ben2013integral}*{Theorem 1.1.3}, there is an equivalence of categories
$$ \DCoh(\X_0) \cong \Fun^{ex}_{\Perf(\bA^1)^{\otimes}} (\Perf(k), \Perf(\X)) $$
It is not difficult to check that the $\DCoh(\Lambda)^\otimes$-module structure on the left-hand side of this equivalence corresponds to the natural $\Fun^{ex}_{\Perf(\bA^1)^\otimes} (\Perf(k), \Perf(k))$-module structure on the right-hand side. If $\A$ is a module category for some symmetric monoidal $\infty$-category $\cC^\otimes$, and $\cA = \sod{\A_i ; i\in I}$ is a possibly infinite semiorthogonal decomposition indexed by a totally ordered set $I$, we say that the semiorthogonal decomposition is $\cC^\otimes$-linear if $\cC \itimes \cA_i \to \cA$ factors through $\cA_i$, in which case it does so uniquely up to contractible choices.

\begin{lem}
Let $\cC^\otimes$ be a symmetric monoidal stable $\infty$-category, and let $\cB$ and $\cA$ be $\cC^\otimes$-module categories with $\cB$ compact in $\cC^\otimes\Mod(\op{Cat}_\infty^{perf})$. If $\cA = \sod{\cA_i ; i \in I}$ is a $\cC^\otimes$-linear semiorthogonal decomposition, then $\Fun_{\cC^\otimes}^{ex}(\cB,\cA_i) \to \Fun^{ex}_{\cC^\otimes}(\cB,\cA)$ is a fully faithful functor, and identifying the former with its essential image in the latter, we have a semiorthogonal decomposition
\[
\Fun_{\cC^\otimes}^{ex}(\cB,\cA) = \sod{\Fun_{\cC^\otimes}^{ex}(\cB,\cA_i) ; i\in I}
\]
\end{lem}
\begin{proof}
The fact that $\cB$ is compact as a $\cC^\otimes$-module category allows us to commute $\Fun_{\cC^\otimes}^{ex}(\cB,-)$ with filtered colimits and therefore reduce to the case of a finite index set $I$. Then by an inductive argument it suffices to prove the claim in the case where we have a two term semiorthogonal decomposition $\cA = \sod{\cA_0,\cA_1}$. If we let $\iota_i : \cA_i \hookrightarrow \cA$ denote the inclusion, and we let $\iota_1^R$ (respectively $\iota_0^L$) denote the right (respectively left) adjoint whose existence is guaranteed by the semiorthogonal decomposition. One can check that the composition functor $\iota_1^R \circ (-) : \Fun_{\cC^\otimes}^{ex}(\cB,\cA) \to \Fun_{\cC^\otimes}^{ex}(\cB,\cA_1)$ is a right adjoint to the composition functor $\iota_1 \circ (-) : \Fun_{\cC^\otimes}^{ex}(\cB,\cA_1) \to \Fun_{\cC^\otimes}^{ex}(\cB,\cA)$, and likewise $\iota_0^L \circ(-)$ is a left adjoint to $\iota_0 \circ (-)$. It is also straightforward to check that the canonical maps $\iota_0^L \circ \iota_0 \circ(-) \to \op{id}$ and $\op{id} \to \iota_1^R \circ \iota_1 \circ(-)$ are equivalences, and $\Map(\iota_1 \circ F, \iota_0 \circ G)$ is contractible for any functors $F \in \Fun_{\cC^\otimes}^{ex}(\cB,\cA_1)$ and $G \in \Fun_{\cC^\otimes}^{ex}(\cB,\cA_0)$. The claim follows.
\end{proof}

An immediate corollary of this is the following:


\begin{lem} \label{lem:main_SOD_lemma} Let $(\X,W)$ be an LG-model. Suppose that $\Perf(\X)$ admits a $\Perf(\mathbb{A}^1)^\otimes$-linear semiorthogonal decomposition $\sod{\A_i ; i\in I}$. Then $\MF(\X,W)$ admits a semiorthogonal decomposition by $k\ls{\beta}$-linear subcategories
\[
\MF(\X,W) = \sod{k\ls{\beta} \itimes_{k\ps{\beta}} \Fun^{ex}_{\Perf(\mathbb{A}^1)^\otimes} (\Perf(k),\A_i)}.
\]
\end{lem}
\begin{proof}
By the previous lemma applied to $\A =\Perf(\X)$, we obtain a semiorthogonal decomposition $\PreMF(\X,W) = \sod{\Fun^{ex}_{\Perf(\mathbb{A}^1)^\otimes} (\Perf(k),\A_i)}$. The semiorthogonal decomposition is $k\ps{\beta}$-linear because $\beta$ acts via endo-functors of $\Perf(k)$ as a $\Perf(\bA^1)^\otimes$-module category. Finally, the localization functor from $k\ps{\beta}$-linear categories to $k\ls{\beta}$-linear categories commutes with filtered colimits, so one gets the desired semiorthogonal decomposition of $\MF(\X,W) := k\ls{\beta} \itimes_{k\ps{\beta}} \PreMF(\X,W)$ by base changing semiorthogonal decompositions.
\end{proof}

\subsection{Motivic decompositions and degeneration for \texorpdfstring{$\MF$}{MF}}

In this section we prove an analog of \autoref{thm:motivic_main} for the $k\ls{\beta}$-linear dg-category $\MF(X/G,W)$. For any $\cA \in \op{LinCat}_{k\ls{\beta}}^{sm}$, let $\bG(\cA) \subset \op{LinCat}_{k\ls{\beta}}^{sm}$ denote the smallest $\infty$-subcategory containing $\cA$ that is closed under splitting countable semiorthogonal decompositions in the following sense: for any $\cC \in \op{LinCat}_{k\ls{\beta}}^{sm}$ which has a $\bZ$-indexed semiorthogonal decomposition $\cC = \langle \cC_i \rangle_{i \in \bZ}$, $\cC \in \bG(\cA)$ if and only if $\cC_i \in \bG(\cA)$ for all $\cA$.

\begin{thm} \label{thm:motivic_main_lg}
Let $G$ be an algebraic group. Let $X$ be a smooth $G$-quasiprojective $k$-scheme with a semi-complete KN stratification, and let $W : X/G \to \bA^1$ be a morphism. Then there is a smooth projective-over-affine Deligne-Mumford stack $\Y$ with a map $W' : \Y \to \bA^1$ such that
\[\MF(X/G,W) \in \bG(MF(\Y,W')) \subset \op{LinCat}_{k\ls{\beta}}^{sm}.\]
Furthermore if $\Perf(\Crit(W)/G)$ is a proper dg-category, then the pair $(\Y,W')$ can be chosen so that $\Y$ is a smooth variety and $W' \colon \Y \to \bA^1$ is projective.
\end{thm}

Note that by \autoref{cor:KN_strat_properness}, $\Perf(\Crit(W)/G)$ is a proper dg-category if and only if the induced KN stratification on $\Crit(W)$ is complete.

\begin{rem}
There is a slightly cleaner formulation of \autoref{thm:motivic_main_lg} using $k\ls{\beta}$-linear additive noncommutative motives, analogous to \autoref{thm:motivic_main}. \cite{tabuada2008higher} constructs additive noncommutative motives for dg-categories which are linear over a commutative ring, and \cite{blumberg2013universal} constructs additive noncommutative motives over the sphere spectrum. The methods of \cite{blumberg2013universal} appear to apply verbatim to construct the $\infty$-category of additive noncommutative motives over an arbitrary $E_\infty$-algebra $R$, such as $k\ls{\beta}$, but in the interest of space we have formulated \autoref{thm:motivic_main_lg} to avoid developing this additional machinery.
\end{rem}

Before proving the theorem, we note the following analogues of \autoref{lem:KN_stratification} and \autoref{lem:rational_morphism}:

\begin{lem} \label{lem:KN_stratification_lg}
If $X$ is a smooth $G$-scheme with a KN stratification, and $W : X/G \to \bA^1$ is an LG-model, then $\MF(X/G,W) \in \bG\left(\MF(X^{ss}/G,W) \oplus \bigoplus_i \MF(Z_i/L_i,W|_{Z_i/L_i}) \right)$ in $\op{LinCat}_{k\ls{\beta}}^{sm}$.
\end{lem}
\begin{proof}
The main semiorthogonal decomposition of \cite{halpern2015derived} extends to categories of singularities by \autoref{lem:main_SOD_lemma}, and hence the argument of \autoref{lem:KN_stratification} applies verbatim to $\MF(X/G,W)$.
\end{proof}

\begin{lem} \label{lem:rational_morphism_lg}
Let $\pi : \Y \to \X$ be a rational morphism of finite-type $k$-stacks, i.e., $R\pi_\ast \cO_\Y \simeq \cO_\X$, and let $W : \X \to \bA^1$ be a morphism. Assume that $\X$ is smooth and $\pi_\ast$ preserves $\DCoh$. Then $\MF(\X,W)$ is a semiorthogonal factor of $\MF(\Y,W)$.
\end{lem}

\begin{proof}
The functors $\pi_\ast$ and $\pi^\ast$ are are $\Perf(\bA^1)$-linear, and it follows from \autoref{lem:main_SOD_lemma} that the semiorthogonal decomposition of $\Perf(\Y)$ in the proof of \autoref{lem:rational_morphism} induces a semiorthogonal decomposition of $\MF(\Y,W)$.
\end{proof}

\begin{proof}[Proof of \autoref{thm:motivic_main_lg}]

The proof of \autoref{thm:motivic_main} mostly applies verbatim, with the following substitutions: \autoref{lem:KN_stratification_lg} in place of \autoref{lem:KN_stratification}; \autoref{lem:rational_morphism_lg} in place of \autoref{lem:rational_morphism}; the properness of the dg-category $\Perf(\Crit(W)/G)$ in place of the properness of the dg-category $\Perf(U/G)$; and the result \autoref{prop:LG_motivic} below in place of \cite{bergh2016geometricity}*{Theorem 6.6} to reduce from the case of a projective-over-affine DM stack to a quasi-projective scheme.


The only parts of the proof of \autoref{thm:motivic_main} which require modification in the current context are those which have to do with the properness of the dg-category $\Perf(\Crit(W)/G)$, and this only affects two cases of the inductive proof. We thus re-write these cases, indicating the necessary modifications:


\medskip
\noindent{\textit{Case $X^{ss} \neq \emptyset$, $X^{ss} \neq X^s$, and $\op{codim}(X^{ss} \setminus X^s, X^{ss}) \geq 2$:}}
\medskip

It suffices by the inductive hypothesis to prove the claim for $(X^{ss}/G,W)$. We apply the inductive partial resolution procedure of \cite{kirwan1985partial}: Let $Y \subset X^{ss}$ be the locus of points whose stabilizer has maximal dimension. Then $Y$ is a smooth closed subvariety, and the blow up $X' := \op{Bl}_{Y}(X^{ss})$ has a KN stratification induced by a relatively ample bundle such that $(X')^{ss}$ has lower-dimensional stabilizers.

Consider a point $x \in Y$, let $R\subset G$ be the (reductive) stabilizer subgroup. The $G$-invariance of $W$ implies that $(dW)_x  \in (\Omega^1_{X^{ss},x})^R \subset \Omega^1_{X^{ss},x}$. Because $(X^{ss})^G \subset Y$ the restriction map $(\Omega^1_{X^{ss},x})^R \to \Omega^1_{Y,x}$ is injective, and hence the pullback map $(\Omega^1_{X^{ss},x})^R \to \Omega^1_{X',y}$ is injective for any point $y$ in the fiber of $x$ under $p : X'=\op{Bl}_{Y} X^{ss} \to X^{ss}$. It follows that for any $y \in p^{-1}(x)$, $(dW)_x = 0$ if and only if $(d(W|_{X'}))_y = 0$, and therefore $\Crit(W|_{X'}) = p^{-1}(\Crit(W))$.

Now by \autoref{lem:rational_morphism_lg} and \autoref{lem:KN_stratification_lg}, if $Z_i$ are the centers of the KN stratification of $X'$, then
\[
\MF(X^{ss}/G,W) \in \bG\left(\MF((X')^{ss}/G,W|_{X'}) \oplus \bigoplus_i \MF(Z_i/L_i,W|_{Z_i/L_i}) \right).
\]
Note that if $\Perf(\Crit(W)/G)$ is a proper dg-category, then so is $\Perf(\Crit(W|_{X'})/G)$, using criterion (2) of \autoref{lem:properness} and the fact that $\Crit(W|_{X'}) \to \Crit(W)$ is proper. It follows from \autoref{lem:KN_strat_properness} that in this case $\Perf(\Crit(W|_{(X')^{ss}})/G)$ and $\Perf(\Crit(W|_{Z_i})/L_i)$ are proper dg-categories as well. The claim holds for the categories $\MF(Z_i/L_i,W|_{Z_i/L_i})$ by the inductive hypothesis, so it suffices to show the claim for $\MF((X')^{ss}/G,W|_{X'})$. Because the dimension of the stabilizer groups of $(X')^{ss}$ are strictly smaller than that of the stabilizer groups in $X^{ss}$, we can replace $X^{ss}$ with $(X')^{ss}$ and iterate this construction until $X^{ss} = X^{s}$, which is handled in a previous case.

\medskip
\noindent \textit{Case $X^{ss} \neq \emptyset$ but $\op{codim}(X^{ss}\setminus X^{s},X^{ss}) \leq 1$:}
\medskip

Here the only modification needed to the proof of \autoref{thm:motivic_main} is the observation that $\Crit(W|_{U\times Y}) = \Crit(W|_U) \times Y$, so $\Perf(\Crit(W|_{U \times Y}))$ is a proper dg-category if $\Perf(\Crit(W|_U)/G)$ is.

\end{proof}



In the proof of \autoref{thm:motivic_main_lg}, we used the following proposition to reduce from the case of a DM stack to the case of a variety. The proof is an adaptation of the methods of \cite{bergh2016geometricity} to categories of matrix factorizations.

\begin{prop} \label{prop:LG_motivic}
Let $W : \X \to \bA^1$ be an LG-model such that $\X$ is a smooth quasi-projective (in the sense of \cite{kresch2009geometry}) DM stack and $\Crit(W)$ is proper. Then there is a smooth variety $Y$ and a projective morphism $W':Y \to \bA^1$ such that $\MF(\X,W) \in \bG(\MF(Y,W')) \subset \op{LinCat}^{sm}_{k\ls{\beta}}$.
\end{prop}

\begin{proof}
As observed in \cite{kresch2009geometry}, $\X = X/G$ for some reductive $G$ and $G$-equivariant locally closed embedding $X \hookrightarrow (\bP^N)^s$ for some linearized action of $G$ on a projective space. Considering the graph of the map $W : X \to \bA^1$ gives a $G$-equivariant locally closed embedding $X \hookrightarrow (\bP^N \times \bA^1)^s$, where $G$ acts trivially on the $\bA^1$ factor. Using Kirwan's resolution algorithm \cite{kirwan1985partial}, as in \cite{kresch2009geometry}, we modify $(\bP^N \times \bA^1)^{ss}$ by a series of smooth blow ups away from $(\bP^N \times \bA^1)^{s}$ such that the resulting semistable locus is a DM stack $\Y$ which is smooth and projective over $\bA^1$. Thus taking the closure of $\cX$ in $\Y$ and resolving singularities, we can produce a Nagata compactification $\bar{W} : \bar{\X} \to \bA^1$ of the original morphism $W : \X \to \bA^1$ such that $\bar{\cX}$ is projective over $\bA^1$ and smooth over the ground field.

Because $\X$ and $\bar{\X}$ are DM and the critical locus $\Crit(W)$ is proper, it follows that $\Crit(W)$ is a union of connected components of $\Crit(\bar{W})$. Therefore, the category $\MF(\bar{\X},\bar{W})$ splits as a direct sum of the subcategory consisting of objects supported on $\Crit(W)$ and the subcategory consisting of objects supported on other components of $\Crit(\bar{W})$ \cite{Preygel}*{Prop.~4.1.6}. So it suffices to prove the claim for $\MF(\bar{\X},\bar{W})$, i.e. we may assume that the potential itself $W : \X \to \bA^1$ is projective rather than assuming that just $\Crit(W)$ is proper.

First we reduce to the case where $\X$ has generically trivial stabilizer. Because $\X$ is a global quotient stack, we may find a vector bundle $V$ over $\X$ on which the automorphism groups act faithfully. Then $\bP(V \oplus \cO_\X) \to \X$ is a rational morphism, so by \autoref{lem:rational_morphism_lg} the claim for $\bP(V \oplus \cO_\X)$ implies the claim for $\X$, and the generic stabilizer of the former is trivial. Furthermore, the function $W$ restricted to $\bP(V \oplus \cO_\X)$ will still be proper.

Next we consider Bergh's destackification, constructed in \cite{bergh2014functorial}
$$\xymatrix{ & \X' \ar[dl]_{\pi} \ar[dr]^f & \\ X' & & \X},$$
where $X'$ is the coarse moduli space and is smooth, the morphism $\pi$ is a composition of root stacks along smooth divisors and $f$ is a composition of root stacks along smooth divisors and blow ups along smooth centers. Note that the morphism $f$ is rational, so by \autoref{lem:rational_morphism_lg} it suffices to prove the claim for the composition $W' : \X' \to \X \to \bA^1$. Note that $W'$ will still be proper.

Finally, because $X'$ is the coarse space of $\X'$, the function $W'$ descends uniquely to $W' : X' \to \bA^1$, and this map is still proper. Because $\pi : \X' \to X'$ is a composition of root stacks along smooth divisors, it suffices to prove the following claim: if $\Y$ is a smooth DM stack and $\X \to \Y$ is a root stack along a smooth divisor in $\Y$ and $W : \Y \to \bA^1$ is a proper map, then the claim of the proposition for $\Y$ implies the claim for $\X$.

Let $\cD \hookrightarrow \Y$ be the Cartier divisor used to form the root construction, and consider the diagram
$$\xymatrix{\cD \times B\mu_{r} \ar[r]^i \ar[d]^{\pi_{\cD}} & \X \ar[d]^\pi \\ \cD \ar[r] & \Y}.$$
Then \cite{ishii2011special}*{Proposition 6.1} shows that we have a semiorthogonal decomposition
$$\Perf(\X) = \sod{i_\ast \pi_\cD^\ast (\Perf(\cD)) \otimes M^{r-1} , \ldots, i_\ast \pi_\cD^\ast (\Perf(\cD)) \otimes M , \pi^\ast(\Perf(\Y)) },$$
where $M$ is the universal invertible sheaf coming from the root stack construction. Given a map $W : \Y \to \bA^1$ and a module $F \in \Perf(\bA^1)$, it is evident from the description of this semiorthogonal decomposition that each subcategory is preserved by the functor $\pi^\ast W^\ast(F) \otimes (-)$. If follows that this is a semiorthogonal decomposition of $\Perf(\bA^1)^\otimes$-module categories, and thus by \autoref{lem:main_SOD_lemma} there are induced semiorthogonal decompositions
\begin{gather*}
\Perf(\X_0) = \sod{i_\ast \pi_{\cD_0}^\ast (\Perf(\cD_0)) \otimes M^{r-1} , \ldots, i_\ast \pi_{\cD_0}^\ast (\Perf(\cD_0)) \otimes M , \pi^\ast(\Perf(\Y_0)) } \\
\DCoh(\X_0) = \sod{i_\ast \pi_{\cD_0}^\ast (\DCoh(\cD_0)) \otimes M^{r-1} , \ldots, i_\ast \pi_{\cD_0}^\ast (\DCoh(\cD_0)) \otimes M , \pi^\ast(\DCoh(\Y_0)) },
\end{gather*}
where the subscript $0$ refers to the derived zero fiber of $W$ and its restrictions to $\cD$, $\X$, and $\cD \times B\mu_r$. Thus we have a semiorthogonal decomposition of $k\ls{\beta}$-linear categories
$$\MF(\X,W \circ \pi) \simeq \sod{\MF(\cD,W|_\cD),\ldots,\MF(\cD,W|_\cD),\MF(\Y,W)}$$
and the claim of the proposition for $\X$ follows from the claim for $\Y$ and induction, because $\cD$ is smooth of one lower dimension.

\end{proof}

\subsubsection{The degeneration property for quotient stacks}

For any $\A \in \op{LinCat}_{k\ls{\beta}}^{sm}$, we may form the $k\ls{\beta}$-linear Hochschild complex $\hoch^{k\ls{\beta}}(\A)$, which is a module over the CDGA $\Lambda\ls{\beta}=k\ls{\beta}[B]/(B^2)$. $\hoch^{k\ls{\beta}}(\A)$ is computed by a Barr complex
\begin{equation} \label{eqn:hochschild}
\hoch^{k\ls{\beta}}(\A) := \bigoplus_{o_1,o_2,\cdots, o_n}\Hom(o_1,o_2) \otimes_{k\ls{\beta}} \Hom(o_2,o_3) \otimes_{k\ls{\beta}} \cdots \otimes_{k\ls{\beta}} \operatorname{Hom}(o_n,o_1),
\end{equation}
where $o_i$ are objects of $\A$, the differential is given by the usual formula for the Hochschild complex of a $k$-linear dg-category \cite{mitchell}, and $B$ acts by the Connes differential. Using the $\Lambda\ls{\beta}$-module structure, we may form the associated complexes $\hoch^{k\ls{\beta},(n)} (\A)$, $\hoch^{k\ls{\beta},-} (\A)$, and $\hoch^{k\ls{\beta},per} (\A)$ just as in the $k$-linear case (see Section \ref{sect:hochschild}).

\begin{defn} \label{def:Z2_degen_property} A $k\ls{\beta}$-linear dg-category $\A$ is said to have the \emph{$k\ls{\beta}$-linear degeneration property} if $H_\ast(\hoch^{k\ls{\beta},(n)}(\A))$ is a flat $k[u] / (u^n)$-module for all $n \geq 1$.  \end{defn}

\begin{lem} \label{lem:degeneration_sums_lg}
Let $\cA \in \op{LinCat}^{sm}_{k\ls{\beta}}$ satisfy the $k\ls{\beta}$-linear degeneration property, then any other category $\cC \in \bG(\cA)$ satisfies the $k\ls{\beta}$-linear degeneration property.
\end{lem}
\begin{proof}
The proof of \autoref{lem:degeneration_sums} applies verbatim to show that if $\cC = \langle \cC_i \rangle$ is a $\bZ$-indexed semiorthogonal decomposition of $k\ls{\beta}$-linear dg-categories, then $\cC$ satisfies the degeneration property if and only if each of the $\cC_i$ do. This implies that the full $\infty$-subcategory of $\op{LinCat}_{k\ls{\beta}}^{sm}$ consisting of categories that satisfy the $k\ls{\beta}$-linear degeneration property contains $\cA$ and is closed under splitting countable semiorthogonal decompositions, so it contains $\bG(\cA)$ by definition.
\end{proof}

Our main degeneration result for categories of matrix factorizations is the following:
\begin{prop} \label{prop:MFdeg}
Let $W : X/G \to \bA^1$ be an LG-model, where $X$ is a smooth $G$-quasi-projective scheme which admits a semi-complete KN stratification. If $\Perf(\Crit(W)/G)$ is a proper dg-category, then the $k\ls{\beta}$-linear degeneration property holds for $\MF(X/G,W)$.
\end{prop}
\begin{proof}
By \autoref{thm:motivic_main_lg} and \autoref{lem:degeneration_sums_lg}, this reduces to the $k\ls{\beta}$-linear degeneration property for $\MF(Y,W)$ where $Y$ is a smooth scheme and $W : Y \to \bA^1$ is a projective map. This amounts to the degeneration of the $W$-twisted Hodge-de Rham complex by \cite{efimov_cyclic}*{Theorem 1.3}, which is established in \cite{ogus2007nonabelian}*{Theorem 4.22}.
\end{proof}

\begin{ex}
The hypotheses of the proposition are satisfied if $X$ is projective-over-affine and $\Gamma(X,\cO_{\Crit(W)})^G$ is finite dimensional -- see \autoref{ex:projective_over_affine}.
\end{ex}

We also have the following, which was left as an assumption in the original version of this paper:
\begin{cor} \label{cor:LG_degeneration}
If $\X$ is a smooth quasi-projective DM stack and $W : \X \to \bA^1$ is a map such that $\Crit(W)$ is proper, then the $k\ls{\beta}$-linear degeneration property holds for $\MF(\X,W)$.
\end{cor}
\begin{proof}
\autoref{prop:LG_motivic} reduces this to the case of a projective morphism $W : Y \to \bA^1$, with $Y$ smooth, which as noted above follows from \cite{ogus2007nonabelian}*{Theorem 4.22}.
\end{proof}

\begin{rem} 
We briefly discuss more concrete descriptions of $\hoch^{k\ls{\beta}}(\MF(\X,W))$ when $\X$ is a smooth separated Deligne-Mumford stack. Namely, for any affine $U$ with an \'{e}tale map $U \to \X$, let $\Omega_\bullet(U,W|_U)$ denote the Tate construction on $\bigoplus \Omega^i (U)[i]$ with respect to the $\Sone$-action given by $-dW \wedge$, i.e. $\Omega_\bullet(U,W|_U)$ is the $\Lambda\ls{\beta}$-module $\bigoplus \Omega^i_U\ls{\beta}[i]$ with differential $-\beta \cdot dW\wedge$. Letting $B$ act on $\Omega_\bullet(U,W|_U)$ via the de Rham differential as usual, $\Omega_\bullet(-,W)$ defines a sheaf of $\Lambda\ls{\beta}$-modules on the small \'{e}tale site $\X^{\tinyop{aff}}_{et}$. We define the global de Rham complex to be the $\Lambda\ls{\beta}$-module $$\Omega_\bullet (\X,W) := R\Gamma \left( \X^{\tinyop{aff}}_{et}, \Omega_\bullet(-,W) \right).$$
Assume for simplicity that $\Crit(W)=\operatorname{Crit}(W)$. Then by combining the approach of \autoref{prop:Toen_de_rham} with \cite{Preygel}*{Theorem 8.2.6}, it is not difficult to show that there is a natural isomorphism of $\Lambda\ls{\beta}$-modules
$$C^{k\ls{\beta}}_\bullet(\MF(\X,W)) \simeq \Omega_{\bullet}(I^{cl}_\X,W),$$
which induces an equivalence $C^{k\ls{\beta},per}_\bullet(\MF(\X,W)) \simeq \Omega_\bullet(I^{cl}_\X,W)^{\Tate}$.  Thus, \autoref{cor:LG_degeneration} is equivalent to the statement that the $k\ls{\beta} \otimes \Lambda$-module $\Omega_\bullet (\X,W)$ has the degeneration property. This constitutes a slight generalization of the degeneration results of \cite{ogus2007nonabelian} to DM stacks.
\end{rem}

\subsection{Graded Landau-Ginzburg models}

We will use these results to establish a large class of examples of $k$-linear dg-categories for which the usual $k$-linear degeneration property holds.

\begin{defn} \label{def:graded_lg}
A graded LG-model is a map $W:\X \to \mathbb{A}^1/\Gm$, where $\X$ is a smooth algebraic $k$-stack whose automorphism groups at geometric points are affine, and $\Gm$ acts on $\mathbb{A}^1$ with weight one.
\end{defn}

Note that the data of a graded LG-model is equivalent to specifying an invertible sheaf $L$ on $X/G$, which is classified by the composition $\X \to \mathbb{A}^1/\Gm \to B\Gm$, and a section $W \in \Gamma(X/G,L)$. Denote by $\phi : \X' \to \X$ the $\Gm$-torsor over $\X$ associated to $L$, i.e., $\X' = \op{Tot}_\X(L^\dual) \setminus 0$. To any graded LG-model, we use the term \emph{associated LG-model} to denote the pair $(\X',\phi^*W : \X' \to \bA^1)$. We will see in \autoref{prop:forget_Z_gradings} that, in a precise sense, the graded LG-model is a refinement of its associated LG-model.\footnote{Note also that given an LG-model $(\X,W)$, we can forget the data of the trivialization of $L$ to obtain a graded LG-model. This will correspond to forgetting the $k\ls{\beta}$-linear structure on the category $\MF(\X,W)$.}


In the setting of graded LG-models, for $F \in \DCoh(\X_0)$, there is a distinguished triangle: $$ F \otimes L^{-1}[1] \to i^*i_* F \to F$$ giving rise to a natural transformation $\beta_L : -\otimes L[-2] \to \operatorname{id} $ which in the non-graded case, where $L \simeq \cO_\X$, is the natural transformation $\operatorname{id}[-2] \to \operatorname{id}$ induced by the $k\ps{\beta}$-linear structure of \autoref{defn:MFcat}. To make this natural transformation more explicit, we introduce an analog of the Koszul-Tate resolution \eqref{eq:KoszulTate}, $\operatorname{KT}(M) \cong M$, and construct a very concrete natural transformation $\beta_L : \operatorname{KT}(-) \otimes L[-2] \to \operatorname{KT}(-)$.

First consider the dg-algebra $\cA := \cO_\X \oplus {L}^{-1} \cdot B_\X$ where $B_\X$ is a formal variable of degree $1$, and $dB_\X = W$, i.e., the differential is trivial on $\cO_\X$ and acts on the second summand by the map $L^{-1} \to \cO_\X$ that is the defining section $W$. An $\cA$-module consists of an $\cO_\X$-module $M$ along with a map $B_\X : L^{-1} \otimes M \to M[-1]$ of $\cO_\X$ modules satisfying the Leibnitz rule $d B_\X(s \otimes m) = W s m - B_\X(s \otimes dm)$, where $s$ is a local section of $L^{-1}$ and $m$ is a homogeneous local section of $M$. The $\infty$-category of $\cA$-modules is equivalent to the $\infty$-category $\DCoh(\X_0)$.

Given a module $M$ over $\cA$, we can form the analogue of $M[B]$ above, which we denote
\[\widehat{M} :=  M\oplus (L^{-1} \otimes M \cdot B),\]
where the differential is just $d_M \oplus L^{-1} \otimes d_M$. We have two operators $B^M_\X: L^{-1} \otimes M \to M [-1]$, defined by the $\cA$-module structure of $M$,
and $B: L^{-1} \otimes M \to L^{-1} \otimes M$, which is just the identity. Using these two operators we can form a canonical action of $\mathcal{A}$ on $\widehat{M}$, where the action of $B_\X$ on the first component is just $B^M_\X + B$, and the action on the second component is $L^{-1} \otimes B^M_\X$.


Let $Q(r): L^{-r} \otimes \widehat{M} \to L^{-r+1} \otimes \widehat{M}$ be the map which sends the first component of $L^{-r} \otimes \widehat{M}$ to the second component of $L^{-r+1} \otimes \widehat{M}$ by the identity. Then we define the Koszul-Tate resolution as
\[\operatorname{KT}(M):= \left(\bigoplus_{r \geq 0} L^{-r} \otimes \widehat{M} [-2r]; d = d_{\widehat{M}} + \sum_{r>0} Q(r) \right).\]
After pulling back to $\X'$, this complex can be identified with $M[B] \otimes (k\ls{\beta} / \beta k\ps{\beta})$, so it follows that the canonical map $\operatorname{KT}(M) \to M$, which annihilates all $r>0$ summands and $B$, is a quasi-isomorphism. We now define the natural map
\[
\beta_L : \operatorname{KT}(M) \otimes L [-2] \to \operatorname{KT}(M)
\]
to be the map of $\cA$-modules which annihilates the $r=0$ summand and identifies the $(r+1)^{\rm st}$ summand of $\operatorname{KT}(M) \otimes L[-2]$ with the $r^{\rm th}$ summand of $\operatorname{KT}(M)$.

\begin{defn}
We define the graded singularity category $\DSing(\X,W)$ to be the idempotent completion of the dg-category with the same objects as $\DCoh(\X_0)$ but with morphisms between $M,N$ given by 
$$\Hom_{\DSing(\X,W)}(M,N):=\hocolim_p \Hom_{\DCoh(\X_0)}(M,N\otimes L^{-p})[2p],$$
where the homotopy colimit is formed with respect to the natural maps $\beta_L : L^{-p} \otimes N[2p] \to L^{-p-1} \otimes N[2p+2]$.
\end{defn}

\begin{ex} \label{ex:graded_lg_model} Let $Y$ be a smooth variety over $k$, $\mathcal{E}$ a vector bundle over $Y$ and let $s \in \Gamma(\mathcal{E})$ be a regular section. We have an action of $\Gm$ on $\op{Tot}(\mathcal{E}^{\vee})$ by scaling in the fibers. The function $s$ therefore determines a mapping  $W_s: \X=\op{Tot}(\mathcal{E}^{\vee})/\Gm \to \mathbb{A}^1/\Gm$. The main theorem of \cite{Isik} gives an equivalence of ($\bZ$-graded) dg-categories $\DSing(\op{Tot}(\mathcal{E}^{\vee})/\Gm,W_s) \simeq \DCoh(s^{-1}(0))$. The construction of this equivalence works equally well when $Y$ has an action of a linear algebraic group $G$, $\cE$ is a $G$-equivariant locally free sheaf, and $s$ is $G$-invariant (see for instance \cite{hirano}). This gives an equivalence of $k$-linear dg-categories
\[
\DSing(\op{Tot}(\mathcal{E}^{\vee}) / (G \times \Gm), W_s) \simeq \DCoh(s^{-1}(0)/G).
\]
In this case the associated LG-model is simply the map $W' : \op{Tot}(\cE^\dual)/G \to \bA^1$ obtained by forgetting the $\Gm$ action.
\end{ex}

For a $\bZ$-graded $k$-linear dg-category $\cC$, we may tensor with $k\ls{\beta}$, thereby collapsing the grading on $\Hom(E,F)$ to a $\bZ/2\bZ$-grading. The following proposition describes the relationship between a graded LG-model and its associated LG-model.

\begin{prop}\label{prop:forget_Z_gradings}
Let $W : \X \to \mathbb{A}^1/\Gm$ be a graded LG-model, and let $W':\X' \to \mathbb{A}^1$ be the associated LG-model. Then we have a canonical equivalence of $\bZ/2\bZ$-graded dg-categories $\DSing(\X,W)\itimes_k k\ls{\beta} \simeq \MF(\X',W')$.
\end{prop}

\begin{lem} \label{lem:equivalent_degeneration_properties}
Let $\pi : \X \to \Y$ be a smooth affine morphism of QCA stacks. Then objects of the form $\pi^\ast F$ for some $F \in \DCoh(\Y)$ split generate $\DCoh(\X)$, i.e. $\DCoh(\X)$ is the smallest subcategory containing these objects which is closed under shifts, cones, and retracts.
\end{lem}

\begin{proof}
First note that the analogous claim holds for perfect stacks using $\Perf$ instead of $\DCoh$ and assuming only that $\pi$ is affine. Indeed, the fact that the pushforward functor $\pi_\ast : \QC(\X) \to \QC(\Y)$ is conservative implies that objects of the form $\pi^\ast F$ with $F \in \Perf(\Y)$ split generate $\QC(\X)$.

In order to conclude the same for $\DCoh$, we must imitate this argument for the categories $\IC$, i.e., show that $\pi_{\IC,\ast}$ is conservative. The pushforward $\pi_{\IC,\ast}$ again has a left adjoint $\pi^\ast_{\IC}$ which preserves $\DCoh$ and agrees with the usual pullback functor there. Because $\pi_{\IC,\ast}$ satisfies base change with respect to the shriek pullback \cite{gaitsgory2013indcoherent}*{5.2.5}, and $\IC$ satisfies fppf descent with respect to shriek pullback, it suffices to show this when $\Y = Y$ is an affine derived scheme and hence $\X = X$ is as well. In this case \cite{gaitsgory2013indcoherent}*{Proposition 4.5.3} shows that the essential image of the functor $\QC(\X) \widehat{\otimes}_{\QC(\Y)} \IC(\Y) \to \IC(\X)$ generates the latter category (in fact this functor is an equivalence). It follows that $\IC(\X)$ is generated by objects of the form $E \otimes \pi^{\IC,\ast}(F)$ for $E \in \Perf(\X)$ and $F \in \DCoh(\Y)$. Furthermore, because $\X$ is affine the category $\Perf(\X)$ is split generated by $\cO_\X$, so $\IC(\X)$ is generated by objects of the form $\pi^{\IC,\ast}(F)$ with $F \in \DCoh(\Y)$. It follows by adjunction that $\pi_{\IC,\ast}$ is conservative.
\end{proof}

\begin{proof}[Proof of \autoref{prop:forget_Z_gradings}]
Pullback along $\phi : \X' \to \X$ defines a functor of dg-categories $\DCoh(\X_0) \to \DCoh(\X'_0)$ which intertwines the action of $\beta_L$ and $\beta_{L|_{\X'}}$. Making use of the canonical trivialization $L|_{\X'} \simeq \cO_{\X'}$, we get a canonical map
\[
\phi^\ast : \hocolim_p \Hom_{\DCoh(\X_0)}(M,N\otimes L^{-p})[2p] \to \hocolim_p \Hom_{\DCoh(\X'_0)}(\phi^\ast(M),\phi^\ast(N))[2p].
\]
The latter can be identified with
\[
\Hom_{\PreMF(\X',W)}(\phi^\ast(M),\phi^\ast(N)) \otimes_{k\ps{\beta}} k\ls{\beta} \cong \Hom_{\MF(\X',W)}(\phi^\ast(M),\phi^\ast(N)).
\]
Because $\DSing(\X,W)$ is generated by objects of $\DCoh(\X_0)$, this extends to a dg-functor on idempotent completions $\phi^*: \DSing(\X,W) \to \MF(\X',W')$. From the universal property of the base-change category, $\phi^\ast$ admits an essentially unique $k\ls{\beta}$-linear extension
\[
\phi^\ast_{k\ls{\beta}} : \DSing(\X,W) \itimes_k k\ls{\beta} \to \MF(\X',W),
\]
and we will show that this functor is an isomorphism.

Concretely, $\phi_{k\ls{\beta}}^*$ maps an object $F \in \DCoh(\X_0)$, regarded as a generator for $\DSing(\X,W) \itimes_k k\ls{\beta}$, to the object $\phi^\ast(F) \in \DCoh(\X'_0)$, regarded as a generator of $\MF(\X',W')$. \autoref{lem:equivalent_degeneration_properties} implies that objects of the form $\phi^\ast(F)$ with $F \in \DCoh(\X_0)$ generate $\MF(\X',W)$, so it suffices to show that $\phi^\ast_{k\ls{\beta}}$ is fully faithful on the full subcategory of $\DSing(\X,W) \itimes_k k\ls{\beta}$ spanned by objects of $\DCoh(\X_0)$.

For $M,N \in \DCoh(\X_0)$, we consider $$\Hom_{\DSing(\X,W)}(M,N)\otimes k\ls{\beta} := \bigoplus_q \hocolim_{p-q} \Hom_{\DCoh(\X_0)}(M,N \otimes L^{q-p})[2(p-q)][2q].$$ Commuting colimits and reshuffling indices, this is isomorphic to $$ \cong \hocolim_p \bigoplus_q \Hom_{\DCoh(\X_0)}(M,L^q\otimes N \otimes L^{-p})[2p]$$  In this presentation, the operator $\beta_L$ corresponds to the isomorphism between $\Hom_{\DCoh(\X_0)}(M,L^q\otimes N \otimes L^{-p})$ and $\Hom_{\DCoh(\X_0)}(M,L^{q+1}\otimes N \otimes L^{-(p+1)})$. The $\Hom$-complex is in turn isomorphic to
$$\cong  \hocolim_p \Hom_{\QC(\X_0)}(M,\phi_\ast(\cO_{\X_0'}) \otimes N \otimes L^{-p})[2p],$$
where this last isomorphism uses the identification $\phi_\ast(\cO_{\X_0'}) \simeq \bigoplus_{n\in \bZ} L^n$. It also uses the fact that $M$ is coherent and $\phi_\ast(\cO_{\X_0'}) \otimes N \otimes L^{-p}$ is homologically bounded above, so that we may commute $\Hom_{\QC(\X_0)}(M,-)$ with this infinite direct sum. Using the projection formula and the adjunction between $\phi_\ast$ and $\phi^\ast$, we finally obtain a natural equivalence $$\cong \hocolim_p \Hom_{\DCoh(\X'_0)}(\phi^*(M),\phi^*(N)\otimes \phi^*(L)^{-p})[2p] $$ 
the multiplication by $\beta$ on $\Hom_{\DSing(\X,W)}(M,N) \otimes k\ls{\beta}$ now corresponds to the canonical isomorphisms: 
$$\Hom_{\DCoh(\X'_0)}(\phi^*(M),\phi^*(N)\otimes \phi^*(L)^{-p}) \to \Hom_{\DCoh(\X'_0)}(\phi^*(M),\phi^*(N)\otimes \phi^*(L)^{-(p+1)})$$
which arise from the canonical trivialization  $\mathcal{O}_{X'_0} \to \phi^*(L^{-1})$.
This operator is identified with the operator $\beta$ in \autoref{defn:MFcat} and therefore we have naturally identified the $\Hom$-complex in $\DSing(\X,W) \itimes_k k\ls{\beta}$ with $ \Hom_{\MF(\X',W')}(\phi^*(M),\phi^*(N))$ as required.




\end{proof}


We can now establish our main result on the degeneration property for graded singularity categories:

\begin{prop} \label{prop:graded_lg_degeneration}
Let $W : X/G \to \bA^1 / \Gm$ be a graded LG-model with associated LG-model $W' : X'/G \to \bA^1$. If $X'$ admits a semi-complete KN stratification and $\Perf(\Crit(W')/G)$ is a proper dg-category, then $\DSing(X/G,W)$ satisfies the $k$-linear degeneration property.
\end{prop}

\begin{lem} \label{lem:graded_degeneration}
Let $\cC$ be a $k$-linear (i.e., $\bZ$-graded) dg-category. Then the degeneration property for $\cC$ is equivalent to the $k\ls{\beta}$-linear degeneration property for $\cC\itimes_k k\ls{\beta}$.
\end{lem}

\begin{proof}
When $\mathcal{D}= \mathcal{C}\itimes k\ls{\beta}$, then the Barr complex \eqref{eqn:hochschild} computing $\hoch^{k\ls{\beta}}(\mathcal{D})$ is quasi-isomorphic to the subcomplex in which all of the objects $o_i$ lie in the generating set of objects of $\cD$ of the form $E \otimes k\ls{\beta}$, with $E \in \cC$. It follows that
$$\hoch^{k\ls{\beta}}(\mathcal{D}) \cong \hoch (\mathcal{C})\otimes k\ls{\beta}$$
canonically as dg-$\Lambda\ls{\beta}$-modules. We therefore have that $$\hoch^{k\ls{\beta},(n)} (\mathcal{D}) \cong \hoch^{(n)} (\mathcal{C})\otimes k\ls{\beta} $$ on the level of chain complexes as well. The result follows since as a complex $k\ls{\beta} \cong \bigoplus_{n \in \bZ} k[2n]$, so the homology $H_*( \hoch^{(n)} (\mathcal{C})\otimes k\ls{\beta}) \cong  H_*(\hoch^{(n)} (\mathcal{C}))\otimes k\ls{\beta}$ will be flat over $k[u]/u^n$ if and only if $H_*(\hoch^{(n)} (\mathcal{C}))$ is flat over the same ring. 
 \end{proof}
 
\begin{proof}[Proof of \autoref{prop:graded_lg_degeneration}]
By \autoref{lem:graded_degeneration} and \autoref{prop:forget_Z_gradings}, the degeneration property for $\DSing(X/G,W)$ is equivalent to the $k\ls{\beta}$-linear degeneration property for $\MF(X'/G,W')$, so the result follows from \autoref{prop:MFdeg}.
\end{proof}

\begin{rem}
Using \autoref{ex:graded_lg_model}, \autoref{prop:graded_lg_degeneration} implies noncommutative Hodge-de Rham degeneration for the derived category of certain complete intersections $s^{-1}(0)/G$, where $s$ is a $G$-invariant section of a locally free sheaf $\cE$ on a smooth $G$-scheme $Y$. In fact, the hypotheses of \autoref{amplif:zero_locus} imply that in the associated LG-model $W : \op{Tot}(\cE^\dual) / G \to \bA^1$, $\Crit(W)$ admits a complete KN-stratification, and thus $\Perf(\Crit(W)/G)$ is a proper dg-category by \autoref{cor:KN_strat_properness}. \autoref{prop:graded_lg_degeneration} therefore provides an alternate proof of \autoref{cor:zero_locus_degeneration}.
\end{rem}

\section{Computations of Hochschild invariants}
\label{sect:explicit}

\subsection{Generalities on Hochschild invariants}


In this section we identify the Hochschild homology with functions on the derived inertia stack (or loop stack) $\X \times^L_{\X \times \X} \X$, which we denote by $I_\X$, and we give an explicit description when $\X$ is a quotient stack. 

\begin{prop}  \label{prop:Hoch_Perf}
Let $\X$ be a smooth algebraic $k$-stack which is perfect, i.e., $\X$ is quasi-compact with affine diagonal, and $\QC(\X) \cong \op{Ind}(\Perf(\X))$, and let $\Delta : \X \to \X \times \X$ be the diagonal. Then we have an identification
$$\hoch(\Perf(\X)) \cong \RGamma(I_\X, \mathcal{O}_{I_{\X}}).$$
\end{prop}
\begin{proof}
Morita theory for perfect stacks \cite{ben2010integral} identifies the identity functor with $\Delta_* \mathcal{O}_X $ in the category  $\operatorname{QCoh}(\X\times \X)$. To compute $\hoch$, we use the Morita invariant definition of the Hochschild homology of a compactly generated dg-category as the trace of the the identity functor. Thus we must compute the trace
$$\hat{\operatorname{tr}}: \operatorname{QCoh}(\X \times \X) \to \operatorname{QCoh}(\operatorname{Spec} k).$$
On sheaves of the form $\pi_1^*(P_1) \otimes \pi_2^*(P_2)$, with $P_1,P_2 \in \Perf(\X)$, we have that the trace is given by
\begin{align*}
\operatorname{tr}(\pi_1^*(P_1) \otimes \pi_2^*(P_2)) &:=    \operatorname{RHom}(\iHom(P_1,\mathcal{O}_\X),P_2)\\
&\cong  R\Gamma (\X,  \Delta^*(\pi_1^*(P_1) \otimes \pi_2^*(P_2) ) ).
\end{align*}
Since the category $\QC(\X \times \X)$ is the colimit completion of sheaves of this form, we have that for an arbitrary object $F \in \QC(\X\times \X)$, the trace can be computed by
\begin{equation}\label{eqn:trace_formula}
F \to \RGamma(\X, \Delta^*(F)).
\end{equation}
It follows that we have an isomorphism $$\hoch(\Perf(\X)) \cong  \RGamma(\X, \Delta^*\Delta_*\mathcal{O}_X).$$ 
\end{proof}

Now let $X/G$ be a global quotient stack. We consider the scheme $\mathcal{P}:=G \times X \times X$. Denote by $\bar{\Delta}: G \times X\to \mathcal{P}$ the map $(g,x) \mapsto (g,x,x)$, and by $\Gamma: G\times X \to \mathcal{P}$ the map $(g,x) \mapsto (g,x,g \cdot x)$.  Both are closed immersions, and we will also use the notation $\bar{\Delta}$ and $\Gamma$ to denote the corresponding subschemes of $\cP$.

\begin{lem} \label{lem:derivedint}
Let $X/G$ be a smooth quotient stack. Both $\Gamma$ and $\bar{\Delta}$ are equivariant with respect to the $G$ action on $\cP$ which sends $h \cdot (g,x_1,x_2) \to (hgh^{-1},hx_1,hx_2)$. We have $$\hoch(\Perf(X/G)) \cong R\Gamma(X,\mathcal{O}_{\bar{\Delta}} \otimes^L \mathcal{O}_{\Gamma})^G.$$
 \end{lem}

\begin{proof} By \autoref{prop:Hoch_Perf} we must compute the derived global sections of the structure sheaf of the derived inertia stack. First note the alternate presentation for the stack $X/G \simeq G \times X/ G^2$, where the $G^2$ action in the second presentation is given by
$$(h_1,h_2) \cdot (g,x)=(h_2 g h_1^{-1}, h_1 x).$$
In this presentation the diagonal $X/G \to X/G \times X/G$ corresponds to the $G^2$-equivariant map $G \times X \to X \times X$ given by $(g,x) \mapsto (x,gx)$.

Let $G^2$ act on $\cP$ by
$$(h_1,h_2) \cdot (g,x_1,x_2) = (h_2 g h_1^{-1},h_1 x_1, h_2 x_2).$$
Then $\Gamma$ is $G^2$-equivariant, and using the presentation above we see that the diagonal factors as the closed immersion $\Gamma : G \times X/ G^2 \to \cP / G^2$ followed by the projection $\cP / G^2 \to X \times X / G^2$, which is smooth and affine. It follows that the derived inertia stack is the derived intersection of $p_1^{-1} \Gamma$ and $p_2^{-1} \Gamma$ in $\cP \times_{X^2} \cP / G^2$.

Now $\cP \times_{X^2} \cP \simeq G \times G \times X  \times X$ with $G^2$-action given by 
$$(h_1,h_2) \cdot (g_1,g_2,x_1,x_2) = (h_2 g_1 h_1^{-1}, h_2 g_2 h_1^{-1} ,h_1 x_1, h_2 x_2)$$
The projections $p_1,p_2 : \cP \times_{X^2} \cP \to \cP$ are given by forgetting $g_2$ and $g_1$ respectively. We claim that $\cP \times_{X^2} \cP / G^2 \simeq \cP / G$, where $G$ acts on $\cP$ as in the statement of the lemma. Indeed we can present $\cP/G$ as the quotient of $G \times \cP$ by the $G^2$-action
$$(h_1,h_2) \cdot (g_1,g_2,x_1,x_2) = (h_2 g_1 h_1^{-1}, h_1 g_2 h_1^{-1},h_1 x_1,h_1 x_2),$$
and we have a $G^2$-equivariant isomorphism $G \times \cP \to \cP \times_{X^2} \cP$ given by
$$(g_1,g_2,x_1,x_2) \mapsto (g_1, g_1g_2,x_1,g_1x_2)$$
The resulting isomorphism $\cP / G \to \cP\times_{X^2} \cP / G^2$ is given by the map $(g,x_1,x_2) \mapsto (1,g,x_1,x_2)$ which is equivariant with respect to the diagonal homomorphism $G \to G^2$.

In order to complete the proof, we must identify the closed substacks $p_1^{-1} (\Gamma / G^2)$ and $p_2^{-1}(\Gamma/G^2)$ in $\cP \times_{X^2} \cP / G^2$ under the isomorphism with $\cP/G$. The first is the closed subscheme $p_1^{-1} \Gamma \cap (\{1\}\times \cP) = \bar{\Delta}$, regarded as a $G$-equivariant closed subscheme of $\cP$, and the second is $p_2^{-1}(\Gamma) \cap (\{1\} \times \cP) = \Gamma$ as a $G$-equivariant closed subscheme of $\cP$.
\end{proof}  

The case when $X$ is a vector space, $\mathbb{V}$, and $G$ acts on $\mathbb{V}$ via a linear action, is of interest in two-dimensional gauge theory. In this case we make the above derived intersection explicit using a Koszul resolution. Denote by $\alpha: G \times \mathbb{V} \to \mathbb{V}$ the action morphism $(g,v) \mapsto g \cdot v$. We choose linear coordinates on $\mathbb{V}$ and identify $\mathbb{V} \times \mathbb{V}$ with $\operatorname{Spec}(k[x_i,y_i])$.

The Koszul complex for the regular sequence $\operatorname{K}_{\mathbb{V} \times \mathbb{V}} (x_i-y_i)$ gives a resolution of the diagonal on $\mathbb{V} \times \mathbb{V}$. An important point is that, in this case, this resolution is $G$-equivariant with respect to the diagonal $G$-action because the $G$-action on $\mathbb{V}$ is linear. Then $$\operatorname{K}_{G \times \mathbb{V} \times \mathbb{V}}(x_i-y_i) \to \mathcal{O}_{\bar{\Delta}}$$ is a resolution of  $\mathcal{O}_{\bar{\Delta}}$ over $\cP$.

\begin{cor} $\hoch(\Perf(\mathbb{V}/G)) \cong (K_ {G \times \mathbb{V}} (x_i-\alpha^*(x_i)))^G $ \end{cor}
\begin{proof}
By the above lemma, $\hoch(\Perf(\mathbb{V}/G))$ is isomorphic to
$$( \operatorname{K}_{G\times \mathbb{V} \times \mathbb{V}}(x_i-y_i) \otimes_{\mathcal{O}_{\cP}} \mathcal{O}_{\Gamma} )^G \simeq (\Gamma^\ast K_{G \times \mathbb{V} \times \mathbb{V}}(x_i-y_i))^G.$$
Pulled back to $G \times \mathbb{V}$, the function $\Gamma^*(y_i)=\alpha^*(x_i)$, hence $\Gamma^\ast K_{G\times \mathbb{V} \times \mathbb{V}}(x_i-y_i) = K_{G\times \mathbb{V}}(x_i-\alpha^\ast(x_i))$.
\end{proof}

\begin{rem} We expect that one could describe similar models for the Hochschild homology of general gauged linear sigma models  $(\mathbb{V}/G,W)$. The key step needed to do this would be to generalize \cite{Preygel}*{Theorem 4.2.3} to quotient (or more generally QCA) stacks. 

\end{rem} 



\subsection{An HKR theorem and quotients of affine varieties}
In \cite{block1994equivariant}, Block and Getzler construct for any compact smooth $M$-manifold $X$ an explicit model for the $M$-equivariant cyclic homology of the algebra $C^\infty(X)$ using differential forms on $X$. Our goal is to translate their construction into algebraic geometry and establish their version of the equivariant Hochschild-Kostant-Rosenberg theorem when $X=\operatorname{Spec}(A)$ is smooth and affine, $G$ is a reductive group, and $X/G$ is formally proper, which is equivalent to the condition that $A^G$ is finite dimensional over $k$. Our proof is an application of \autoref{thm:motivic_main} and \autoref{thm:lattice_conjecture}. For simplicity, we let $k=\mathbb{C}$ throughout this section.

To compute the derived intersection appearing in \autoref{lem:derivedint}, we may use the bar resolution $B(A)$ of $A$ as an $A-A$ bimodule. Namely, $$ B_n(A):=A\otimes A^{\otimes n} \otimes A $$ where the differential can be described as the sum $b=\Sigma_i(-1)^i \partial_i$, where
\[\partial_i(a'_0\otimes a_1 \otimes \cdots a_n \otimes a''_0):= \left\{\begin{array}{ll} a'_0 a_1 \otimes \cdots a_n \otimes a''_0,  & i=0 \\ a'_0 \otimes \cdots \otimes a_ia_{i+1} \otimes \cdots \otimes a''_0,  & i \neq 0, \neq n \\ a'_0 \otimes \cdots \otimes  a_{n}a''_0, & i=n \end{array}\right.\]
Our notation is meant to highlight the fact that the first and last variables in the bar complex play a distinguished role from the other $a_i$. We then have that $\cO_G \otimes B(A)$ is a resolution of $\mathcal{O}_{\bar{\Delta}}$ which we may restrict to $\Gamma$. The result is a complex where the $n$-th graded piece is
$$C_{n}(A,G):=\cO_G \otimes A^{\otimes {n+1}}=\Gamma(G\times X^{n+1}, \cO_{G \times X^{n+1}}).$$
For any $c \in \Gamma(\cO_{G \times X^{n+1}})$, the differentials $\partial_i$ above now take the form
\[
\partial_ic(g,x_0,x_1,\cdots,x_{n-1}):=
\left\{ \begin{array}{ll}
c(g,x_0,x_0,x_1,\cdots,x_{n-1})  &  i=0 \\
c(g,x_0,\cdots, x_i,x_i,\cdots,x_{n-1}) ,  & i \neq 0, \neq n \\
c(g,x_0,\cdots, x_{n-1},g \cdot x_0), & i=n
\end{array}\right.
\]

We define the $\Lambda$-module $C_{\bullet,G}(A) := \hoch(A,G)^G$, and note the following corollary of \autoref{lem:derivedint}.
\begin{cor} $\hoch(\Perf(X/G)) \cong C_{\bullet,G}(A)$ \end{cor}

Let $C_{\bullet,G}(A)_{g}^\wedge$ denote the completion of $C_{\bullet,G}(A)$ as a complex of modules over the representation ring $\Rep(G)\otimes \bC=\Gamma(\cO_G)^G$ at the conjugacy class $[g]$. Let us work for the moment with a fixed normal element $g \in G$. Let $Y=\operatorname{Spec}(A)^g$ denote the fixed point locus of $g$ and $B=\Gamma(\mathcal{O}_Y)$. The letter $Z$ will designate the centralizer of $g$ and $\mathfrak{z}$ denotes its Lie algebra, and normality of $g$ ensures that $Z$ is the complexification of $Z_c := Z \cap M$ for a maximal compact subgroup $M \subset G$. We have embeddings $j: Z \to G$ and $k: Y \to X$.

\begin{lem}
When $\operatorname{Spec}(A)/G$ is formally proper, the natural restriction map gives rise to an isomorphism $k^*:C_{\bullet,G}(A)_{g}^\wedge \to C_{\bullet,Z}(B)_{g}^\wedge$.
\end{lem} 
\begin{proof}
Note that because $\op{Spec}(A)/G$ is formally proper, each $C_{n,G}(A)$ is a coherent $\Rep(G) \otimes \bC$-module, so completion commutes with taking homology in this case, and it suffices to prove the result on the level of homology. It is known \cite{freed2008twisted}*{Proposition 3.10} that the map $k^*: K^\ast_M(X^{an},\mathbb{C})_g^\wedge \to K_{Z_c}(Y^{an},\mathbb{C})_g^\wedge$ is an isomorphism. The comparison maps 
$$
\xymatrix{
K^\ast_M(X^{an},\mathbb{C}) \ar[r]^\simeq \ar[d]^{k^*} & \operatorname{H}_\ast (C_{\bullet,G}^{per}(A)) \ar[d]^{k^*}\\
K^\ast_{Z_c}(Y^{an},\mathbb{C}) \ar[r]^\simeq & \operatorname{H}_\ast (C_{\bullet,Z}^{per}(B))
}
$$
are maps of $\Rep(G) \otimes \bC$ modules so we conclude that the map $k^*: \operatorname{H}_\ast (C^{per}_{\bullet,G}(A))_g^\wedge \to  \operatorname{H}_\ast (C^{per}_{\bullet,Z}(B))_g^\wedge $ is an isomorphism as well. By \autoref{thm:Hodge_main}, the vector spaces $\operatorname{H}_\ast(C^{per}_{\bullet,G}(A))$ and $\operatorname{H}_\ast(C^{per}_{\bullet,Z}(B))$ admit compatible Hodge structures. We observe that the Hodge decompositions
$$ \operatorname{H}_\ast(C^{per}_{\bullet,G}(A)) \cong \bigoplus_{n \equiv \ast \mod 2} \op{H}_n(C_{\bullet,G}(A)) $$ $$\operatorname{H}_\ast(C^{per}_{\bullet,Z}(B)) \cong \bigoplus_{n \equiv \ast \mod 2} \op{H}_\ast(C_{\bullet,Z}(B))$$ 
are decompositions of $\Rep(G)$ and $\Rep(Z)$ modules respectively. This follows because both the Hodge filtration and the conjugate filtration are filtrations of $\Rep(G)$ modules as can be seen for example by examining the explicit model for $C^{per}_{\bullet,G}(A)$. The lemma now follows by taking completions of these decompositions. 
\end{proof}

Next we construct a model for $C_{\bullet,Z}(B)_g^\wedge$ based on algebraic differential forms $\Omega^n_Y$, regarded as a projective $B$-module. Recall that the Cartan differential
$$\mathfrak{i} : \Sym(\fz^\ast) \otimes \Omega^n_Y \to \Sym(\fz^\ast) \otimes \Omega_Y^{n-1}$$
is the unique extension of the contraction map $\Omega_Y^n \to \fz^\ast \otimes \Omega_Y^{n-1}$ to a differential satisfying the Liebniz rule. Alternatively, regarding $\omega \in \Sym(\fz^\ast) \otimes \Omega_Y^n$ as a section of a quasi-coherent sheaf over $\fz$, we have
$$(\mathfrak{i} \omega)(z)=\mathfrak{i}_z \omega(z).$$
We thus have a chain complex, in fact a CDGA,
$$\Omega^\bullet_Y[\fz^\ast] = \left( \bigoplus_n \Sym(\fz^\ast) \otimes \Omega_Y^n[n] , \mathfrak{i} \right).$$
Note that $\mathfrak{i}$ is $Z$-equivariant, and that it descends to the quotient $\Sym(\fz)^\ast / \fm^k$. Thus we can define
\begin{gather*}
\Omega^\bullet_Y\ps{\mathfrak{z}^{\ast}}_k := \left(\bigoplus_n \operatorname{Sym}(\mathfrak{z}^{\ast})/\mathfrak{m}^k \otimes \Omega_Y^n[n],\mathfrak{i} \right) \\
\Omega^\bullet_Y\ps{\mathfrak{z}^{\ast}}^{Z}_k :=\left ( \bigoplus_n (\operatorname{Sym}(\mathfrak{z}^{\ast})/\mathfrak{m}^k \otimes \Omega_Y^n[n])^Z,\mathfrak{i} \right)
\end{gather*}

\begin{prop} \label{prop:HKR}
The comparison map of \autoref{construct:HKR} below is a quasi-isomorphism of $\Lambda$-modules
$$ \HKR_g^\wedge: C_{\bullet,Z}(B)_{g}^\wedge \to \varprojlim_k (\Omega^\bullet_Y\ps{\mathfrak{z}^{\ast}}^Z_k).$$ 
Hence when $\op{Spec}(A)/G$ is formally proper, we have a quasi-isomorphism of $\Lambda$-modules
$$\HKR_g^\wedge \circ k^*: C_{\bullet,G}(A)^\wedge_{g} \to \varprojlim_k (\Omega^\bullet_Y\ps{\mathfrak{z}^*}^Z_k).$$
\end{prop}

Let $Z_{(k)}$ denote the $k$-th infinitesimal neighborhood of the identity in $Z$.  The exponential map provides a compatible system of isomorphisms $\Exp_k: \op{Spec}(\Sym(\fz^\ast)/\fm^k)) \to Z_{(k)}$. Note that under this equivalence, $\Gm$ acts algebraically on $Z_{(k)}$ by scaling, and this action actually extends to an action of the monoid $\mathbb{A}^1$. This is encoded algebraically via a coaction map $\Sym(\fz^\ast)/\fm^k \to \Sym(\fz^\ast)/\fm^k \otimes \bC[t]$.
\begin{construction} \label{construct:HKR}
For any $b \in B$, the coaction of $\cO_Z$ on $B$, the exponential map $\Exp_k$, and the $\Gm$-action on $\fz$ define an element
$$\Exp_k(-t \cdot z) \cdot b \in B \otimes \Sym(\fz^\ast)/\fm^k \otimes \bC[t].$$
We define $C_{n,k}(B,Z) := \cO_{Z_{(k)}} \otimes B^{n+1}$, i.e. the reduction of $\hoch(B,Z)$ modulo $\fm^k$, and introduce the map $\HKR_{g,k}: C_{n,k}(B,Z) \to \Omega^n_Y\ps{\mathfrak{z}^{\ast}}_k$ given by
\begin{equation} \label{eqn:HKR}
\psi \otimes b'_0 \otimes b_1 \otimes \cdots b_n \mapsto \psi(g \cdot \Exp_k(z)) \int_{\Delta_n} b'_0 d(\exp_k(-t_1 z) \cdot b_1)\wedge \cdots \wedge d(\Exp_k(-t_n z) \cdot b_n) dt_1dt_2 \cdots dt_n.
\end{equation}
Here $d(-)$ denotes the $\Sym(\fz^\ast)/\fm^k \otimes \bC[t]$-linear extension of the exterior derivative
$$d : B \otimes \Sym(\fz^\ast)/\fm^k \otimes \bC[t] \to \Omega_Y^1 \otimes \Sym(\fz^\ast)/\fm^k \otimes \bC[t].$$
The integrand is regarded as an element of $\Sym(\fz^{\ast})/\fm^k \otimes \Omega_Y^n \otimes \bC[t_1,\ldots,t_n]$, and the integral over the standard $n$-simplex $\Delta_n$ is regarded formally as a linear map $\bC[t_1,\ldots,t_n] \to \bC$. This formula is identical to the one used in \cite{block1994equivariant}, so it follows formally from the computations there that $\HKR_{g,k}$ is a chain map (See for instance \cite{block1994equivariant}*{Theorem 3.2}). This map is $Z$-equivariant, so it restricts to a chain maps
\begin{gather*}
\HKR_{g,k} : C_{\bullet,k}(B,Z)^Z \to \Omega^\bullet_Y\ps{\fz^\ast}_k^Z, \text{ and} \\
\HKR_g^\wedge := \varprojlim_k \HKR_{g,k} : C_{\bullet,Z}(B)_g^\wedge \to \varprojlim_k \Omega^\bullet_Y\ps{\fz^\ast}_k^Z.
\end{gather*}
\end{construction}

\begin{proof}[Proof of \autoref{prop:HKR}]
By the compatibility of the HKR maps with translation by the central element $g$, it suffices to consider the case $g=\operatorname{id}$. The maps $\HKR_{\op{id},k} : C_{\bullet,k}(B,Z) \to \Omega^\bullet_Y\ps{\fz^\ast}_k$ are a compatible family of maps of bounded complexes with coherent homology over $B \otimes \Sym(\fz^\ast) / \fm^k$. $\HKR_{\op{id},1}$ is the classical HKR map

$$b'_0 \otimes \cdots \otimes b_n \mapsto \frac{1}{n!} b'_0 db_1 \cdots db_n,$$
which is an equivalence of $\Lambda$-modules. Hence by Nakayama's lemma each $\HKR_{\op{id},k}$ is a quasi-isomorphism, and the same is true after taking $Z$-invariants. Hence $\HKR_{\op{id}}^\wedge$ is a quasi-isomorphism. The final statement of the proposition combines this with the previous lemma.
 \end{proof} 


\bibliographystyle{plain}
\bibliography{hodge_references}

\end{document}